%% file: main.tex
\newif\ifHAL 
\HALtrue

\ifHAL 
\documentclass[11pt, a4paper, twoside, oneside]{article}
\else
\documentclass[a4paper,fleqn]{cas-sc}
\fi

\input{preamble.tex}
\input{macros.tex} 
\ifHAL 
\graphicspath{{/}} 
\else
\graphicspath{{images/}} 
\fi
 
\begin{document}


\ifHAL
\title{Elasto-acoustic wave propagation in geophysical media using hybrid high-order methods on general meshes}
\makeatletter
\newcommand{\sep}{, }
\newcommand{\@keywords}{} 
\newenvironment{keywords}{%
    \renewcommand{\@keywords} 
}{} 
\renewcommand{\maketitle}{%
\begin{center}
\huge \@title \\[0.5cm]
\large
Romain Mottier\footnote[2]{CEA, DAM, DIF, F-91297 Arpajon, France and CERMICS, ENPC, Institut Polytechnique de Paris, F-77455 Marne-la-Vall\'ee cedex 2, and Centre Inria de Paris, 48 rue Barrault, CS 61534, F-75647 Paris, France. Email adress: romain.mottier@enpc.fr}, 
$\qquad$ Alexandre Ern \footnote[3]{CERMICS, ENPC, Institut Polytechnique de Paris, F-77455 Marne-la-Vall\'ee cedex 2, and Centre Inria de Paris, 48 rue Barrault, CS 61534, F-75647 Paris, France. E-mail address: alexandre.ern@enpc.fr},
$\qquad$ Laurent Guillot \footnote[4]{CEA,DAM, DIF, F-91297 Arpajon, France. Email adress: laurent.guillot@cea.fr},\\[0.5cm]
\end{center}
}
\makeatother
\fi
\maketitle 
\begin{abstract}
Hybrid high-order (HHO) methods are discretization schemes characterized by several interesting properties such as local conservation, geometric flexibility and high-order accuracy. Here, HHO methods are studied for the space semi-discretization of coupled elasto-acoustic waves in the time domain using a first-order formulation. Explicit and singly diagonal implicit Runge--Kutta (ERK $\&$ SDIRK) schemes are used for the time discretization. We show that an efficient implementation of explicit (resp. implicit) time schemes calls for a static condensation of the face (resp. cell) unknowns. Crucially, both static condensation procedures only involve block-diagonal matrices. Then, we provide numerical estimates for the CFL stability limit of ERK schemes and present a comparative study on the efficiency of explicit versus implicit schemes. Our findings indicate that implicit time schemes remain competitive in many situations. Finally, simulations in a 2D realistic geophysical configuration are performed, illustrating the geometrical flexibility of the HHO method: both hybrid (triangular and \rev{quadrilateral}) and nonconforming (with hanging nodes) meshes are easily handled, delivering results of comparable accuracy to a reference spectral element software based on tensorized elements. \ifHAL \\[0.125cm] \fi 
\ifHAL
{
\small
\noindent \textbf{Mathematics Subjects Classification.} 65M12, 65M60, 74J10, 74S05, 35L05.\\[0.125cm]
\noindent \textbf{Keywords.} Hybrid high-order methods (HHO), Elasto-acoustic coupling, Wave equations, Runge--Kutta schemes, Efficient implementation, General meshes.\\[0.125cm]
}
\fi
\end{abstract}
\ifHAL
\vspace{-1cm}
\fi

\section{Introduction} 

The accurate computation of mechanical wave propagation is crucial in numerous scientific fields, such as medical imaging \cite{SZABO_2014}, non-destructive testing \cite{LD_2009}, or geophysics \cite{LW_1995}. In geophysics for instance, seismic (\textit{i.e.}, elastic), hydroacoustic or infrasonic waves are used to probe the Earth at the local (petroleum and marine exploration) and global scales. Owing to a broad suite of recording instruments deployed in various Earth's environments, these wavefields jointly provide insights into the Earth's internal structure \cite{TOMO_2024}, seismic source properties, or help monitoring natural hazards \cite{AJDH_2010}. The success of these endeavors strongly depends on the ability to compute and understand how an elasto-acoustic wavefield propagates through media that are heterogeneous and different in nature (solid, fluid). Closed-form or semi-analytical solutions are available only for the simplest geometrical configurations \cite{AR_2002}, while asymptotic methods often rely on assumptions that may not be fulfilled in more complex scenarios (\textit{e.g.}, ray theory or single-scattering approximations). Consequently, numerical methods that accurately simulate the partial differential equations describing elasto-acoustic waves are required to simulate wave propagation in realistic media of complex geometry.

Advances in high-performance computing have triggered the emergence of various numerical methods. The finite difference method (FDM) is still widely used \cite{MKG_2014}, evolving from regular-grid to staggered-grid approaches with improved accuracy \cite{VIRIEUX_1986}. Despite its efficiency on rectilinear meshes, FDM has long struggled with complex geometries and strong material contrasts. Extensions to more intricate domains have been developed \cite{PS_2017}, but requiring an increased computational effort. Additionally, higher-order FDM introduce larger stencils, negatively impacting parallel scalability.

A natural alternative to FDM are continuous finite element methods (FEM), which accommodate unstructured meshes and intricate boundaries. However, low-order FEM is computationally inefficient \cite{MARFURT_1984}, and its high-order extension can suffer from aliasing issues. The spectral element method (SEM), introduced in \cite{PATERA_1984} for fluid dynamics and later adapted to time-domain wave equations \cite{KV_1998}, leverages the accuracy of spectral methods owing to high-order polynomial approximations based on Lagrangian basis functions whose nodes coincide with the tensorized Gauss-Lobatto-Legendre (GLL) quadrature points. This choice ensures a diagonal mass matrix and efficient explicit time-stepping while retaining high-order accuracy. SEM have been successfully applied to coupled elasto-acoustic wave propagation \cite{KBT_2000, CV_2004} with a description of the fluid velocity using a potential field. However, SEM essentially relies on a space discretization with quadrilateral/hexahedral elements, which poses meshing challenges for realistic 3D applications. Different extensions to simplicial meshes \cite{CJRT_2001, MVS_2006} have demonstrated good dispersion properties \cite{MR_2012}, yet they may introduce additional computational overheads, difficulties in increasing the order of accuracy, or a loss of mass matrix diagonalization \cite{MVS_2006}. Moreover, handling severe wavespeed discontinuities necessitates multiscale discretizations in space (nonconforming meshes with \textit{h-} or \textit{p-}adaptivity) or time (multiple-time stepping), leading to increased algorithmic complexity, higher computational costs \cite{CCV_2003}, and potential instabilities \cite{DJ_2005}.

These limitations have sparked interest in discontinuous methods, such as finite volume (FV) and discontinuous Galerkin (dG) methods. Both methods rely on local information exchanges through fluxes at mesh interfaces, thereby relaxing the classical FEM $C^0$-continuity constraint, and accommodating more naturally meshes with hanging nodes for instance. Some dG methods adopt a second-order formulation in time (\textit{e.g.}, interior penalty methods \cite{GSS_2006, DSW_2008, ABM_2020}), whereas others (including FV) adopt a first-order formulation in time, solving for both primal (pressure in acoustics, velocity in elastodynamics) and dual (velocity in acoustics, stress in elastodynamics) variables. Several dG variants exist, depending on the fluxes (some based on Godunov-type upwind fluxes \cite{DK_2006, WSBG_2010}, others on central fluxes \cite{MG_2015} or modified upwind \cite{YDPPW_2016}) and on the polynomial setting (\textit{e.g.,} dG-SEM combinations \cite{WSBG_2010, ABM_2020_SEM}). While dG methods achieve FEM-like convergence rates, their high computational costs due to the substantial amount of  degrees of freedom generally remain a drawback.

Thus, there is room for alternative numerical strategies. Among others, like the Distributional Finite-Difference Method (DFDM) \cite{MASSON_2022}, or the Virtual Element Method (VEM) \cite{DFMSV_2021}, one promising approach is the family of hybrid nonconforming methods, which leverage ideas from both FEM and dG. Hybrid nonconforming methods can achieve higher-order accuracy through a consistent gradient reconstruction based on an additional unknown field defined on the mesh skeleton, hence the term "hybrid". Compared to dG methods, hybrid nonconforming methods have several advantages: (i) They include a built-in and robust stabilization that is less sensitive to parameter tuning than the one used in dG methods; (ii) In the context of elliptic problems, they achieve $H^1$-convergence of order $(k+1)$ and $L^2$-convergence of order $(k+2)$, \textit{i.e.} one order higher than dG methods. For the wave equation in first-order form, they can achieve $L^2$-convergence of order $(k+2)$ under a post-processing for the primal variable \cite{CQ_2014} (although it may be insufficient in elastodynamics \cite{TVG_2017}). Moreover, these methods allow for a natural enforcement of coupling conditions by means of the unknowns attached to the mesh skeleton.

The most salient \rev{example} of hybrid nonconforming method is  the Hybridizable dG method (HDG). This method was first proposed in \cite{NPC_2011} for the wave equation in first-order form, following the seminal work \cite{CGL_2009} on HDG methods for elliptic problems (see also \cite{CFHJSS_2018} for a second-order formulation in time). HDG methods allow for efficient implicit time-stepping by leveraging static condensation (\textit{i.e.} the local elimination of the cell unknowns). Efficient explicit time-stepping is also possible whenever the stabilization bilinear form leads to a block-diagonal face-face matrix. Explicit and implicit HDG schemes for the wave equation are compared in \cite{KSMW_2016}, where explicit schemes were found more efficient for various configurations. Different conclusions are drawn in the present work, allowing at least to say that the explicit and implicit approaches are both relevant. \rev{We also mention \cite{ACSVS:25} where the coupling of acoustic and elastic methods is studied using HDG methods and a Laplace transform in time.}
  
Another compelling instance in the class of hybrid nonconforming methods is the Hybrid High-Order (HHO) method. The HHO method approximates the solution within the mesh cells and its trace on the mesh skeleton, both by polynomial unknowns. The setting is of equal-order when the cell and face polynomials are of the same order $k\ge0$, whereas the setting is of mixed-order when the cell polynomials are one order higher than the face polynomials. The two devising principles of HHO methods are: (i) a gradient reconstruction in each mesh cell from the local cell and face unknowns; (ii) a stabilization that weakly enforces locally the matching of the trace of the cell unknowns with the corresponding face unknowns. Although the devising principles in HHO differ from those in HDG, it has been shown in \cite{CDE_2016} (see also \cite{CEP_2021}) that tight connections exist between both approaches, as well as with other hybrid nonconforming methods such as Weak Galerkin and nonconforming VEM \rev{(respectively devised in \cite{WY_2013} and \cite{ALM_2016})}. HHO methods were introduced in \cite{DE_2015} for linear elasticity and in \cite{DEL_2014} for linear diffusion. Since then, HHO methods have been applied to a broad range of problems. We refer the reader to the two textbooks \cite{CEP_2021, DD_2020} and to \cite{BDE_2022, BDE_2022bis, SEJD_2023} for the development of HHO methods applied to the (purely) acoustic or elastic wave equation in first- and second-order form (see also \cite{ES_2024, EK_2024} for further analysis results).

The approximation of coupled elasto-acoustic waves using HHO methods in space and Runge--Kutta (RK) schemes in time has been started in \cite{MEKG_2025}. Therein, the devising principles of the method were established, and an error analysis in the space semi-discrete setting was proposed. Moreover, numerical results, mostly on academic test cases, were presented to showcase the stability and accuracy properties of the proposed method. Let us briefly summarize the main findings. On the one hand, if explicit RK (ERK) time-stepping schemes are considered, it is recommended to scale the stabilization by an $\cal{O}(1)$-weight (to alleviate the CFL restriction) and to use a plain Least-Squares stabilization. Because of this, it is natural to consider an equal-order setting for the cell and face unknowns. On the other hand, for implicit RK (IRK) schemes, one can consider tighter stabilizations to obtain higher-order convergence rates. Namely, one can scale the stabilization by an $\cal{O}(\frac{1}{h})$-weight (where $h$ is the mesh size) and employ any of the high-order HHO stabilization operators. As discussed in \cite{CDE_2018} for elliptic problems, the mixed-order setting with Lehrenfeld--Sch\"oberl stabilization appears to be more computationally effective. 

In this work, we further develop the above methodology, addressing important computational aspects and assessing the approach on more realistic geophysical settings. The HHO formulations for the explicit and implicit time schemes are those discussed above. Our first main contribution in the paper is to recast the fully discrete problem in algebraic form, so as to show that an efficient implementation of implicit (resp. explicit) time schemes calls for a static condensation of the cell (resp. face) unknowns. Crucially, both static condensation procedures only involve block-diagonal matrices. Our second main contribution is to estimate numerically the CFL stability limit for explicit time schemes and to provide a comparative study on the efficiency of explicit versus implicit time schemes. This is performed on test cases with manufactured solutions and in bilayered problems with an initial Ricker wavelet. As mentioned above, our findings indicate that implicit time schemes remain competitive in many situations. Our third main contribution is to consider a realistic test case including waves trapped in a sedimentary bassin, and to show that the proposed methodology delivers accurate solutions, while offering the capability of handling general meshes with hanging nodes. Reference solutions for this latter test case are obtained using a specral element software using tensorized elements. 

The paper is organized as follows. Section 2 introduces the coupled elasto-acoustic model and presents its first-order weak formulation. Section 3 presents the HHO space semi-discretization. Section 4 describes the fully discrete problem, incorporating both explicit and implicit RK schemes and showing the relevance of static condensation procedures in both cases. Section 5 (i) investigates numerically the CFL stability limit for ERK and compares the efficiency (error vs. CPU time) of both explicit and implicit time schemes ; (ii) verifies the accuracy of the proposed methodology on a Ricker wavelet test case (with realistic material properties) for which a semi-analytical solution is available, and performs a further comparative study of explicit vs. implicit time schemes; (iii) showcases the performances of the proposed method on a two-dimensional test case inspired by recent seismoacoustic observations, focusing on the emission of acoustic energy due to the interaction between seismic waves and a sedimentary basin \cite{HLVHCPMRB_2018, AASSE_2018}. Numerical results obtained on quadrilateral, simplicial, hybrid, and nonconforming meshes (with hanging nodes) are compared to those produced by a reference spectral element solver based on tensorized elements, thus highlighting the versatility of the HHO method to handle general meshes in complex geological media.

\section{Model problem}\label{sec::model_problem}

This section presents the model problem. We use boldface (resp. blackboard) fonts for vectors (resp. tensors), and for vector-valued (resp. tensor-valued) fields and spaces composed of such fields. For a bounded and uniformly positive weight $\kappa$, we define the $\kappa$-weighted $L^2$-inner product as $(u, v)_{L^2(\kappa; \Omega)} := \int_{\Omega} \kappa u v ~\mathrm{d}\Omega$ for all $u, v \in L^2(\Omega)$. A similar notation is used for inner products involving vector- and tensor-valued fields. 

\subsection{Strong formulation}

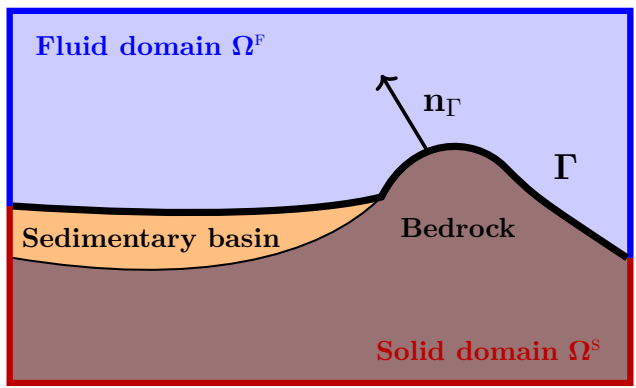
\begin{figure}[!htb]
\centering
\resizebox{0.5\textwidth}{!}{%
\begin{tikzpicture}[scale=0.85]
\draw[line width=1.75mm, ceared] (-1,0.85) -- (-1,-2) -- (9,-2) -- (9,0);
\draw[line width=1.75mm, blue] (9,0) -- (9,4) -- (-1,4) -- (-1,0.85);
\fill[orange!50] (-1,0.85) .. controls (2,1.5) and (4,1.5) .. (5,1.8)  
.. controls (5.5,2) and (6.5,2) .. (7,1.5)
-- (7,1) .. controls (6.5,1.2) and (5.5,1) .. (5,1)  
.. controls (4,0) and (2,-0.5) .. (-1,0) ;
\fill[reddishgray] (-1,0) .. controls (2,-0.5) and (4,0) .. (5,1)  
.. controls (5.5,2) and (6.5,2) .. (7,1.5) 
.. controls (7.5,1) and (9,0) .. (9,0) 
-- (9,-2) -- (-1,-2) -- cycle;
\fill[blue!20] (-1,0.85) -- (-1,4) -- (9,4) -- (9,0) .. controls (7.5,1) and (7,1.5) .. (7,1.5)
.. controls (6.5,2) and (5.5,2) .. (5,1) .. controls (4,0.75) and (2,0.65) .. (-1,0.85);
\draw[thick] (-1,0.85) .. controls (2,0.65) and (4,0.75) .. (5,1);  
\draw[thick] (-1,0) .. controls (2,-0.5) and (4,0) .. (5,1)
.. controls (5.5,2) and (6.5,2) .. (7,1.5)  
.. controls (7.5,1) .. (9,0);
\draw[line width=1mm] (-1,0.85) .. controls (2,0.65) and (4,0.75) .. (5,1) 
.. controls (5.5,2) and (6.5,2) .. (7,1.5)  
.. controls (7.5,1) .. (9,0);  
\draw[->, line width=0.5mm] (5.75,1.75) -- (5,3);
\node at (6.25,0.5) {\textbf{Bedrock}};
\node at (1.25,0.3) {\textbf{Sedimentary basin}};
\node at (1.25,3.5) {\color{blue} \textbf{Fluid domain} $\mathbf{\domain{f}}$};
\node at (6.75,-1.5){\color{ceared} \textbf{Solid domain} $\mathbf{\domain{s}}$};
\node at (8,1.5) {\Large $\mathbf{\G}$};
\node at (6,2.5) {\Large $\mathbf{n}_\G$};
\end{tikzpicture}}
\caption{Solid subdomain $\domain{s}$ composed of the sedimentary basin and the bedrock, fluid subdomain $\domain{f}$, and  unit normal $\bd{n}_\Gamma$ along the interface $\Gamma$ (highlighted with a thick black line)}.
\label{fig:model}
\end{figure}

Let $J := (0, T_{\rm{f}})$ denotes the time interval of interest, where $T_{\rm{f}} > 0$ represents the final time, and let $\Omega$ be a polyhedral domain (open, bounded, connected, Lipschitz set) in $\mathbb{R}^d$, with $d \in \{2, 3\}$. We partition $\Omega$ into two disjoint, polyhedral subdomains, $\domain{f}$ and $\domain{s}$ (of boundaries $\partial \domain{f}$ and $\partial \domain{s}$), representing the solid and fluid media, respectively. These subdomains share the interface $\G := \partial \domain{f}\cap \partial \domain{s}$, and the unit normal vector $\bd{n}_{\G}$ on $\G$ conventionally points from $\domain{s}$ to $\domain{f}$. A simplified representation of this setting is displayed in \hyperref[fig:model]{\Cref{fig:model}} illustrating the geophysical application considered in Section 5.

The linear acoustic wave propagation is defined by the evolution of its scalar pressure field $p \left[ \rm{Pa} \right]$ and velocity field $\bd{m} \left[ \frac{\rm{m}}{\rm{s}} \right]$, which are described by the following system of partial differential equations in $J \times \domain{f}$:
\begin{subequations}\label{continuous_acoustic}
\begin{align}
\rho^\sc{f} \partial_t \bd{m} - \nabla p & = \bd{0},\label{continuous_acoustic_eq_1}\\
\frac{1}{\kappa} \partial_t p - \nabla \cdot \bd{m} & = f^\sc{f},
\label{continuous_acoustic_eq_2}
\end{align}
\end{subequations}
where $\rho^\sc{f} \left[ \frac{\rm{kg}}{\rm{m}^3} \right]$ is the fluid density, $\kappa \left[ \rm{Pa} \right]$ is the fluid bulk modulus, and $f^\sc{f} \left[ \rm{Hz} \right]$ is the source term. The equation \eqref{continuous_acoustic_eq_2} is derived by considering a perfect fluid, where the stress tensor is given by $\bbm{s}^{\sc{f}} := p \bbm{Id}$, instead of the classical convention $\bbm{s}^{\sc{f}} := - p \bbm{Id}$, where $\bbm{Id}$ is the identity tensor. The material properties $\rho^\sc{f}$ and $\kappa$ are typically piecewise constant in $\domain{f}$, although it is not mandatory (they can smoothly evolve in space). The speed of the acoustic waves is given by $c_{\sc{p}}^\sc{f} := \sqrt{\kappa / \rho^\sc{f}} \left[ \frac{\rm{m}}{\rm{s}} \right]$. The initial conditions are specified as $p(0) = p_0$ and $\bd{m}(0) = \bd{m}_0$ with initial data $p_0$ and $\bd{m}_0$. 

Let $\nabla_{\rm{sym}} := \frac{1}{2}(\nabla + \nabla^\dagger)$ denote the symmetric gradient operator. The elastic wave equations, which describe the evolution of the linearized Cauchy stress tensor $\bbm{s} \left[\rm{Pa}\right]$ and the velocity field $\bd{v} \left[ \frac{\rm{m}}{\rm{s}} \right]$,  consist of the following partial differential equations in $J \times \domain{s}$:
\begin{subequations}\label{continuous_elastic_eq}
\begin{align}
\bbm{C}^{-1} \partial_t \bbm{s} - \nabla_{\rm{sym}} \bd{v} & = \bb{0},\label{continuous_elastic_1} \\
\rho^\sc{s} \partial_t \bd{v} - \nabla \cdot \bbm{s} & = \bd{f}^\sc{s},
\label{continuous_elastic_2}
\end{align}
\end{subequations}
where $\rho^\sc{s} \left[ \frac{\rm{kg}}{\rm{m}^3} \right]$ is the solid density, $\bbm{C} \left[ \rm{Pa} \right]$ is the fourth-order Hooke tensor, and $\bd{f}^\sc{s} \left[ \frac{\rm{N}}{\rm{m^3}} \right]$ is the source term. The material properties $\bbm{C}$ and $\rho^{\sc{s}}$ can be considered piecewise constant in $\domain{s}$. In the framework of isotropic elasticity, $\bbm{C}$ depends on the Lamé parameters $\lambda \left[ \rm{Pa} \right]$ and $\mu \left[ \rm{Pa} \right]$ as $\bbm{C}_{i j k l} := \lambda \delta_{i j} \delta_{k l} + \mu (\delta_{i k} \delta_{j l} + \delta_{i l} \delta_{j k})$, for all $i,j,k,l \in \{1,...,d\}$, where $\delta$ is the Kronecker symbol. This setting gives rise to two wave speeds $ \left[ \frac{\rm{m}}{\rm{s}} \right]$ corresponding to two types of body waves, far from material interfaces:
\begin{subequations}\label{cs}
\begin{alignat}{3}
c_{\sc{p}}^\sc{s} & := \sqrt{\left(\lambda + 2 \mu\right) / \rho^\sc{s}} & \qquad & \text{for compressional (P-) waves},\label{cs:a}\\
c_{\sc{s}}^\sc{s} & := \sqrt{\mu / \rho^\sc{s}} & \qquad & \text{for shear (S-) waves.}\label{cs:b}
\end{alignat}
\end{subequations}
We do not consider the incompressible limit where $\frac{\lambda}{\mu} \gg 1$, so that the two wave speeds in \eqref{cs} are expected to be of similar magnitude. The initial conditions are specified as $\bd{v}(0) = \bd{v}_0$ and $\bbm{s}(0) = \bbm{s}_0$, with initial data $\bd{v}_0$ and $\bbm{s}_0$. 

For simplicity, we impose homogeneous Dirichlet boundary conditions on $p$ on $\partial \domain{f} \setminus \G$ and on $\bd{v}$ on $\partial \domain{s} \setminus \G$, which lead to perfect wave reflection at the domain boundaries. The computational domain and time are chosen appropriately, so as to avoid any spurious reflection that could interfere with the main phenomena we want to simulate.
\rev{Neumann boundary conditions on the normal fluid velocity and solid stress can be enforced in a straightforward manner if the boundary of the domain is composed of planar faces. If the domain has a curved boundary, one possibility to retain high-order is to use an unfitted HHO method \cite{BCDE_2021}; another one is to enrich the local basis functions on curved faces by using suitable non-polynomial functions \cite{YEMM_2024}.}

The interface conditions on $J \times \G$ are expressed as follows:
\begin{subequations}\label{continuous_coupling}
\begin{align}
\label{coupling1}
\bd{v} \cdot \bd{n}_\G & = \bd{m} \cdot \bd{n}_\G, \\
\bbm{s} \cdot \bd{n}_\G & = p \, \bd{n}_\G,
\label{coupling2}
\end{align}
\end{subequations}
where the first equation enforces the continuity of the normal velocity component across the interface (kinematic condition on a perfect fluid-solid boundary where tangential slip is allowed, but there is no separation or interpenetration), and the second equation (dynamic condition) expresses the balance of forces.

\subsection{Weak formulation}

We define the following functional spaces:  
\begin{subequations}
\begin{alignat}{6}
\bd{{M}^{\sc{f}}} & := \bd{L}^2(\domain{f}), \qquad && P^{\sc{f}} := \left\{p \in H^1(\domain{f}) : p|_{\partial \domain{f} \backslash \G} = 0\right\}, \\  
\bbm{{S}}^{\sc{s}} & := \bbm{L}^2_{\rm{sym}}(\domain{s}), && \bd{V}^{\sc{s}} := \left\{ \bd{v} \in \bd{H}^1(\domain{s}) : \bd{v}|_{\partial \domain{s} \backslash \G} = \bd{0} \right\}.
\end{alignat}   
\end{subequations}
The coupled elasto-acoustic wave problem then consists of finding $(\bd{m}, p) : J \to \bd{{M}^{\sc{f}}} \times P^{\sc{f}}$ and $(\bbm{s}, \bd{v}) : J \to \bbm{{S}}^{\sc{s}} \times \bd{V}^{\sc{s}}$ such that, for all $t \in J$, the following weak equations hold:
\ifHAL
\begin{enumerate}[label=\roman*)]
\else
\begin{enumerate}[a)]
\fi
\item For the acoustic wave equations, for all test functions $(\bd{r}, q) \in \bd{{M}^{\sc{f}}} \times P^{\sc{f}}$, 
\begin{subequations} 
\label{weak_form_acoustic_eq}
\begin{align}
(\partial_t \bd{m}(t), \bd{r})_{\bd{L}^2(\rho^\sc{f};\domain{f})} - (\nabla p(t), \bd{r})_{\bd{L}^2(\domain{f})} &= 0, \label{aweak1} \\
(\partial_t p(t), q)_{L^2(\frac{1}{\kappa};\domain{f})} + (\bd{m}(t), \nabla q)_{\bd{L}^2(\domain{f})} + (\bd{v}(t) \cdot \bd{n}_{\G}, q)_{L^2(\G)} &= (f^\sc{f}(t), q)_{L^2(\domain{f})}. 
\label{aweak2}
\end{align}
\end{subequations} 
\item For the elastic wave equations, for all test functions $(\bbm{b}, \bd{w}) \in \bbm{{S}}^{\sc{s}} \times \bd{V}^{\sc{s}}$,
\begin{subequations} 
\label{weak_form_elastic_eq}
\begin{align}
(\partial_t \bbm{s}(t), \bbm{b})_{\bbm{L}^2(\bbm{C}^{-1};\domain{s})} - (\nabla_{\rm{sym}} \bd{v}(t), \bbm{b})_{\bbm{L}^2(\domain{s})} &= 0, \label{eweak1} \\
(\partial_t \bd{v}(t), \bd{w})_{\bd{L}^2(\rho^\sc{s};\domain{s})} 
+ (\bbm{s}(t), \nabla_{\rm{sym}} \bd{w})_{\bbm{L}^2(\domain{s})} - (p(t) \bd{n}_{\G}, \bd{w})_{\bd{L}^2(\G)} &= (\bd{f^\sc{s}}(t), \bd{w})_{\bd{L}^2(\domain{s})}. 
\label{eweak2}
\end{align}
\end{subequations}  
\end{enumerate}
The coupling conditions \eqref{coupling1} and \eqref{coupling2} are directly incorporated into \eqref{aweak2} and \eqref{eweak2}, respectively. Thus, they are weakly imposed.

\section{HHO space semi-discretization}\label{sec::HHO_discretization}

In this section, we introduce a simple description of the space semi-discretization of (\ref{weak_form_acoustic_eq})-(\ref{weak_form_elastic_eq}) in the HHO framework. Namely, we take a quick look at the way the method deals with general meshes, and how (\ref{weak_form_acoustic_eq})-(\ref{weak_form_elastic_eq}) are approximated with HHO operators. Our goal is to provide the set of ODEs that drive the time-evolution of the space semi-discrete system in algebraic form, using a generic choice for the local polynomial bases.

\subsection{Mesh}

Let $\T$ be a polygonal mesh of the domain $\Omega$, constructed so as to align with the partition of $\Omega$ into the subdomains $\domain{f}$ and $\domain{s}$. The mesh cells are supposed to have planar faces, constructed so as to align with the different boundaries of the subdomains $\domain{f}$ and $\domain{s}$. Thus, we generate a piecewise affine approximation of the interface $\G$ and of the subdomains $\domain{s}$ and $\domain{f}$. We assume that the mesh is fine enough to neglect the geometric error introduced by this approximation. Consistently with this assumption, we also assume that all material properties are piecewise constant on $\T$. For latter use, we define the two sub-meshes, $\T^{\sc{f}}$ and $\T^{\sc{s}}$, which cover exactly the subdomains $\domain{f}$ and $\domain{s}$, respectively. 

The set of mesh faces $\cal{F}$ is partitioned as $\F := \F^\circ \cup \F^\partial$, where $\F^\circ$ collects all the interior faces, including those lying on the interface $\G$, and $\F^\partial$ collects all the boundary faces on $\partial\Omega$. The set of interior faces is further split as $\F^\circ := \F^{\circ \sc{f}} \cup \F^{\circ \sc{s}} \cup \F^\G$, whereas the set of boundary faces is split as $\F^\partial := \F^{\partial \sc{f}} \cup \F^{\partial \sc{s}}$, with obvious notation. To simplify later expressions, we define $\F^\sF := \F^{\circ \sF} \cup \F^{\partial \sc{f}} \cup \F^\G$ and $\F^\sS := \F^{\circ \sS} \cup \F^{\partial \sc{s}} \cup \F^\G$ for the sets of faces lying in $\overline{\domain{f}}$ and $\overline{\domain{s}}$, respectively. 

We consider general meshes which can contain hanging nodes (see \hyperref[fig:mesh]{\Cref{fig:mesh}}). These meshes can be handled in a seamless fashion by any discretization method supporting polyhedral cells, as HHO, dG and VEM methods. For instance, in the example shown in \hyperref[fig:mesh]{\Cref{fig:mesh}}, elements sharing an edge with the solid-fluid interface and containing hanging nodes from the other mesh (filled elements) are no longer considered as triangles in the fluid region or as quadrilaterals in the solid region, but are treated as quadrilaterals and hexagons, respectively. This approach greatly simplifies the meshing of the domain, particularly at the fluid-solid interface, since it allows for independent meshing of the two subdomains. 
\ifHAL
\begin{figure}[!htb]
\centering
\resizebox{\textwidth}{!}{
\begin{tikzpicture}
\draw[fill=ceared!25] (0, 0) rectangle (5, -1);
\draw[fill=blue!25] (0.75, 0) -- (1.25, 0) -- (1, 0.25) -- cycle;
\draw[fill=blue!25] (1.75, 0) -- (2.25, 0) -- (2, 0.25) -- cycle;
\draw[fill=blue!25] (2.75, 0) -- (3.25, 0) -- (3, 0.25) -- cycle;
\draw[fill=blue!25] (3.75, 0) -- (4.25, 0) -- (4, 0.25) -- cycle;
\foreach \x in {0, 1, 2, 3, 4} {
\foreach \y in {0} {
\draw[ceared] (\x, \y) rectangle (\x+1, \y-1);
}
}
\foreach \x in {0, 1, 2, 3, 4} {
\foreach \y in {-1} {
\draw[ceared] (\x, \y) rectangle (\x+1, \y-1);
}
}
\foreach \x in {0, 0.5, 1, 1.5, 2, 2.5, 3, 3.5, 4, 4.5} {
\foreach \y in {0.25} {
\draw[blue] (\x, \y) -- (\x+0.5, \y) -- (\x+0.25, \y-0.25) -- cycle;
}
}
\draw[line width=0.35mm] (-0.75,0) -- (5.75,0);
\foreach \x in {0, 1, 2, 3, 4} {
\foreach \y in {0} {
\fill[ceared] (\x+1, \y-1) circle (0.05);
\fill[ceared] (\x, \y) circle (0.05);
}
}
\foreach \x in {0, 1, 2, 3, 4} {
\foreach \y in {-1} {
\fill[ceared] (\x+1, \y-1) circle (0.05);
\fill[ceared] (\x, \y) circle (0.05);
}
}
\fill[ceared] (5, 0) circle (0.05);
\fill[ceared] (0, -1) circle (0.05);
\fill[ceared] (0, -2) circle (0.05);
\foreach \x in {0, 0.5, 1, 1.5, 2, 2.5, 3, 3.5, 4, 4.5} {
\foreach \y in {0.25} {
\fill[blue] (\x, \y) circle (0.05);
\fill[blue] (\x+0.5, \y) circle (0.05);
\fill[blue] (\x+0.25, \y-0.25) circle (0.05);
}
}
\node at (6,0) {\small $\G$};
\node at (5.5,0.35) {\color{blue} \small $\Tf$};
\node at (5.5,-0.35) {\color{ceared} \small $\Ts$};
\end{tikzpicture}
}
\caption{Generic example of a mesh with hanging nodes at the interface $\Gamma$. Filled elements are reinterpreted: filled squares as hexagons and filled triangles as quadrilaterals. Unfilled elements are treated for what they are: squares and triangles.}
\label{fig:mesh}
\end{figure}
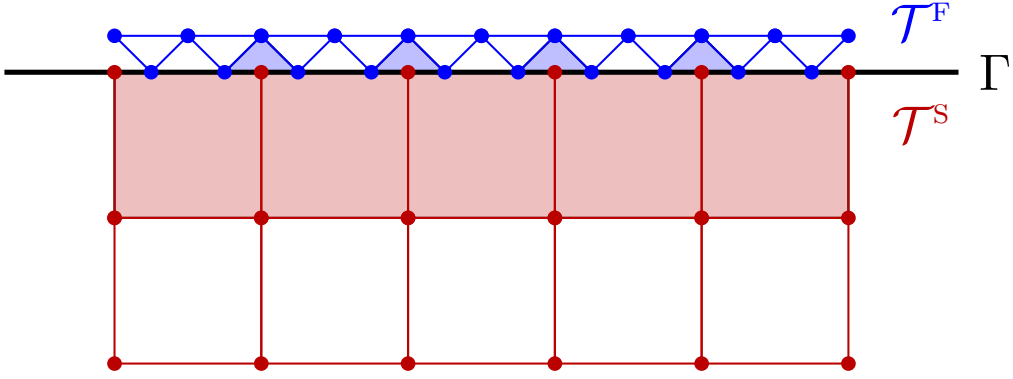
\fi

\subsection{Degrees of freedom}

Recall that the coupled problem \eqref{weak_form_acoustic_eq}-\eqref{weak_form_elastic_eq} relies on primal ($p$ and $\bd{v}$) and on dual ($\bd{m}$ and $\bbm{s}$) variables. Here, the primal variables are discretized using the HHO method, whereas the dual variables are discretized via a standard dG approach. The dG variables are thus piecewise polynomials of degree $k$ (only a cell component), whereas the HHO variables consist of a pair: a cell component and a face component. The cell component is a piecewise polynomial of degree $k^\prime \in \{k, k+1\}$, whereas the face component is a piecewise polynomial of degree $k$. The HHO discretization is referred to as equal-order when $k^\prime = k$ and as mixed-order when $k^\prime = k+1$. A polynomial degree $k \geq 1$ is required in the solid subdomain to control the rigid body motions \cite{DE_2015}. For simplicity, we consider the same polynomial degree $k$ in both subdomains, although it is not mandatory (we can take $k \geq 0$ in the fluid subdomain).

A generic mesh cell is denoted by $T \in \T$, with its diameter $h_T$, its outward unit normal $\bd{n}_T$, and the set of its boundary faces $\cal{F}_{\partial T}$. A generic mesh face is denoted by $F \in \F$. We define the local polynomial spaces $\bbm{P}_{\rm{sym}}^k(T)$, $\bd{P}^k(T)$ and  $P^k(T)$ (resp. $\bd{P}^k(F)$ and  $P^k(F)$) as the restrictions to the mesh cell $T$ (resp. the mesh face $F$) of symmetric tensor-, vector- and scalar-valued $d$-variate polynomials of degree at most $k$ (resp. $(d-1)$-variate polynomials of degree at most $k$). Altogether, we consider the following dofs \rev{(where the cross refers to a Cartesian product)}: 
\ifHAL
\begin{enumerate}[left=0pt, label=\textbullet]
\else
\begin{enumerate}[\textbullet]
\fi
\item from a local point of view, we have 
\begin{subequations}
\small 
\begin{alignat}{5}
& p_T \in P^{k^\prime}(T), & \qquad\! & p_{\partial T} \in P^{k}(\F_{\partial T}) := \underset{F \in \F_{\partial T}}{\bigtimes} P^{k}(F), & \qquad\! & \bd{m}_{T} \in \bd{P}^k(T), \qquad & \forall T \in \Tf,\\
& \bd{v}_{T} \in \bd{P}^{k^{\prime}}(T), & \qquad\! & \bd{v}_{\partial T} \in \bd{P}^{k}(\F_{\partial T}) := \underset{F \in \F_{\partial T}}{\bigtimes} \bd{P}^{k}(F), & \qquad\! & \bbm{s}_{T} \in \bbm{P}_{\rm{sym}}^k(T), \qquad & \forall T \in \Ts.
\end{alignat}
\end{subequations}
\item from a global point of view, we have\\[-0.5cm]
\begin{subequations}
\small 
\begin{alignat}{5}
& p_{\Tf} \in P^{k^\prime}(\Tf) := \underset{T \in \Tf}{\bigtimes} P^{k^\prime}(T), & \quad\!\! & p_{\Ff} \in P^{k}(\Ff) := \underset{F \in \Ff}{\bigtimes} P^{k}(F), & \quad\!\! & \bd{m}_{\Tf} \in \bd{M}^{k}(\Tf) := \underset{T \in \Tf}{\bigtimes} \bd{P}^k(T),
\\
& \bd{v}_{\Ts} \in \bd{{V}}^{k^\prime}(\Ts) := \underset{T \in \Ts}{\bigtimes} \bd{P}^{k^{\prime}}(T), & \quad\!\! & \bd{v}_{\Fs} \in \bd{{V}}^{k}(\Fs) := \underset{F \in \Fs}{\bigtimes} \bd{P}^k(F), & \quad\!\! & \bbm{s}_{\Ts} \in \bbm{S}^k(\Ts) := \underset{T \in \cal{T}^{\sc{s}}}{\bigtimes} \bbm{P}_{\rm{sym}}^k(T).
\end{alignat}
\end{subequations}
\end{enumerate}
In \hyperref[dofs]{\Cref{dofs}}, we illustrate the dofs involved in the space discretization of a 2D domain near a solid-fluid interface, for hexagonal meshes and polynomial degrees $k^\prime = k = 1$. There are $9$ cell unknowns in each fluid cell ($3$ for the primal variable and $6$ for the dual variable), and there are $2$ unknowns for each face $F \in \F^{\circ \sc{f}}$. There are $15$ cell unknowns in each solid cell ($6$ for the primal variable and $9$ for the dual variable), and there are $4$ unknowns for each face $F \in \F^{\circ \sc{s}}$. Finally, there are $6$ unknowns for each face $F \in \F^{\G}$ ($2$ from the fluid part and $4$ from the solid part).
\ifHAL
\begin{figure}[!htb]
\centering
\hspace{0.25cm}
\resizebox{0.9\textwidth}{!}{%
\begin{tikzpicture}[scale=2]
\node at (0.0,1.25) {\color{blue} \textbf{Fluid domain} $\mathbf{\domain{f}}$};
\node at (-1.75,0.5) {\color{blue} \textbf{Cell} $T \in \Tf$};
\node at (0.0,-2.75) {\color{ceared} \textbf{Solid domain} $\mathbf{\domain{s}}$};
\begin{scope}[xshift=-2cm, yshift=-0.25cm]
\draw[black, line width=2pt, fill=blue] (-1.025, -3.25) circle (2.5pt);
\filldraw[black, line width=1.5pt] (-0.75, -3.25) node[shape=diamond, draw, fill=blue, inner sep=0pt, minimum size=10pt] {};
\node at (-0,-3.25) {\textbf{Scalar dofs}};
\begin{scope}[xshift=2.5cm]
\draw[black, line width=2pt, fill=ceared] (-1, -3.4) circle (2.5pt);
\draw[black, line width=2pt, fill=ceared] (-1.2, -3.4) circle (2.5pt);
\draw[black, line width=2pt, fill=blue] (-1, -3.1) circle (2.5pt);
\draw[black, line width=2pt, fill=blue] (-1.2, -3.1) circle (2.5pt);
\filldraw[black, line width=1.5pt] (-0.71, -3.25) node[shape=diamond, draw, fill=ceared, inner sep=0pt, minimum size=10pt] {};
\filldraw[black, line width=1.5pt] (-0.615, -3.25) node[shape=diamond, draw, fill=ceared, inner sep=0pt, minimum size=10pt] {};
\node at (0.25,-3.25) {\textbf{Vector dofs}};
\end{scope}
\begin{scope}[xshift=5.25cm]
\draw[black, line width=2pt, fill=ceared] (-1.025, -3.37) circle (2.5pt);
\draw[black, line width=2pt, fill=ceared] (-1.025, -3.175) circle (2.5pt);
\draw[black, line width=2pt, fill=ceared] (-1.195, -3.175) circle (2.5pt);
\node at (-0.1,-3.2) {\textbf{Symmetric}};
\node at (-0.1,-3.375) {\textbf{tensor dofs}};
\end{scope}
\end{scope}
\node at (-2.5,-0.5) {\color{orange} \textbf{Solid material 1}};
\node at (-2.5,-0.75) {\color{orange} \textbf{Cell} ${T^1 \in \Ts}$};
\node at (2.5,-0.5) {\color{gray} \textbf{Solid material 2}};
\node at (2.5,-0.75) {\color{gray} \textbf{Cell} ${T^2 \in \Ts}$};
\def\a{1} 
\def\b{0.87} 
\def\c{0.5}  
\begin{scope}
\draw[blue, line width=2pt, fill=blue!25] (0,\a) -- (\b,\c) -- (\b,-\c) -- (0,-\a) -- (-\b,-\c) -- (-\b,\c) -- cycle;
\end{scope}
\begin{scope}[shift={(0.875,-1.505)}]
\draw[ceared, line width=2pt, fill=gray!75] (0,\a) -- (\b,\c) -- (\b,-\c) -- (0,-\a) -- (-\b,-\c) -- (-\b,\c) -- cycle;
\end{scope}
\begin{scope}[shift={(-0.875,-1.505)}]
\draw[ceared, line width=2pt, fill=orange!50] (0,\a) -- (\b,\c) -- (\b,-\c) -- (0,-\a) -- (-\b,-\c) -- (-\b,\c) -- cycle;
\end{scope}
\draw[black, line width=2.5pt] (0,-\a) -- (-2*\b,0.0025) -- (-2.5*\a,0.0025) ;  
\draw[black, line width=2.5pt] (0,-\a) -- (2*\b,0.0025) -- (2.5*\a,0.0025);  
\begin{scope}[scale=0.4]
\begin{scope}[yshift=-1.5cm, xshift=-0.9cm]
\draw[black, line width=2pt, fill=blue] ( 0.9, 2.0) circle (5pt);
\draw[black, line width=2pt, fill=blue] ( 1.7, 2.5) circle (5pt);
\draw[black, line width=2pt, fill=blue] ( 0.1, 2.5) circle (5pt);
\node at (1, 3) {\Large \color{blue} $p_{T}$};  
\node at (4, 1.5) {\Large \color{blue} $p_{{\partial T}}$};  
\end{scope} 
\begin{scope}[yshift=-3.5cm, xshift=-0.9cm]
\draw[black, line width=2pt, fill=blue] ( 1.1, 3.15)   circle (5pt);
\draw[black, line width=2pt, fill=blue] ( 0.7, 3.15)   circle (5pt);
\draw[black, line width=2pt, fill=blue] ( 1.7, 2.55)  circle (5pt);
\draw[black, line width=2pt, fill=blue] ( 2.1, 2.55)  circle (5pt);
\draw[black, line width=2pt, fill=blue] ( 0.1, 2.55)  circle (5pt);
\draw[black, line width=2pt, fill=blue] (-0.3, 2.55) circle (5pt);
\node at (1,2) {\Large \color{blue} $\bd{m}_{T}$};  
\end{scope} 
\filldraw[black, line width=1.5pt] (-0.75, 2.075) node[shape=diamond, draw, fill=blue, inner sep=0pt, minimum size=10pt, rotate=30] {};
\filldraw[black, line width=1.5pt] (-1.5, 1.65) node[shape=diamond, draw, fill=blue, inner sep=0pt, minimum size=10pt, rotate=30] {};
\filldraw[black, line width=1.5pt] (0.75, 2.075) node[shape=diamond, draw, fill=blue, inner sep=0pt, minimum size=10pt, rotate=65] {};
\filldraw[black, line width=1.5pt] (1.5, 1.65) node[shape=diamond, draw, fill=blue, inner sep=0pt, minimum size=10pt, rotate=65] {};
\filldraw[black, line width=1.5pt] (2.175, -0.5) node[shape=diamond, draw, fill=blue, inner sep=0pt, minimum size=10pt] {};
\filldraw[black, line width=1.5pt] (2.175, 0.5) node[shape=diamond, draw, fill=blue, inner sep=0pt, minimum size=10pt] {};
\filldraw[black, line width=1.5pt] (-2.175, -0.5) node[shape=diamond, draw, fill=blue, inner sep=0pt, minimum size=10pt] {};
\filldraw[black, line width=1.5pt] (-2.175, 0.5) node[shape=diamond, draw, fill=blue, inner sep=0pt, minimum size=10pt] {};
\begin{scope}[yshift=-4.75cm]
\begin{scope}[xshift=1.25cm]
\filldraw[black, line width=1.5pt] (-0.85, 2.5) node[shape=diamond, draw, fill=blue, inner sep=0pt, minimum size=10pt, rotate=35] {};
\filldraw[black, line width=1.5pt] (-0.45, 2.7) node[shape=diamond, draw, fill=ceared, inner sep=0pt, minimum size=10pt, rotate=35] {};
\filldraw[black, line width=1.5pt] (-0.25, 2.8) node[shape=diamond, draw, fill=ceared, inner sep=0pt, minimum size=10pt, rotate=35] {};
\filldraw[black, line width=1.5pt] (0.15, 3.05) node[shape=diamond, draw, fill=blue, inner sep=0pt, minimum size=10pt, rotate=35] {};
\filldraw[black, line width=1.5pt] (0.55, 3.25) node[shape=diamond, draw, fill=ceared, inner sep=0pt, minimum size=10pt, rotate=35] {};
\filldraw[black, line width=1.5pt] (0.75, 3.35) node[shape=diamond, draw, fill=ceared, inner sep=0pt, minimum size=10pt, rotate=35] {};
\filldraw[black, line width=1.5pt] (2.7, 2.5) node[shape=diamond, draw, fill=ceared, inner sep=0pt, minimum size=10pt, rotate=65] {};
\filldraw[black, line width=1.5pt] (2.5, 2.6) node[shape=diamond, draw, fill=ceared, inner sep=0pt, minimum size=10pt, rotate=65] {};
\filldraw[black, line width=1.5pt] (1.8, 3.05) node[shape=diamond, draw, fill=ceared, inner sep=0pt, minimum size=10pt, rotate=65] {};
\filldraw[black, line width=1.5pt] (1.6, 3.15) node[shape=diamond, draw, fill=ceared, inner sep=0pt, minimum size=10pt, rotate=65] {};
\filldraw[black, line width=1.5pt] (-1.25, 0.30)  node[shape=diamond, draw, fill=ceared, inner sep=0pt, minimum size=10pt] {};
\filldraw[black, line width=1.5pt] (-1.25, 0.55)  node[shape=diamond, draw, fill=ceared, inner sep=0pt, minimum size=10pt] {};
\filldraw[black, line width=1.5pt] (-1.25, 1.40)  node[shape=diamond, draw, fill=ceared, inner sep=0pt, minimum size=10pt] {};
\filldraw[black, line width=1.5pt] (-1.25, 1.65)  node[shape=diamond, draw, fill=ceared, inner sep=0pt, minimum size=10pt] {};
\filldraw[black, line width=1.5pt] (3.12, 0.30)  node[shape=diamond, draw, fill=ceared, inner sep=0pt, minimum size=10pt] {};
\filldraw[black, line width=1.5pt] (3.12, 0.55)  node[shape=diamond, draw, fill=ceared, inner sep=0pt, minimum size=10pt] {};
\filldraw[black, line width=1.5pt] (3.12, 1.40)  node[shape=diamond, draw, fill=ceared, inner sep=0pt, minimum size=10pt] {};
\filldraw[black, line width=1.5pt] (3.12, 1.65)  node[shape=diamond, draw, fill=ceared, inner sep=0pt, minimum size=10pt] {};
\filldraw[black, line width=1.5pt] (-0.75, -0.55) node[shape=diamond, draw, fill=ceared, inner sep=0pt, minimum size=10pt, rotate=55] {};
\filldraw[black, line width=1.5pt] (-0.55, -0.65) node[shape=diamond, draw, fill=ceared, inner sep=0pt, minimum size=10pt, rotate=55] {};
\filldraw[black, line width=1.5pt] (0.25, -1.125) node[shape=diamond, draw, fill=ceared, inner sep=0pt, minimum size=10pt, rotate=55] {};
\filldraw[black, line width=1.5pt] (0.45, -1.225) node[shape=diamond, draw, fill=ceared, inner sep=0pt, minimum size=10pt, rotate=55] {};
\filldraw[black, line width=1.5pt] (2.475, -0.65) node[shape=diamond, draw, fill=ceared, inner sep=0pt, minimum size=10pt, rotate=35] {};
\filldraw[black, line width=1.5pt] (2.675, -0.55) node[shape=diamond, draw, fill=ceared, inner sep=0pt, minimum size=10pt, rotate=35] {};
\filldraw[black, line width=1.5pt] (1.65, -1.125) node[shape=diamond, draw, fill=ceared, inner sep=0pt, minimum size=10pt, rotate=35] {};
\filldraw[black, line width=1.5pt] (1.45, -1.225) node[shape=diamond, draw, fill=ceared, inner sep=0pt, minimum size=10pt, rotate=35] {};
\draw[black, line width=2pt, fill=ceared] ( 1.10, 1.4)   circle (5pt);
\draw[black, line width=2pt, fill=ceared] ( 0.70, 1.4)   circle (5pt);
\draw[black, line width=2pt, fill=ceared] ( 1.8,  1.9)  circle (5pt);
\draw[black, line width=2pt, fill=ceared] ( 2.2,  1.9)  circle (5pt);
\draw[black, line width=2pt, fill=ceared] ( 0.0,  1.9)  circle (5pt);
\draw[black, line width=2pt, fill=ceared] (-0.4,  1.9) circle (5pt);
\node at (1, 2.5) {\Large \color{ceared} $\bd{v}_{T^2}$};
\draw[black, line width=2pt, fill=ceared] (0.70,  0.75) circle (5pt);    
\draw[black, line width=2pt, fill=ceared] (1.10,  0.75) circle (5pt);    
\draw[black, line width=2pt, fill=ceared] (1.10,  0.35) circle (5pt);  
\draw[black, line width=2pt, fill=ceared] ( 0.0,  0.25) circle (5pt);    
\draw[black, line width=2pt, fill=ceared] (-0.4,  0.25) circle (5pt);    
\draw[black, line width=2pt, fill=ceared] ( 0.0, -0.15) circle (5pt);  
\draw[black, line width=2pt, fill=ceared] ( 2.2,  0.25) circle (5pt);    
\draw[black, line width=2pt, fill=ceared] ( 1.8,  0.25) circle (5pt);    
\draw[black, line width=2pt, fill=ceared] ( 2.2, -0.15) circle (5pt);  
\node at (1,-0.5) {\Large \color{ceared} $\bbm{s}_{T^2}$};
\node at (4.25,1) {\Large \color{ceared} $\bd{v}_{\partial T^2}$};
\end{scope}
\end{scope}
\begin{scope}[yshift=-4.75cm]
\begin{scope}[xshift=-3.15cm]
\filldraw[black, line width=1.5pt] (2.85, 2.45) node[shape=diamond, draw, fill=blue, inner sep=0pt, minimum size=10pt, rotate=65] {};
\filldraw[black, line width=1.5pt] (2.45, 2.65) node[shape=diamond, draw, fill=ceared, inner sep=0pt, minimum size=10pt, rotate=65] {};
\filldraw[black, line width=1.5pt] (2.25, 2.75) node[shape=diamond, draw, fill=ceared, inner sep=0pt, minimum size=10pt, rotate=65] {};
\filldraw[black, line width=1.5pt] (1.85, 3.0) node[shape=diamond, draw, fill=blue, inner sep=0pt, minimum size=10pt, rotate=65] {};
\filldraw[black, line width=1.5pt] (1.45, 3.25) node[shape=diamond, draw, fill=ceared, inner sep=0pt, minimum size=10pt, rotate=65] {};
\filldraw[black, line width=1.5pt] (1.25, 3.35) node[shape=diamond, draw, fill=ceared, inner sep=0pt, minimum size=10pt, rotate=65] {};
\filldraw[black, line width=1.5pt] (-0.75, 2.5) node[shape=diamond, draw, fill=ceared, inner sep=0pt, minimum size=10pt, rotate=35] {};
\filldraw[black, line width=1.5pt] (-0.55, 2.6) node[shape=diamond, draw, fill=ceared, inner sep=0pt, minimum size=10pt, rotate=35] {};
\filldraw[black, line width=1.5pt] (0.25, 3.05) node[shape=diamond, draw, fill=ceared, inner sep=0pt, minimum size=10pt, rotate=35] {};
\filldraw[black, line width=1.5pt] (0.45, 3.15) node[shape=diamond, draw, fill=ceared, inner sep=0pt, minimum size=10pt, rotate=35] {};
\filldraw[black, line width=1.5pt] (-1.215, 0.30)  node[shape=diamond, draw, fill=ceared, inner sep=0pt, minimum size=10pt] {};
\filldraw[black, line width=1.5pt] (-1.215, 0.55)  node[shape=diamond, draw, fill=ceared, inner sep=0pt, minimum size=10pt] {};
\filldraw[black, line width=1.5pt] (-1.215, 1.40)  node[shape=diamond, draw, fill=ceared, inner sep=0pt, minimum size=10pt] {};
\filldraw[black, line width=1.5pt] (-1.215, 1.65)  node[shape=diamond, draw, fill=ceared, inner sep=0pt, minimum size=10pt] {};
\filldraw[black, line width=1.5pt] (-0.75, -0.55) node[shape=diamond, draw, fill=ceared, inner sep=0pt, minimum size=10pt, rotate=55] {};
\filldraw[black, line width=1.5pt] (-0.55, -0.65) node[shape=diamond, draw, fill=ceared, inner sep=0pt, minimum size=10pt, rotate=55] {};
\filldraw[black, line width=1.5pt] (0.25, -1.125) node[shape=diamond, draw, fill=ceared, inner sep=0pt, minimum size=10pt, rotate=55] {};
\filldraw[black, line width=1.5pt] (0.45, -1.225) node[shape=diamond, draw, fill=ceared, inner sep=0pt, minimum size=10pt, rotate=55] {};
\filldraw[black, line width=1.5pt] (2.475, -0.65) node[shape=diamond, draw, fill=ceared, inner sep=0pt, minimum size=10pt, rotate=35] {};
\filldraw[black, line width=1.5pt] (2.675, -0.55) node[shape=diamond, draw, fill=ceared, inner sep=0pt, minimum size=10pt, rotate=35] {};
\filldraw[black, line width=1.5pt] (1.65, -1.125) node[shape=diamond, draw, fill=ceared, inner sep=0pt, minimum size=10pt, rotate=35] {};
\filldraw[black, line width=1.5pt] (1.45, -1.225) node[shape=diamond, draw, fill=ceared, inner sep=0pt, minimum size=10pt, rotate=35] {};
\draw[black, line width=2pt, fill=ceared] ( 1.10, 1.4)   circle (5pt);
\draw[black, line width=2pt, fill=ceared] ( 0.70, 1.4)   circle (5pt);
\draw[black, line width=2pt, fill=ceared] ( 1.8,  1.9)  circle (5pt);
\draw[black, line width=2pt, fill=ceared] ( 2.2,  1.9)  circle (5pt);
\draw[black, line width=2pt, fill=ceared] ( 0.0,  1.9)  circle (5pt);
\draw[black, line width=2pt, fill=ceared] (-0.4,  1.9) circle (5pt);
\node at (1, 2.5) {\Large \color{ceared} $\bd{v}_{T^1}$};
\draw[black, line width=2pt, fill=ceared] (0.70,  0.75) circle (5pt);    
\draw[black, line width=2pt, fill=ceared] (1.10,  0.75) circle (5pt);    
\draw[black, line width=2pt, fill=ceared] (1.10,  0.35) circle (5pt);  
\draw[black, line width=2pt, fill=ceared] ( 0.0,  0.25) circle (5pt);    
\draw[black, line width=2pt, fill=ceared] (-0.4,  0.25) circle (5pt);    
\draw[black, line width=2pt, fill=ceared] ( 0.0, -0.15) circle (5pt);  
\draw[black, line width=2pt, fill=ceared] ( 2.2,  0.25) circle (5pt);    
\draw[black, line width=2pt, fill=ceared] ( 1.8,  0.25) circle (5pt);    
\draw[black, line width=2pt, fill=ceared] ( 2.2, -0.15) circle (5pt);  
\node at (1,-0.5) {\Large \color{ceared} $\bbm{s}_{T^1}$};
\node at (-2.25,1) {\Large \color{ceared} $\bd{v}_{\partial T^1}$};
\node at (8.5,5.5) {\huge \color{black} $\bd{\G}$};
\end{scope}
\end{scope}
\end{scope}
\end{tikzpicture}}
\caption{Elasto-acoustic dofs along a solid-fluid interface, with hexagonal mesh cells and lowest equal-order discretization $(k^\prime = k = 1)$.}
\label{dofs}
\end{figure}
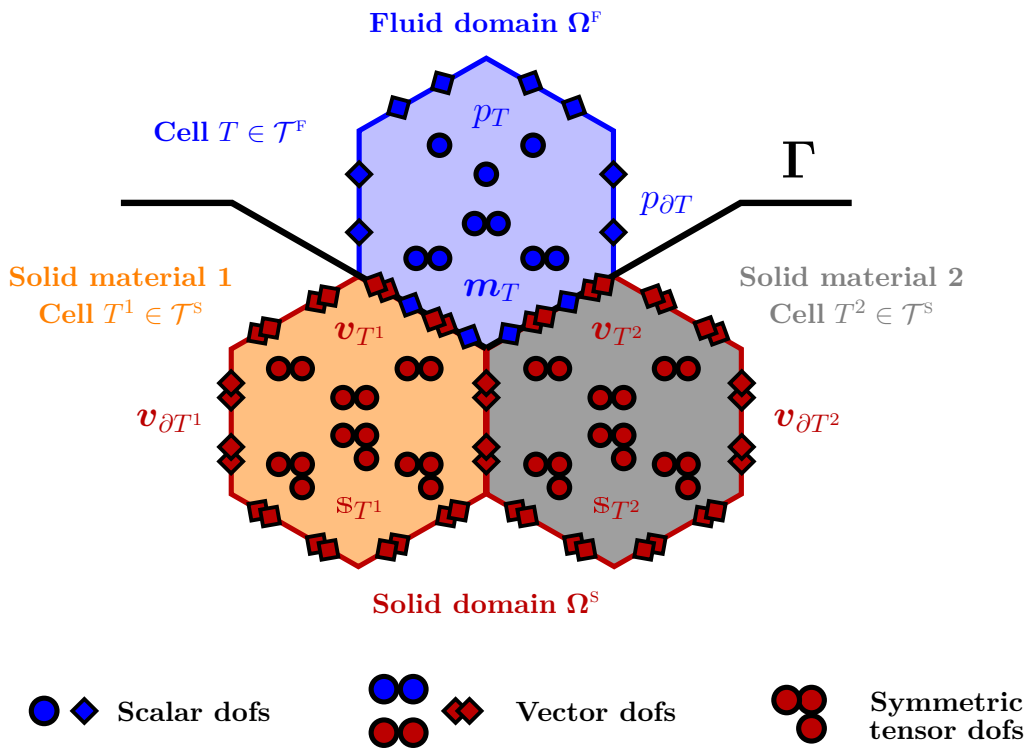  
\fi    
\subsection{Space semi-discrete formulation and algebraic realization} 

The space semi-discretization of the model problem \eqref{weak_form_acoustic_eq}-\eqref{weak_form_elastic_eq} consists in finding $(p_{\Tf}(t),$ $ p_{\Ff}(t), \bd{m}_{\Tf}(t)) \in P^{k^\prime}(\Tf) \times P^{k}(\Ff) \times \bd{M}^{k}(\Tf)$ in the fluid subdomain and $(\bd{v}_{\Ts}(t), \bd{v}_{\Fs}(t),$ $\bbm{s}_{\Ts}(t)) \in \bd{V}^{k^\prime}(\Ts) \times \bd{V}^{k}(\Fs) \times \bbm{S}^{k}(\Ts)$ in the solid subdomain, such that, for all $t \in \overline{J}$,
\begin{subequations}
\label{HHO1}
\begin{align}
(\partial_t \bd{m}_{\Tf}(t) & , \bd{r}_{\Tf})_{{\bd{L}^2(\rho^\sc{f};\domain{f})}} - (\bd{g}_{\Tf}(p_{\Tf}(t),p_{\Ff}(t)), \bd{r}_{\Tf})_{{\bd{L}^2(\domain{f})}} = 0,\\
(\partial_t p_{\Tf}(t) & ,q_{\Tf})_{{L^2(\frac{1}{\kappa}; \domain{f})}} + (\bd{m}_{\Tf}(t), \bd{g}_{\Tf}(q_{\Tf},q_{\Ff}))_{{\bd{L}^2(\domain{f})}} \nonumber\\
& + s^{\sc{f}}((p_{\Tf}(t),p_{\Ff}(t)),(q_{\Tf},q_{\Ff})) + (\bd{v}_{\Fs}(t) {\cdot} \bd{n}_{\G}, q_{\Ff})_{L^2(\G)} = (f^{\sc{f}}(t),q_{\Tf})_{L^2(\domain{f})}, 
\end{align}
\end{subequations}  
for all $\bd{r}_{\Tf} \in \bd{M}^{k}(\Tf)$ and all $(q_{\Tf}, q_{\Ff}) \in P^{k^\prime}(\Tf) \times P^{k}(\Ff)$, and 
\begin{subequations}
\label{HHO2} 
\begin{align}
(\partial_t & \bbm{s}_{\Ts}(t), \bbm{b}_{\Ts})_{\bbm{L}^2(\bbm{C}^{-1};\domain{s})} - (\bbm{g}^{\rm{sym}}_{\Ts}(\bd{v}_{\Ts}(t),\bd{v}_{\Fs}(t)), \bbm{b}_{\Ts})_{\bbm{L}^2(\domain{s})} = 0,\\
(\partial_t & \bd{v}_{\Ts}(t), \bd{w}_{\Ts})_{\bd{L}^2(\rho^\sc{s};\domain{s})} + (\bbm{s}_{\Ts}(t), \bbm{g}_{\Ts}^{\rm{sym}}(\bd{w}_{\Ts},\bd{w}_{\Fs}))_{\bbm{L}^2(\domain{s})} \nonumber\\
& + s^{\sc{s}}((\bd{v}_{\Ts}(t),\bd{v}_{\Fs}(t)), (\bd{w}_{\Ts},\bd{w}_{\Fs})) -(p_{\Ff}(t) \bd{n}_{\G},\bd{w}_{\Fs})_{\bd{L}^2(\G)} = (\bd{f^{\sc{s}}}(t), \bd{w}_{\Ts})_{\bd{L}^2(\domain{s})},
\end{align}
\end{subequations}   
for all $\bbm{b}_{\Ts} \in \bbm{S}^k(\Ts)$ and all $(\bd{w}_{\Ts},\bd{w}_{\Fs}) \in \bd{V}^{k^\prime}(\Ts) \times \bd{V}^{k}(\Fs)$. Here, $\bd{g}_{\Tf}$, $\bbm{g}^{\rm{sym}}_{\Ts}$, $s^{\sc{f}}$, and $s^{\sc{s}}$ are the global (symmetric) gradient reconstruction operators and global stabilization bilinear forms whose expressions are given by the right-hand sides of equations \eqref{rec}-\eqref{stab} below. The main ideas are the following: the (symmetric) gradient reconstruction allows to consistently and locally build the (symmetric) gradient of the primal variable using both its cell and face dofs; the stabilization bilinear form weakly enforces locally the matching of the trace of the cell dofs with the face dofs. The reader is referred to \cite{MEKG_2025} for further insight into the HHO formulation.

To derive the algebraic realization of \eqref{HHO1}-\eqref{HHO2}, we consider the following generic bases:
\begin{subequations}
\begin{alignat}{5}
& \{\varphi_i\}_{1 \leq i \leq N_{\Tf}^{k^\prime}}, & \qquad & \{\psi_i \}_{1 \leq i \leq N_{\Ff}^{k}}, & \qquad & \{\bd{\zeta}_k\}_{1 \leq k \leq M_{\Tf}^{k}},\\
& \{\bd{\phi}_i\}_{1 \leq i \leq L_{\Ts}^{k^\prime}}, & \qquad & \{\bd{\theta}_i \}_{1 \leq i \leq L_{\Fs}^{k}}, & \qquad & \{\bbm{Y}_k\}_{1 \leq k \leq H_{\Ts}^{k}},
\end{alignat}
\end{subequations}
with the corresponding dimensions
\begin{subequations}
\begin{alignat}{5}
& N_{\Tf}^{k^\prime} := \dim(P^{k^\prime}(\Tf)), & \qquad & N_{\Ff}^{k} := \dim(P^{k}(\Ff)), & \qquad & M_{\Tf}^{k} := \dim(\bd{M}^{k}(\Tf)),\\
& L_{\Ts}^{k^\prime} := \dim(\bd{{V}}^{k^\prime}(\Ts)), & \qquad & L_{\Fs}^{k} := \dim(\bd{{V}}^{k}(\Fs)), & \qquad & H_{\Ts}^{k} := \dim(\bbm{S}^k(\Ts)).
\end{alignat}
\end{subequations}
The basis $\{\bd{\zeta}_k\}_{1 \leq k \leq M_{\Tf}^{k}}$ is build by combining Cartesian basis vectors of $\bbm{R}^d$ with scalar basis functions of $P^{k}(\Tf)$. Similarly, the bases $\{\bd{\phi}_i \}_{1 \leq i \leq L_{\Ts}^{k^\prime}}$, $\{\bd{\theta}_i \}_{1 \leq i \leq L_{\Fs}^{k}}$ and $\{\bbm{Y}_k\}_{1 \leq k \leq H_{\Ts}^{k}}$ are built by combining Cartesian basis vectors of $\bbm{R}^{d}$, $\bbm{R}^{d}$ and $\bbm{R}_{\text{sym}}^{d \times d}$, respectively, with scalar basis functions of $P^{k^\prime}(\Ts)$, $P^{k}(\Fs)$ and $P^{k}(\Ts)$, respectively. For all the scalar-valued polynomial spaces, we use a local modal basis composed of monomials centered at the barycenter of the corresponding geometric object.

\paragraph{Degrees of freedom.} Let $(\dofs{P}{T}{f}(t), \dofs{P}{F}{f}(t)) \in \bbm{R}^{N_{\Tf}^{k^\prime} \times N_{\Ff}^{k}}$ and $\dofs{M}{T}{f}(t) \in \bbm{R}^{M_{\Tf}^{k}}$ represent the time-dependent component vectors of $(p_{\Tf}(t), p_{\Ff}(t)) \in P^{k^\prime}(\Tf) \times P^{k}(\Ff)$ and $\bd{m}_\Tf(t) \in \bd{M}^{k}(\Tf)$ in the fluid subdomain $\domain{f}$, respectively. Let $(\dofs{V}{T}{s}(t), \dofs{V}{F}{s}(t)) \in \bbm{R}^{L_{\Ts}^{k^\prime} \times L_{\Fs}^{k}}$ and $\dofs{S}{T}{s} \in \bbm{R}^{H_{\Ts}^{k}}$ represent the time-dependent component vectors of $(\bd{v}_{\Ts}(t), \bd{v}_{\Fs}(t)) \in \bd{{V}}^{k^\prime}(\Ts) \times \bd{{V}}^{k}(\Fs)$ and $\bbm{s}_{\Ts} \in \bbm{S}^k(\Ts)$ in the solid subdomain $\domain{s}$, respectively. The algebraic realization of \eqref{HHO1}-\eqref{HHO2} is, for all $t \in \overline{J}$, 
\begin{equation}
\resizebox{0.925\textwidth}{!}{%
$
\pmatset{4}{12pt} 
\pmatset{5}{10pt} 
\resizebox{!}{13.5ex}{
$
\begin{pmat}[{.|.|}]
\normalizecell{\mass{\rho^{\sc{f}}}{f}} & 0      & 0      & 0      & 0      & 0 \cr
0      & \normalizecell{\mass{\frac{1}{\kappa}}{f}} & 0      & 0      & 0      & 0 \cr \-
0      & 0      & \normalizecell{\mass{\bbm{C}^{-1}}{s}} & 0      & 0      & 0 \cr 
0      & 0      & 0      & \normalizecell{\mass{\rho^{\sc{s}}}{s}} & 0      & 0 \cr \-
0      & 0      & 0      & 0      & 0      & 0 \cr
0      & 0      & 0      & 0      & 0      & 0 \cr
\end{pmat}
$
}
\dfrac{\rm{d}}{\rm{dt}}
\pmatset{4}{12pt} 
\pmatset{5}{10pt} 
\begin{pmat}[{}]
\normalizecell{\dofs{M}{T}{f}} \cr
\normalizecell{\dofs{P}{T}{f}} \cr \-
\normalizecell{\dofs{S}{T}{s}} \cr
\normalizecell{\dofs{V}{T}{s}} \cr \-
\normalizecell{\dofs{P}{F}{f}} \cr 
\normalizecell{\dofs{V}{F}{s}} \cr
\end{pmat} +
\pmatset{4}{12pt} 
\pmatset{5}{10pt} 
\begin{pmat}[{.|.|}]
0 & \normalizecell{\graddiag{T}{f}} & 0 & 0 & \normalizecell{\grad{T}{F}{f}} & 0 \cr
- \normalizecell{\graddiag{T}{f}^{\dagger}} & \normalizecell{\stabdiag{f}{T}} & 0 & 0 & \normalizecell{\stab{f}{T}{F}} & 0 \cr \- 
0 & 0 & 0 & \normalizecell{\straindiag{T}} & 0 & \normalizecell{\strain{T}{F}}\cr
0 & 0 & -\normalizecell{\straindiag{T}^{\dagger}} & \normalizecell{\stabdiag{s}{T}} & 0 & \normalizecell{\stab{s}{T}{F}} \cr \- 
-\normalizecell{\grad{T}{F}{f}^{\dagger}} & \normalizecell{\stab{f}{T}{F}^{\dagger}} & 0 & 0 & \normalizecell{\stabdiag{f}{F}} & \normalizecell{\coupling} \cr
0 & 0 & -\normalizecell{\strain{T}{F}^\dagger} & \normalizecell{\stab{s}{T}{F}^{\dagger}} & -\normalizecell{\coupling^{\dagger}} & \normalizecell{\stabdiag{s}{F}} \cr
\end{pmat}
\pmatset{4}{12pt} 
\pmatset{5}{10pt} 
\begin{pmat}[{}]
\normalizecell{\dofs{M}{T}{f}} \cr
\normalizecell{\dofs{P}{T}{f}} \cr \-
\normalizecell{\dofs{S}{T}{s}} \cr
\normalizecell{\dofs{V}{T}{s}} \cr \-
\normalizecell{\dofs{P}{F}{f}} \cr 
\normalizecell{\dofs{V}{F}{s}} \cr
\end{pmat} =
\pmatset{4}{12pt} 
\pmatset{5}{10pt} 
\begin{pmat}[{}]
\normalizecell{0}  \cr
\normalizecell{\dofs{F}{T}{f}} \cr \-
\normalizecell{0}  \cr 
\normalizecell{\dofs{F}{T}{s}} \cr \-
\normalizecell{0}  \cr
\normalizecell{0}  \cr
\end{pmat} 
$,
}
\label{algebraic_system}
\end{equation}    
where the different blocks are defined as follows. 

\paragraph{Mass matrices.} The mass matrices $\mass{\rho^{\sc{f}}}{f} \in \bbm{R}^{M_{\Tf}^{k} \times M_{\Tf}^{k}}$ and $\mass{\frac{1}{\kappa}}{f} \in \bbm{R}^{N_{\Tf}^{k^\prime} \times N_{\Tf}^{k^\prime}}$ are associated with the inner products in $\bd{L}^2(\rho^\sc{f}; \Omega^\sc{f})$ and $L^2(\frac{1}{\kappa}; \Omega^\sc{f})$, respectively, and the mass matrices $\mass{\bbm{C}^{-1}}{s} \in \bbm{R}^{H_{\Ts}^{k} \times H_{\Ts}^{k}}$ and $\mass{\rho^{\sc{s}}}{s} \in \bbm{R}^{L_{\Ts}^{k^\prime} \times L_{\Ts}^{k^\prime}}$ are associated with the inner products in $\bbm{L}^2_{\rm{sym}}(\bbm{C}^{-1};\Omega^\sc{s})$ and $\bd{L}^2(\rho^\sc{s}; \Omega^\sc{s})$, respectively. Their (classical) expressions are as follows:
\begin{subequations}\label{mass}
\begin{alignat}{6}
& \dofs{M}{T}{f}^{\dagger} \mass{\rho^{\sc{f}}}{f} \rm{R}_{\T^\sc{f}} 
&& = \sum_{T \in \Tf} \rho^{\sc{f}}|_T(\bd{m}_T, \bd{r})_{\bd{L}^2(T)}, \qquad
&& \dofs{P}{T}{f}^{\dagger} \mass{\frac{1}{\kappa}}{f} \rm{Q}_{\T^\sc{f}} 
&& = \sum_{T \in \Tf} \frac{1}{\kappa|_T}({p}_T, q)_{{L}^2(T)}, \\  
& \dofs{S}{T}{s}^{\dagger} \mass{\bbm{C}^{-1}}{s} \rm{B}_{\T^\sc{s}} 
&& = \sum_{T \in \Ts} \bbm{C}^{-1}|_T(\bbm{s}_T, \bbm{b})_{\bbm{L}^2(T)}, \qquad
&& \dofs{V}{T}{s}^{\dagger} \mass{\rho^{\sc{s}}}{s} \rm{W}_{\T^\sc{s}} 
&& = \sum_{T \in \Ts} \rho^{\sc{s}}|_T(\bd{v}_T, \bd{w})_{\bd{L}^2(T)},
\end{alignat}                 
\end{subequations}
for all $\bd{m}_{\Tf}, \bd{r} \in \bd{P}^{k}(\Tf)$ and all $p_{\Tf}, q \in {P}^{k^\prime}(\Tf)$ with component vectors $\dofs{M}{T}{f}$, $\rm{R}_{\T^\sc{f}}$ and $\dofs{P}{T}{f}$, $\rm{Q}_{\T^\sc{f}}$, respectively, and all $\bbm{s}_{\Ts}, \bbm{b} \in \bbm{S}^{k}(\Ts)$ and all $\bd{v}_{\Ts}, \bd{w} \in \bd{P}^{k^\prime}(\Ts)$ with component vectors $\dofs{S}{T}{s}$, $\rm{B}_{\T^\sc{s}}$ and $\dofs{V}{T}{s}$, $\rm{W}_{\T^\sc{s}}$, respectively.

\paragraph{Gradient reconstruction blocks.} The gradient reconstruction blocks $\graddiag{T}{f} \in \bbm{R}^{M_{\Tf}^{k} \times N_{\Tf}^{k^\prime}}$ and $\grad{T}{F}{f} \in \bbm{R}^{M_{\Tf}^{k} \times N_{\Ff}^{k}}$ and the symmetric gradient reconstruction blocks $\straindiag{\T} \in \bbm{R}^{H_{\Ts}^{k} \times L_{\Ts}^{k^\prime}}$ and $\strain{T}{F} \in \bbm{R}^{H_{\Ts}^{k} \times L_{\Fs}^{k}}$ are defined so that they satisfy
\begin{subequations}\label{rec}
\begin{alignat}{6}
\label{rec_a}
& (\graddiag{T}{f}\dofs{\rm{P}}{T}{f} + \grad{T}{F}{f} \dofs{\rm{P}}{F}{f})^{\dagger} \rm{R}_{\T^\sc{f}} 
&& = \sum_{T \in \Tf} \left\{(\nabla p_T, \bd{r})_{\bd{L}^2(T)} - (p_{T} - p_{\partial T}, \bd{r} {\cdot} \bd{n}_T)_{L^2(\partial T)} \right\}, \\  
& (\straindiag{T}\dofs{\rm{V}}{T}{s} + \strain{T}{F} \dofs{\rm{V}}{F}{s})^{\dagger} \rm{B}_{\T^\sc{s}} 
&& = \sum_{T \in \Ts} \left\{ (\nabla_{\rm{sym}} \bd{v}_T, \bbm{b})_{\bbm{L}^2_{\rm{sym}}(T)} - (\bd{v}_T - \bd{v}_{\partial T}, \bbm{b} {\cdot} \bd{n}_T)_{\bd{L}^2(\partial T)} \right\},
\label{rec_e}
\end{alignat}                 
\end{subequations}
for all $(p_{\Tf},p_{\Ff}) \in P^{k^\prime}(\Tf) \times P^{k}(\Ff)$, all $(\bd{v}_{\Ts},\bd{v}_{\Fs}) \in \bd{V}^{k^\prime}(\Tf) \times \bd{V}^{k}(\Ff)$ with component vectors $\dofs{P}{T}{f}$, $\dofs{P}{F}{f}$ and $\dofs{V}{T}{s}$, $\dofs{V}{F}{s}$, respectively, and for all $\bd{r} \in \bd{M}^{k}(\Tf)$, all $\bbm{b} \in \bbm{S}^k(\Ts)$ with component vectors $\rm{R}_{\Tf}$ and $\dofs{B}{T}{s}$, respectively.

\paragraph{Stabilization blocks.} The fluid stabilization blocks $\stabdiag{f}{T} \in \bbm{R}^{N_{\Tf}^{k^\prime} \times N_{\Tf}^{k^\prime}}$, $\stab{f}{T}{F} \in \bbm{R}^{N_{\Tf}^{k^\prime} \times N_{\Ff}^{k}}$, and $\stabdiag{f}{F} \in \bbm{R}^{N_{\Ff}^{k} \times N_{\Ff}^{k}}$ and the solid stabilization blocks $\stabdiag{s}{T} \in \bbm{R}^{L_{\Ts}^{k^\prime} \times L_{\Ts}^{k^\prime}}$, $\stab{s}{T}{F} \in \bbm{R}^{L_{\Ts}^{k^\prime} \times L_{\Fs}^{k}}$, and $\stabdiag{s}{F} \in \bbm{R}^{L_{\Fs}^{k} \times L_{\Fs}^{k}}$ are defined so that they satisfy 
\begin{subequations}\label{stab}
\begin{align}
\label{stab_a}
(\stabdiag{f}{T}\dofs{\rm{P}}{T}{f} + \stab{f}{T}{F} \dofs{\rm{P}}{F}{f})^\dagger \dofs{\rm{Q}}{T}{f} &+ (\stab{f}{T}{F}^{\dagger}\dofs{\rm{P}}{T}{f} + \stabdiag{f}{F}\dofs{\rm{P}}{F}{f})^\dagger \dofs{\rm{Q}}{F}{f} \notag \\
& = \sum_{T \in \Tf} \tau^\sc{f}_{T} (S_{\partial T}(p_T,p_{\partial T}), S_{\partial T}(q_T,q_{\partial T}))_{L^2(\partial T)}, \\
(\stabdiag{s}{T}\dofs{\rm{V}}{T}{s} + \stab{s}{T}{F}
\dofs{\rm{V}}{T}{s})^\dagger \dofs{\rm{W}}{T}{s} &+ (\stab{s}{T}{F}^{\dagger}\dofs{\rm{V}}{T}{s} + \stabdiag{s}{F}\dofs{\rm{V}}{F}{s})^\dagger \dofs{\rm{W}}{F}{s} \notag \\
& = \sum_{T \in \Ts} \tau^\sc{s}_{T} (\bd{S}_{\partial T}(\bd{v}_T,\bd{v}_{\partial T}), \bd{S}_{\partial T}(\bd{w}_T,\bd{w}_{\partial T}))_{\bd{L}^2(\partial T)},
\label{stab_e}
\end{align}
\end{subequations}    
for all $(p_{\Tf},p_{\Ff}) \in P^{k^\prime}(\Tf) \times P^{k}(\Ff)$ and all $(\bd{v}_{\Ts},\bd{v}_{\Fs}) \in \bd{V}^{k^\prime}(\Ts) \times \bd{V}^{k}(\Fs)$ with component vectors $\dofs{P}{T}{f}$, $\dofs{P}{F}{f}$ and $\dofs{V}{T}{s}$, $\dofs{V}{F}{s}$, respectively, and for all $(q_{\Tf}, q_{\Ff}) \in P^{k^\prime}(\Tf) \times P^{k}(\Ff)$ and all $(\bd{w}_{\Ts}, \bd{w}_{\Fs}) \in \bd{V}^{k^\prime}(\Ts) \times \bd{V}^{k}(\Fs)$ with component vectors $\rm{Q}_{\Tf}$, $\rm{Q}_{\F^\sc{f}}$ and $\rm{W}_{\Ts}$, $\rm{W}_{\Fs}$, respectively. The stabilization parameters are defined as follows:
\begin{equation}
\tau^\sc{f}_{T} := (\rho^\sc{f}c_\sc{p}^\sc{f})|_T^{-1}\tih_T^{-\alpha}, \qquad \tau^\sc{s}_{T} := (\rho^\sc{s}c^\sc{s})|_T \tih_T^{-\alpha},
\end{equation}
where the local scaling by $(\rho^\sc{f}c_\sc{p}^\sc{f})|_T^{-1}$ and $(\rho^\sc{s}c^\sc{s})|_T$ follows from physical consistency, the parameter $\alpha \in \{0,1\}$ is used to shift from $\cal{O}(1)$-stabilization ($\alpha=0$) to $\cal{O}(\frac{1}{h})$-stabilization ($\alpha=1$), and $\tih_T := \frac{h_T}{\ell_\Omega}$ is a normalized cell diameter, where the global length scale $\ell_\Omega := \text{diam}(\Omega)$ is introduced for dimensional consistency. \rev{The stabilization parameters can also be scaled by the reciprocal of the polynomial degree.}

Regarding the local stabilization operators $S_{\partial T} : P^{k'}(T) \times P^k(\F_{\partial T}) \rightarrow P^k(\F_{\partial T})$ and $\bd{S}_{\partial T} : \bd{P}^{k'}(T) \times \bd{P}^k(\F_{\partial T}) \rightarrow \bd{P}^k(\F_{\partial T})$, we consider two configurations:
\begin{itemize}
\item The first one, which we use in the equal-order HHO discretization ($k' = k$), corresponds to a plain Least-Squares stabilization:
\begin{equation}
S_{\partial T}(p_T, p_{\partial T}) := p_T|_{\partial T} - p_{\partial T}, \qquad \bd{S}_{\partial T}(\bd{v}_T, \bd{v}_{\partial T}) := \bd{v}_T|_{\partial T} - \bd{v}_{\partial T}.
\end{equation}
This stabilization will be used with the scaling $\alpha = 0$ and explicit time-stepping schemes.
\item The second one, which we use in the mixed-order HHO discretization ($k' = k+1$), corresponds to the Lehrenfeld--Schöberl stabilization \cite{LEHRENFELD_2010}:
\begin{equation}
S_{\partial T}(p_T, p_{\partial T}) := \Pi_{\partial T}^k(p_T|_{\partial T}) - p_{\partial T}, \qquad \bd{S}_{\partial T}(\bd{v}_T, \bd{v}_{\partial T}) := \Pi_{\partial T}^k(\bd{v}_T|_{\partial T}) - \bd{v}_{\partial T},
\end{equation}
where $\Pi_{\partial T}^k$ is the local $L^2(\partial T)$-orthogonal projection onto $P^{k}(\cF_{\partial T})$ or its vector-valued version. This stabilization will be used with the scaling $\alpha = 1$ and implicit time-stepping schemes.
\end{itemize}    

\paragraph{Coupling blocks.} The matrix $\coupling{} \in \bbm{R}^{N_{\Ff}^{k} \times L_{\Fs}^{k}}$ representing the coupling terms is defined so that 
\begin{equation}
\rm{Q}_{\F^\sc{f}}^\dagger \coupling{} \dofs{V}{F}{s} = (\bd{v}_{\Fs} {\cdot} \bd{n}_\Gamma, q_{\cal{F}^\sc{f}})_{{L^2(\Gamma)}},
\end{equation}
for all $q_{\cal{F}^\sc{f}} \in P^{k}(\Ff)$ with component vector $\rm{Q}_{\F^\sc{f}} \in \bbm{R}^{N_{\Ff}^{k}}$ and all $\bd{v}_{\Fs} \in \bd{P}^k(\cal{\Fs})$ with component vector $\dofs{\rm{V}}{F}{s} \in \bbm{R}^{L_{\Fs}^{k}}$. Notice that $\coupling{}$ is block-diagonal with a nonzero block only for all $F \in \cal{F}^{\Gamma}$.

\paragraph{Source terms.} The source terms $\dofs{F}{T}{f}(t) \in \bbm{R}^{N_{\Tf}^{k^\prime}}$ and $\dofs{F}{T}{s}(t) \in \bbm{R}^{L_{\Ts}^{k^\prime}}$ are time-dependent and correspond to the algebraic realization of the linear forms
\begin{subequations} \label{rhs}
\begin{alignat}{4}
& \dofs{F}{T}{f}^{\dagger}(t) \mass{}{f} \dofs{Q}{T}{f} && = \sum_{T \in \Tf} (f^\textsc{f}(t), q)_{L^2(T)}, \\
& \dofs{F}{T}{s}^{\dagger}(t) \mass{}{s} \dofs{W}{T}{s} && = \sum_{T \in \Ts} (\bd{f}^\textsc{s}(t), \bd{w})_{\bd{L}^2(T)},
\end{alignat}
\end{subequations}
for all $q \in P^{k^\prime}(\Tf)$ and all $\bd{w} \in \bd{P}^{k^\prime}(\Ts)$ with component vectors $\dofs{Q}{T}{f}$ and $\dofs{W}{T}{s}$, respectively, and $\mass{}{f}$ and $\mass{}{s}$ are the mass matrices associated with the canonical inner products in $L^2(\Omega^\sc{f})$ and $\bd{L}^2(\Omega^\sc{s})$, respectively.

Equation \eqref{algebraic_system} can be rewritten in a more compact form by regrouping all the cell unknowns on the one hand, and all the face unknowns on the other hand. With obvious notation, we obtain
\begin{equation}
\pmatset{4}{12pt} 
\pmatset{5}{12pt}
\begin{pmat}[{}]
\normalizecell{\mass{}{}} & 0 \cr
0 & 0  \cr 
\end{pmat}
\partial_t
\begin{pmat}[{}]
\rm{U}_\T \cr
\rm{U}_\F  \cr 
\end{pmat} + 
\pmatset{4}{12pt} 
\pmatset{5}{12pt}
\begin{pmat}[{}]
\normalizecell{\cal{K}_{\cal{T}}} & \normalizecell{\cal{K}_{\cal{TF}}} \cr
\normalizecell{{\cal{K}_{\cal{FT}}}} & \normalizecell{\cal{K}_{\cal{F}}}  \cr 
\end{pmat}
\begin{pmat}[{}]
\rm{U}_\T \cr
\rm{U}_\F  \cr 
\end{pmat} = 
\begin{pmat}[{}]
\rm{F}_\T \cr
0  \cr 
\end{pmat},
\label{compact_AR}
\end{equation}
where $\mass{}{}$ is the cell mass matrix, $\cal{K}_{\T}$, $\cal{K}_{\T\F}$, $\cal{K}_{\F\T}=\cal{K}_{\T\F}^\dagger$, $\cal{K}_{\F}$, are the blocks of the stiffness matrix, $\rm{U}_\T$ and $\rm{U}_\F$ are the components of the vector of unknowns and $\rm{F}_\T$ is the cell source term (the face source term is null). A crucial observation is that both $\cal{K}_{\T}$ and $\cal{K}_{\F}$ are block diagonal. This important property will be exploited in Section 4 in the context of implicit and explicit time-stepping schemes. While the block-diagonal structure of $\cal{K}_\T$ is obvious, let us briefly motivate why $\cal{K}_\F$ is also block-diagonal. First, this is the case for all $F \in \F^{\circ \sc{f}} \cup \F^{\circ \sc{s}}$ since the stabilization bilinear forms are local to each interior face for both fluid and solid subdomains, with local diagonal blocks of size $\dim{P^k(F)}$ or $\dim{\bd{P}^k(F)}$, respectively. Furthermore, for all $F \in \F^{\G}$, the couplings are also local to the face, but the diagonal block in $\cal{K}_\F$ is of larger size, namely $\dim{P^k(F)} + \dim{\bd{P}^k(F)}$.

\section{Fully discrete schemes}

The time discretization is based on Runge--Kutta (RK) schemes. Butcher tableaux are a classical way to define RK schemes by means of the coefficients $\left\{a_{ij}\right\}_{1 \leq i, j \leq s}$, $\left\{b_i\right\}_{1 \leq i \leq s}$, and $\left\{c_i\right\}_{1 \leq i \leq s}$, where $s \geq 1$ is the number of stages. For explicit RK schemes with $s$ stages and order $s$ (ERK$(s)$), the matrix $\left\{a_{ij}\right\}_{1 \leq i, j \leq s}$ is strictly lower-triangular. For singly diagonal implicit RK schemes with $s$ stages and order $(s+1)$ (SDIRK$(s,s+1)$), the matrix is lower-triangular with $a_{11} = \ldots = a_{ss} := a_*$. The corresponding Butcher tableaux are as follows:
\begin{equation}
\renewcommand{\arraystretch}{1.3}
\begin{tabular}{c|cccc}
$c_1$ & $a_*$ & $0$ & $\cdots$ & $0$ \\
$c_2$ & $a_{21}$ & $a_*$ & $\ddots$ & $0$ \\
$\vdots$ & $\vdots$ & $\ddots$ & $\ddots$ & $\vdots$ \\
$c_s$ & $a_{s1}$ & $\cdots$ & $a_{s,s-1}$ & $a_*$ \\
\hline
& $b_1$ & $\cdots$ & $b_{s-1}$ & $b_s$
\end{tabular}
\hspace*{2cm}
\begin{tabular}{c|cccc}
$c_1$ & $0$ & $\cdots$ & $\cdots$ & $0$ \\
$c_2$ & $a_{21}$ & $0$ & $\cdots$ & $0$ \\
$\vdots$ & $\vdots$ & $\ddots$ & $\ddots$ & $\vdots$ \\
$c_s$ & $a_{s1}$ & $\cdots$ & $a_{s,s-1}$ & $0$ \\
\hline
& $b_1$ & $\cdots$ & $b_{s-1}$ & $b_s$
\end{tabular}
\end{equation}
To simplify the writing of SDIRK and ERK schemes, we define $a_{s+1,i} := b_i$ for all $i \in \{1,...,s\}$ and set $a_{s+1,s+1} := 0$. 

We now detail the time discretization of \eqref{compact_AR} by SDIRK and ERK schemes. We show that an effective SDIRK implementation is obtained by the local elimination of the cell unknowns, whereas an effective ERK implementation is obtained by the local elimination of the face unknowns. The outcome of this elimination is illustrated in \hyperref[fig::sc]{\Cref{fig::sc}}. Let $(t^n)_{n \in \{0,\ldots,N\}}$ with $t^0=0$ and $t^N=T_\rm{f}$ be the discrete time nodes, and let $\Delta t^n:=t^n-t^{n-1}$ be the time step. For simplicity, we consider that the time-step is constant and drop the superscript $n$.
\begin{figure}[!htb]
\centering
\resizebox{0.97\textwidth}{!}{%
\begin{tikzpicture}[scale=2]


\begin{scope}[xshift=-7cm, yshift=0.2cm]
\draw[black, line width=0.5pt, fill=ceared!25] (2.5,0) -- (3.5,0) -- (3.5,1) -- (2.5,1) -- cycle;
\draw[black, line width=0.5pt, fill=blue!25] (3.5,0) -- (4.5,0) -- (4.5,1) -- (3.5,1) -- cycle;
\draw[black, line width=0.5pt] (4.0, 0.00) -- (4.0, 1.00);
\draw[black, line width=0.5pt] (3.5, 0.50) -- (4.5, 0.50);

\filldraw[black, line width=0.5pt] (3.5, 0.935) node[shape=diamond, draw, fill=ceared, inner sep=0pt, minimum size=5pt] {};
\filldraw[black, line width=0.5pt] (3.5, 0.90) node[shape=diamond, draw, fill=ceared, inner sep=0pt, minimum size=5pt] {};
\filldraw[black, line width=0.5pt] (3.5, 0.80) node[shape=diamond, draw, fill=blue, inner sep=0pt, minimum size=5pt] {};
\filldraw[black, line width=0.5pt] (3.5, 0.70) node[shape=diamond, draw, fill=ceared, inner sep=0pt, minimum size=5pt] {};
\filldraw[black, line width=0.5pt] (3.5, 0.665) node[shape=diamond, draw, fill=ceared, inner sep=0pt, minimum size=5pt] {};
\filldraw[black, line width=0.5pt] (3.5, 0.56) node[shape=diamond, draw, fill=blue, inner sep=0pt, minimum size=5pt] {};
\filldraw[black, line width=0.5pt] (3.5, 0.435) node[shape=diamond, draw, fill=ceared, inner sep=0pt, minimum size=5pt] {};
\filldraw[black, line width=0.5pt] (3.5, 0.40) node[shape=diamond, draw, fill=ceared, inner sep=0pt, minimum size=5pt] {};
\filldraw[black, line width=0.5pt] (3.5, 0.30) node[shape=diamond, draw, fill=blue, inner sep=0pt, minimum size=5pt] {};
\filldraw[black, line width=0.5pt] (3.5, 0.20) node[shape=diamond, draw, fill=ceared, inner sep=0pt, minimum size=5pt] {};
\filldraw[black, line width=0.5pt] (3.5, 0.165) node[shape=diamond, draw, fill=ceared, inner sep=0pt, minimum size=5pt] {};
\filldraw[black, line width=0.5pt] (3.5, 0.06) node[shape=diamond, draw, fill=blue, inner sep=0pt, minimum size=5pt] {};

\filldraw[black, line width=0.5pt] (2.5, 0.7) node[shape=diamond, draw, fill=ceared, inner sep=0pt, minimum size=5pt] {};
\filldraw[black, line width=0.5pt] (2.5, 0.665) node[shape=diamond, draw, fill=ceared, inner sep=0pt, minimum size=5pt] {};
\filldraw[black, line width=0.5pt] (2.5, 0.335) node[shape=diamond, draw, fill=ceared, inner sep=0pt, minimum size=5pt] {};
\filldraw[black, line width=0.5pt] (2.5, 0.3) node[shape=diamond, draw, fill=ceared, inner sep=0pt, minimum size=5pt] {};

\filldraw[black, line width=0.5pt] (3.2, 1) node[shape=diamond, draw, fill=ceared, inner sep=0pt, minimum size=5pt] {};
\filldraw[black, line width=0.5pt] (3.165, 1) node[shape=diamond, draw, fill=ceared, inner sep=0pt, minimum size=5pt] {};
\filldraw[black, line width=0.5pt] (2.835, 1) node[shape=diamond, draw, fill=ceared, inner sep=0pt, minimum size=5pt] {};
\filldraw[black, line width=0.5pt] (2.8, 1) node[shape=diamond, draw, fill=ceared, inner sep=0pt, minimum size=5pt] {};

\filldraw[black, line width=0.5pt] (3.2, 0) node[shape=diamond, draw, fill=ceared, inner sep=0pt, minimum size=5pt] {};
\filldraw[black, line width=0.5pt] (3.165, 0) node[shape=diamond, draw, fill=ceared, inner sep=0pt, minimum size=5pt] {};
\filldraw[black, line width=0.5pt] (2.835, 0) node[shape=diamond, draw, fill=ceared, inner sep=0pt, minimum size=5pt] {};
\filldraw[black, line width=0.5pt] (2.8, 0) node[shape=diamond, draw, fill=ceared, inner sep=0pt, minimum size=5pt] {};

\filldraw[black, line width=0.5pt] (3.67, 0) node[shape=diamond, draw, fill=blue, inner sep=0pt, minimum size=5pt] {};
\filldraw[black, line width=0.5pt] (3.84, 0) node[shape=diamond, draw, fill=blue, inner sep=0pt, minimum size=5pt] {};
\filldraw[black, line width=0.5pt] (3.67, 0.5) node[shape=diamond, draw, fill=blue, inner sep=0pt, minimum size=5pt] {};
\filldraw[black, line width=0.5pt] (3.84, 0.5) node[shape=diamond, draw, fill=blue, inner sep=0pt, minimum size=5pt] {};
\filldraw[black, line width=0.5pt] (3.67, 1) node[shape=diamond, draw, fill=blue, inner sep=0pt, minimum size=5pt] {};
\filldraw[black, line width=0.5pt] (3.84, 1) node[shape=diamond, draw, fill=blue, inner sep=0pt, minimum size=5pt] {};

\filldraw[black, line width=0.5pt] (4.17, 0) node[shape=diamond, draw, fill=blue, inner sep=0pt, minimum size=5pt] {};
\filldraw[black, line width=0.5pt] (4.34, 0) node[shape=diamond, draw, fill=blue, inner sep=0pt, minimum size=5pt] {};
\filldraw[black, line width=0.5pt] (4.17, 0.5) node[shape=diamond, draw, fill=blue, inner sep=0pt, minimum size=5pt] {};
\filldraw[black, line width=0.5pt] (4.34, 0.5) node[shape=diamond, draw, fill=blue, inner sep=0pt, minimum size=5pt] {};
\filldraw[black, line width=0.5pt] (4.17, 1) node[shape=diamond, draw, fill=blue, inner sep=0pt, minimum size=5pt] {};
\filldraw[black, line width=0.5pt] (4.34, 1) node[shape=diamond, draw, fill=blue, inner sep=0pt, minimum size=5pt] {};
    
\filldraw[black, line width=0.5pt] (4, 0.17) node[shape=diamond, draw, fill=blue, inner sep=0pt, minimum size=5pt] {};
\filldraw[black, line width=0.5pt] (4, 0.33) node[shape=diamond, draw, fill=blue, inner sep=0pt, minimum size=5pt] {};
\filldraw[black, line width=0.5pt] (4, 0.67) node[shape=diamond, draw, fill=blue, inner sep=0pt, minimum size=5pt] {};
\filldraw[black, line width=0.5pt] (4, 0.83) node[shape=diamond, draw, fill=blue, inner sep=0pt, minimum size=5pt] {};

\end{scope}

\begin{scope}[xshift=-4.75cm, yshift=0.2cm]
    \draw[black, line width=0.5pt, fill=ceared!25] (2.5,0) -- (3.5,0) -- (3.5,1) -- (2.5,1) -- cycle;
    \draw[black, line width=0.5pt, fill=blue!25] (3.5,0) -- (4.5,0) -- (4.5,1) -- (3.5,1) -- cycle;
    \draw[black, line width=0.5pt] (4.0, 0.00) -- (4.0, 1.00);
    \draw[black, line width=0.5pt] (3.5, 0.50) -- (4.5, 0.50);
    
    \filldraw[black, line width=0.5pt] (3.5, 0.935) node[shape=diamond, draw, fill=ceared, inner sep=0pt, minimum size=5pt] {};
    \filldraw[black, line width=0.5pt] (3.5, 0.90) node[shape=diamond, draw, fill=ceared, inner sep=0pt, minimum size=5pt] {};
    \filldraw[black, line width=0.5pt] (3.5, 0.80) node[shape=diamond, draw, fill=blue, inner sep=0pt, minimum size=5pt] {};
    \filldraw[black, line width=0.5pt] (3.5, 0.70) node[shape=diamond, draw, fill=ceared, inner sep=0pt, minimum size=5pt] {};
    \filldraw[black, line width=0.5pt] (3.5, 0.665) node[shape=diamond, draw, fill=ceared, inner sep=0pt, minimum size=5pt] {};
    \filldraw[black, line width=0.5pt] (3.5, 0.56) node[shape=diamond, draw, fill=blue, inner sep=0pt, minimum size=5pt] {};
    \filldraw[black, line width=0.5pt] (3.5, 0.435) node[shape=diamond, draw, fill=ceared, inner sep=0pt, minimum size=5pt] {};
    \filldraw[black, line width=0.5pt] (3.5, 0.40) node[shape=diamond, draw, fill=ceared, inner sep=0pt, minimum size=5pt] {};
    \filldraw[black, line width=0.5pt] (3.5, 0.30) node[shape=diamond, draw, fill=blue, inner sep=0pt, minimum size=5pt] {};
    \filldraw[black, line width=0.5pt] (3.5, 0.20) node[shape=diamond, draw, fill=ceared, inner sep=0pt, minimum size=5pt] {};
    \filldraw[black, line width=0.5pt] (3.5, 0.165) node[shape=diamond, draw, fill=ceared, inner sep=0pt, minimum size=5pt] {};
    \filldraw[black, line width=0.5pt] (3.5, 0.06) node[shape=diamond, draw, fill=blue, inner sep=0pt, minimum size=5pt] {};
    
    \filldraw[black, line width=0.5pt] (2.5, 0.7) node[shape=diamond, draw, fill=ceared, inner sep=0pt, minimum size=5pt] {};
    \filldraw[black, line width=0.5pt] (2.5, 0.665) node[shape=diamond, draw, fill=ceared, inner sep=0pt, minimum size=5pt] {};
    \filldraw[black, line width=0.5pt] (2.5, 0.335) node[shape=diamond, draw, fill=ceared, inner sep=0pt, minimum size=5pt] {};
    \filldraw[black, line width=0.5pt] (2.5, 0.3) node[shape=diamond, draw, fill=ceared, inner sep=0pt, minimum size=5pt] {};
    
    \filldraw[black, line width=0.5pt] (3.2, 1) node[shape=diamond, draw, fill=ceared, inner sep=0pt, minimum size=5pt] {};
    \filldraw[black, line width=0.5pt] (3.165, 1) node[shape=diamond, draw, fill=ceared, inner sep=0pt, minimum size=5pt] {};
    \filldraw[black, line width=0.5pt] (2.835, 1) node[shape=diamond, draw, fill=ceared, inner sep=0pt, minimum size=5pt] {};
    \filldraw[black, line width=0.5pt] (2.8, 1) node[shape=diamond, draw, fill=ceared, inner sep=0pt, minimum size=5pt] {};
    
    \filldraw[black, line width=0.5pt] (3.2, 0) node[shape=diamond, draw, fill=ceared, inner sep=0pt, minimum size=5pt] {};
    \filldraw[black, line width=0.5pt] (3.165, 0) node[shape=diamond, draw, fill=ceared, inner sep=0pt, minimum size=5pt] {};
    \filldraw[black, line width=0.5pt] (2.835, 0) node[shape=diamond, draw, fill=ceared, inner sep=0pt, minimum size=5pt] {};
    \filldraw[black, line width=0.5pt] (2.8, 0) node[shape=diamond, draw, fill=ceared, inner sep=0pt, minimum size=5pt] {};
    
    \filldraw[black, line width=0.5pt] (3.67, 0) node[shape=diamond, draw, fill=blue, inner sep=0pt, minimum size=5pt] {};
    \filldraw[black, line width=0.5pt] (3.84, 0) node[shape=diamond, draw, fill=blue, inner sep=0pt, minimum size=5pt] {};
    \filldraw[black, line width=0.5pt] (3.67, 0.5) node[shape=diamond, draw, fill=blue, inner sep=0pt, minimum size=5pt] {};
    \filldraw[black, line width=0.5pt] (3.84, 0.5) node[shape=diamond, draw, fill=blue, inner sep=0pt, minimum size=5pt] {};
    \filldraw[black, line width=0.5pt] (3.67, 1) node[shape=diamond, draw, fill=blue, inner sep=0pt, minimum size=5pt] {};
    \filldraw[black, line width=0.5pt] (3.84, 1) node[shape=diamond, draw, fill=blue, inner sep=0pt, minimum size=5pt] {};
    
    \filldraw[black, line width=0.5pt] (4.17, 0) node[shape=diamond, draw, fill=blue, inner sep=0pt, minimum size=5pt] {};
    \filldraw[black, line width=0.5pt] (4.34, 0) node[shape=diamond, draw, fill=blue, inner sep=0pt, minimum size=5pt] {};
    \filldraw[black, line width=0.5pt] (4.17, 0.5) node[shape=diamond, draw, fill=blue, inner sep=0pt, minimum size=5pt] {};
    \filldraw[black, line width=0.5pt] (4.34, 0.5) node[shape=diamond, draw, fill=blue, inner sep=0pt, minimum size=5pt] {};
    \filldraw[black, line width=0.5pt] (4.17, 1) node[shape=diamond, draw, fill=blue, inner sep=0pt, minimum size=5pt] {};
    \filldraw[black, line width=0.5pt] (4.34, 1) node[shape=diamond, draw, fill=blue, inner sep=0pt, minimum size=5pt] {};
        
    \filldraw[black, line width=0.5pt] (4, 0.17) node[shape=diamond, draw, fill=blue, inner sep=0pt, minimum size=5pt] {};
    \filldraw[black, line width=0.5pt] (4, 0.33) node[shape=diamond, draw, fill=blue, inner sep=0pt, minimum size=5pt] {};
    \filldraw[black, line width=0.5pt] (4, 0.67) node[shape=diamond, draw, fill=blue, inner sep=0pt, minimum size=5pt] {};
    \filldraw[black, line width=0.5pt] (4, 0.83) node[shape=diamond, draw, fill=blue, inner sep=0pt, minimum size=5pt] {};
        
    \draw[->, thick, bend left=15] (3.75,1.15) to (5.25,1.25);
    \draw[<-, thick, bend left=15] (1.75,1.25) to (3.25,1.15);
    \begin{scope}[xshift=2.525cm, yshift=-0.75cm]
    \draw[black, line width=0.5pt, fill=ceared] (0.20,1.4) circle (1.25pt);    
    \draw[black, line width=0.5pt, fill=ceared] (0.29,1.4) circle (1.25pt);   
    \draw[black, line width=0.5pt, fill=ceared] (0.65,1.4) circle (1.25pt);    
    \draw[black, line width=0.5pt, fill=ceared] (0.74,1.4) circle (1.25pt);   
    \draw[black, line width=0.5pt, fill=ceared] (0.42,1.6) circle (1.25pt);
    \draw[black, line width=0.5pt, fill=ceared] (0.51,1.6) circle (1.25pt);  
    \begin{scope}[yshift=1.125cm]
    \draw[black, line width=0.5pt, fill=ceared] (0.20,-0.150) circle (1.25pt);    
    \draw[black, line width=0.5pt, fill=ceared] (0.29,-0.150) circle (1.25pt);    
    \draw[black, line width=0.5pt, fill=ceared] (0.29,-0.245) circle (1.25pt); 
    \draw[black, line width=0.5pt, fill=ceared] (0.42, 0.075) circle (1.25pt);    
    \draw[black, line width=0.5pt, fill=ceared] (0.51, 0.075) circle (1.25pt);    
    \draw[black, line width=0.5pt, fill=ceared] (0.51, -0.02) circle (1.25pt); 
    \draw[black, line width=0.5pt, fill=ceared] (0.65,-0.150) circle (1.25pt);    
    \draw[black, line width=0.5pt, fill=ceared] (0.74,-0.150) circle (1.25pt);    
    \draw[black, line width=0.5pt, fill=ceared] (0.74,-0.245) circle (1.25pt);  
\end{scope}
\begin{scope}[scale=0.75, xshift=0.3875cm, yshift=0.65cm]
    \draw[black, line width=0.5pt, fill=blue] (1.15, 1.45) circle (1.25pt);    
    \draw[black, line width=0.5pt, fill=blue] (1.35, 1.45) circle (1.25pt);    
    \draw[black, line width=0.5pt, fill=blue] (1.25, 1.60) circle (1.25pt);  

    \draw[black, line width=0.5pt, fill=blue] (1.825, 1.45) circle (1.25pt);    
    \draw[black, line width=0.5pt, fill=blue] (2.025, 1.45) circle (1.25pt);    
    \draw[black, line width=0.5pt, fill=blue] (1.925, 1.60) circle (1.25pt);

    \draw[black, line width=0.5pt, fill=blue] (1.825, 0.78) circle (1.25pt);    
    \draw[black, line width=0.5pt, fill=blue] (2.025, 0.78) circle (1.25pt);    
    \draw[black, line width=0.5pt, fill=blue] (1.925, 0.93) circle (1.25pt);

    \draw[black, line width=0.5pt, fill=blue] (1.15, 0.78) circle (1.25pt);    
    \draw[black, line width=0.5pt, fill=blue] (1.35, 0.78) circle (1.25pt);    
    \draw[black, line width=0.5pt, fill=blue] (1.25, 0.93) circle (1.25pt);      
\end{scope}
\begin{scope}[scale=0.75, xshift=0.3875cm, yshift=1.75cm]
    \draw[black, line width=0.5pt, fill=blue] (1.59+0.175, 0.05) circle (1.25pt);    
    \draw[black, line width=0.5pt, fill=blue] (1.68+0.175, 0.05) circle (1.25pt);    
    \draw[black, line width=0.5pt, fill=blue] (1.82+0.175, 0.05) circle (1.25pt);    
    \draw[black, line width=0.5pt, fill=blue] (1.91+0.175, 0.05) circle (1.25pt);    
    \draw[black, line width=0.5pt, fill=blue] (1.70+0.175, 0.20) circle (1.25pt);
    \draw[black, line width=0.5pt, fill=blue] (1.79+0.175, 0.20) circle (1.25pt);

    \draw[black, line width=0.5pt, fill=blue] (1.09,0.05) circle (1.25pt);    
    \draw[black, line width=0.5pt, fill=blue] (1.18,0.05) circle (1.25pt);    
    \draw[black, line width=0.5pt, fill=blue] (1.32,0.05) circle (1.25pt);    
    \draw[black, line width=0.5pt, fill=blue] (1.41,0.05) circle (1.25pt);    
    \draw[black, line width=0.5pt, fill=blue] (1.20,0.20) circle (1.25pt);  
    \draw[black, line width=0.5pt, fill=blue] (1.29,0.20) circle (1.25pt);  
    
    \draw[black, line width=0.5pt, fill=blue] (1.59+0.175,-0.60) circle (1.25pt);    
    \draw[black, line width=0.5pt, fill=blue] (1.68+0.175,-0.60) circle (1.25pt);    
    \draw[black, line width=0.5pt, fill=blue] (1.82+0.175,-0.60) circle (1.25pt);    
    \draw[black, line width=0.5pt, fill=blue] (1.91+0.175,-0.60) circle (1.25pt);    
    \draw[black, line width=0.5pt, fill=blue] (1.70+0.175,-0.45) circle (1.25pt);
    \draw[black, line width=0.5pt, fill=blue] (1.79+0.175,-0.45) circle (1.25pt);

    \draw[black, line width=0.5pt, fill=blue] (1.09,-0.60) circle (1.25pt);    
    \draw[black, line width=0.5pt, fill=blue] (1.18,-0.60) circle (1.25pt);    
    \draw[black, line width=0.5pt, fill=blue] (1.32,-0.60) circle (1.25pt);    
    \draw[black, line width=0.5pt, fill=blue] (1.41,-0.60) circle (1.25pt);    
    \draw[black, line width=0.5pt, fill=blue] (1.20,-0.45) circle (1.25pt);  
    \draw[black, line width=0.5pt, fill=blue] (1.29,-0.45) circle (1.25pt); 
\end{scope}

\end{scope}

\end{scope}

\begin{scope}[xshift=2.25cm]
    \begin{scope}[xshift=-4.75cm, yshift=0.2cm]
        \draw[black, line width=0.5pt, fill=ceared!25] (2.5,0) -- (3.5,0) -- (3.5,1) -- (2.5,1) -- cycle;
        \draw[black, line width=0.5pt, fill=blue!25] (3.5,0) -- (4.5,0) -- (4.5,1) -- (3.5,1) -- cycle;
        \draw[black, line width=0.5pt] (4.0, 0.00) -- (4.0, 1.00);
        \draw[black, line width=0.5pt] (3.5, 0.50) -- (4.5, 0.50);
                    
        \begin{scope}[xshift=2.525cm, yshift=-0.75cm]
        \draw[black, line width=0.5pt, fill=ceared] (0.20,1.4) circle (1.25pt);    
        \draw[black, line width=0.5pt, fill=ceared] (0.29,1.4) circle (1.25pt);   
        \draw[black, line width=0.5pt, fill=ceared] (0.65,1.4) circle (1.25pt);    
        \draw[black, line width=0.5pt, fill=ceared] (0.74,1.4) circle (1.25pt);   
        \draw[black, line width=0.5pt, fill=ceared] (0.42,1.6) circle (1.25pt);
        \draw[black, line width=0.5pt, fill=ceared] (0.51,1.6) circle (1.25pt);  
        \begin{scope}[yshift=1.125cm]
        \draw[black, line width=0.5pt, fill=ceared] (0.20,-0.150) circle (1.25pt);    
        \draw[black, line width=0.5pt, fill=ceared] (0.29,-0.150) circle (1.25pt);    
        \draw[black, line width=0.5pt, fill=ceared] (0.29,-0.245) circle (1.25pt); 
        \draw[black, line width=0.5pt, fill=ceared] (0.42, 0.075) circle (1.25pt);    
        \draw[black, line width=0.5pt, fill=ceared] (0.51, 0.075) circle (1.25pt);    
        \draw[black, line width=0.5pt, fill=ceared] (0.51, -0.02) circle (1.25pt); 
        \draw[black, line width=0.5pt, fill=ceared] (0.65,-0.150) circle (1.25pt);    
        \draw[black, line width=0.5pt, fill=ceared] (0.74,-0.150) circle (1.25pt);    
        \draw[black, line width=0.5pt, fill=ceared] (0.74,-0.245) circle (1.25pt);  
    \end{scope}
    \begin{scope}[scale=0.75, xshift=0.3875cm, yshift=0.65cm]
        \draw[black, line width=0.5pt, fill=blue] (1.15, 1.45) circle (1.25pt);    
        \draw[black, line width=0.5pt, fill=blue] (1.35, 1.45) circle (1.25pt);    
        \draw[black, line width=0.5pt, fill=blue] (1.25, 1.60) circle (1.25pt);  
    
        \draw[black, line width=0.5pt, fill=blue] (1.825, 1.45) circle (1.25pt);    
        \draw[black, line width=0.5pt, fill=blue] (2.025, 1.45) circle (1.25pt);    
        \draw[black, line width=0.5pt, fill=blue] (1.925, 1.60) circle (1.25pt);
    
        \draw[black, line width=0.5pt, fill=blue] (1.825, 0.78) circle (1.25pt);    
        \draw[black, line width=0.5pt, fill=blue] (2.025, 0.78) circle (1.25pt);    
        \draw[black, line width=0.5pt, fill=blue] (1.925, 0.93) circle (1.25pt);
    
        \draw[black, line width=0.5pt, fill=blue] (1.15, 0.78) circle (1.25pt);    
        \draw[black, line width=0.5pt, fill=blue] (1.35, 0.78) circle (1.25pt);    
        \draw[black, line width=0.5pt, fill=blue] (1.25, 0.93) circle (1.25pt);      
    \end{scope}
    \begin{scope}[scale=0.75, xshift=0.3875cm, yshift=1.75cm]
        \draw[black, line width=0.5pt, fill=blue] (1.59+0.175, 0.05) circle (1.25pt);    
        \draw[black, line width=0.5pt, fill=blue] (1.68+0.175, 0.05) circle (1.25pt);    
        \draw[black, line width=0.5pt, fill=blue] (1.82+0.175, 0.05) circle (1.25pt);    
        \draw[black, line width=0.5pt, fill=blue] (1.91+0.175, 0.05) circle (1.25pt);    
        \draw[black, line width=0.5pt, fill=blue] (1.70+0.175, 0.20) circle (1.25pt);
        \draw[black, line width=0.5pt, fill=blue] (1.79+0.175, 0.20) circle (1.25pt);
    
        \draw[black, line width=0.5pt, fill=blue] (1.09,0.05) circle (1.25pt);    
        \draw[black, line width=0.5pt, fill=blue] (1.18,0.05) circle (1.25pt);    
        \draw[black, line width=0.5pt, fill=blue] (1.32,0.05) circle (1.25pt);    
        \draw[black, line width=0.5pt, fill=blue] (1.41,0.05) circle (1.25pt);    
        \draw[black, line width=0.5pt, fill=blue] (1.20,0.20) circle (1.25pt);  
        \draw[black, line width=0.5pt, fill=blue] (1.29,0.20) circle (1.25pt);  
        
        \draw[black, line width=0.5pt, fill=blue] (1.59+0.175,-0.60) circle (1.25pt);    
        \draw[black, line width=0.5pt, fill=blue] (1.68+0.175,-0.60) circle (1.25pt);    
        \draw[black, line width=0.5pt, fill=blue] (1.82+0.175,-0.60) circle (1.25pt);    
        \draw[black, line width=0.5pt, fill=blue] (1.91+0.175,-0.60) circle (1.25pt);    
        \draw[black, line width=0.5pt, fill=blue] (1.70+0.175,-0.45) circle (1.25pt);
        \draw[black, line width=0.5pt, fill=blue] (1.79+0.175,-0.45) circle (1.25pt);
    
        \draw[black, line width=0.5pt, fill=blue] (1.09,-0.60) circle (1.25pt);    
        \draw[black, line width=0.5pt, fill=blue] (1.18,-0.60) circle (1.25pt);    
        \draw[black, line width=0.5pt, fill=blue] (1.32,-0.60) circle (1.25pt);    
        \draw[black, line width=0.5pt, fill=blue] (1.41,-0.60) circle (1.25pt);    
        \draw[black, line width=0.5pt, fill=blue] (1.20,-0.45) circle (1.25pt);  
        \draw[black, line width=0.5pt, fill=blue] (1.29,-0.45) circle (1.25pt); 
    \end{scope}
    
    \end{scope}
    
    \end{scope}
\end{scope}

\node at (-3.5,1.45)  {\color{black} \textsc{Implicit}};
\node at (1,1.45)  {\color{black} \textsc{Explicit}};

\begin{scope}[xshift=1.45cm, yshift=0.15cm]
    \node at (-3.60,-0.2)  {\color{blue} $p_{\Tf}$};
    \draw[black, line width=0.5pt, fill=blue] (-3.95, -0.2) circle (1.25pt);    
    
    \node at (-2.35,-0.2)  {\color{blue} $p_{\Ff}$};
    \filldraw[black, line width=0.5pt] (-2.7, -0.2) node[shape=diamond, draw, fill=blue, inner sep=0pt, minimum size=5pt] {};
     
    \node at (-1.10,-0.2)  {\color{blue} $\bd{v}_{\Tf}$};
    \draw[black, line width=0.5pt, fill=blue] (-1.40,-0.2) circle (1.25pt);    
    \draw[black, line width=0.5pt, fill=blue] (-1.49,-0.2) circle (1.25pt); 
    
    \node at (-3.60,-0.65+0.1)  {\color{ceared} $\bd{v}_{\Ts}$};
    \draw[black, line width=0.5pt, fill=ceared] (-3.90,-0.65+0.1) circle (1.25pt);    
    \draw[black, line width=0.5pt, fill=ceared] (-3.99,-0.65+0.1) circle (1.25pt); 
    
    \node at (-2.35,-0.65+0.1)  {\color{ceared} $\bd{v}_{\Fs}$};
    \filldraw[black, line width=0.5pt] (-2.675, -0.65+0.1) node[shape=diamond, draw, fill=ceared, inner sep=0pt, minimum size=6pt] {};
    \filldraw[black, line width=0.5pt] (-2.725, -0.65+0.1) node[shape=diamond, draw, fill=ceared, inner sep=0pt, minimum size=6pt] {};
        
    \node at (-1.10,-0.65+0.1)  {\color{ceared} $\bbm{s}_{\Ts}$};
    \draw[black, line width=0.5pt, fill=ceared] (-1.40,-0.60+0.1) circle (1.25pt);    
    \draw[black, line width=0.5pt, fill=ceared] (-1.49,-0.60+0.1) circle (1.25pt);   
    \draw[black, line width=0.5pt, fill=ceared] (-1.40,-0.69+0.1) circle (1.25pt); 
\end{scope}

\end{tikzpicture}
}
\caption{Fully coupled unknowns after static condensation for SDIRK schemes (left) and for ERK schemes (right) for the lowest equal-order discretization $(k^\prime=k=1)$.}
\label{fig::sc}
\end{figure}
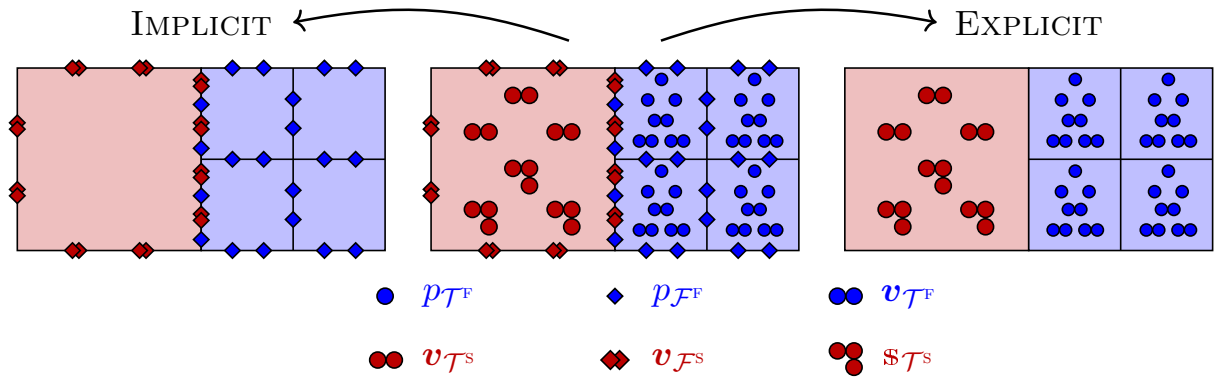

We observe from \hyperref[fig::sc]{\Cref{fig::sc}} that the number of cell dofs appears to be larger than the face dofs. Let us give a more quantitative estimate. Neglecting the number of boundary faces, and assuming that each cell has $n$ faces (e.g., $n=3$ for triangles and $n=4$ for tetrahedra), the asymptotic ratio between the number of cells and the number of faces is $n \cdot \#\sc{Cells} \approx 2 \#\sc{Faces}$. Under this assumption, \hyperref[tab::DOFS]{\Cref{tab::DOFS}} reports the expected percentage of dofs associated with faces for triangular and tetrahedral meshes, and for equal- and mixed-order HHO settings. A first observation is, as expected, that the proportion of face dofs decreases with the polynomial degree. A second observation is that, in the mixed-order case (used with implicit time schemes), the share of cell dofs is at least 75\%. In this setting, static condensation — which eliminates here the cell dofs — allows for a significant reduction in the size of the global linear system. In the equal-order case (used with explicit time schemes), the proportion of cell dofs is smaller but remains dominant (at least 67\%). Here, static condensation eliminates face dofs instead, which still reduces the system size, albeit less substantially. Altogether, we can conclude that the reduction in computational costs by static condensation is much more impactful when using implicit time schemes.
\ifHAL
\begin{table}[!htb]
\centering
\resizebox{0.875\textwidth}{!}{%
\begin{tabular}{c|cc|c|c|c|c|}
\cline{4-7} 
\multicolumn{3}{c|}{} & \multicolumn{1}{c|}{\rule{0pt}{2.75ex} \sc{General $k$}} & \multicolumn{1}{c|}{$k=1$} & \multicolumn{1}{c|}{$k=2$} & \multicolumn{1}{c|}{$k=3$} \\ 
\hline 
\multicolumn{1}{|c}{\multirow{5}{*}{\sc{Equal-order}}} & \multicolumn{1}{|c}{\multirow{2}{*}{$2$D}} & \sc{Acoustic} & \rule{0pt}{3.5ex} $\frac{1}{k+3}$ & $25\%$ & $20\%$ & $17\%$ \\[1ex]
\multicolumn{1}{|c}{} & \multicolumn{1}{|c}{} & \sc{Elastic} & $\frac{6}{5k+16}$ & $29\%$ & $23\%$ & $19\%$ \\[1ex]
\cline{2-7}
\multicolumn{1}{|c}{\multirow{2}{*}{}} & \multicolumn{1}{|c}{\multirow{2}{*}{$3$D}} & \sc{Acoustic} & \rule{0pt}{3.5ex} $\frac{3}{2k+9}$ & $27\%$ & $23\%$ & $20\%$ \\
\multicolumn{1}{|c}{} & \multicolumn{1}{|c}{} & \sc{Elastic} & \rule{0pt}{3.5ex} $\frac{2}{k+5}$ & $33\%$ & $29\%$ & $25\%$ \\[1ex]
\hline   
\multicolumn{1}{|c}{\multirow{5}{*}{\sc{Mixed-order}}} & \multicolumn{1}{|c}{\multirow{2}{*}{$2$D}} & \sc{Acoustic} & \rule{0pt}{3.5ex} $\frac{3(k+1)}{3k^2+14k+13}$ & $20\%$ & $17\%$ & $15\%$\\[0.5ex]
\multicolumn{1}{|c}{} & \multicolumn{1}{|c}{} & \sc{Elastic} & \rule{0pt}{3.5ex} $\frac{6(k+1)}{5k^2+25k+24}$ & $22\%$ & $19\%$ & $17\%$ \\[1ex]
\cline{2-7}
\multicolumn{1}{|c}{\multirow{2}{*}{}} & \multicolumn{1}{|c}{\multirow{2}{*}{$3$D}} & \sc{Acoustic} & \rule{0pt}{3.5ex} $\frac{6(k+1)(k+2)}{4k^3+33k^2+77k+54}$ & $21\%$ & $19\%$ & $17\%$ \\
\multicolumn{1}{|c}{} & \multicolumn{1}{|c}{} & \sc{Elastic} & \rule{0pt}{3.5ex} $\frac{6(k+1)(k+2)}{3k^3+27k^2+66k+48}$ & $25\%$ & $23\%$ & $21\%$ \\[1ex]
\hline  
\end{tabular}
}
\caption{Percentage of face dofs for purely acoustic and elastic cases, in $2$D and $3$D, for equal- and mixed-order settings, for various polynomial degrees $k$.}
\label{tab::DOFS}
\end{table}
\fi

\subsection{Singly Diagonal Implicit Runge-Kutta (SDIRK) schemes}
   
Given $\dofs{U}{T}{}^{n-1}$ from the previous time step or the initial condition, the time discretization of \eqref{compact_AR} by SDIRK$(s,s+1)$ schemes is as follows: We solve sequentially for all $i \in \{1,...,s+1\}$,
\begin{equation}
\resizebox{0.925\textwidth}{!}{%
$
\pmatset{4}{12pt} 
\pmatset{5}{12pt}
\begin{pmat}[{}]
\normalizecell{\mass{}{}} & 0 \cr
0 & 0  \cr 
\end{pmat}
\begin{pmat}[{}]
\rm{U}_{\T}^{n,i} \cr
\rm{U}_{\F}^{n,i}  \cr 
\end{pmat} = 
\pmatset{4}{12pt} 
\pmatset{5}{12pt}
\begin{pmat}[{}]
\normalizecell{\mass{}{}} & 0 \cr
0 & 0  \cr 
\end{pmat}
\begin{pmat}[{}]
\rm{U}_{\T}^{n-1} \cr
\rm{U}_{\F}^{n-1}  \cr 
\end{pmat} +
\Delta t {\displaystyle \sum_{j=1}^i a_{ij}}
\left( 
\begin{pmat}[{}]
\rm{F}_{\T}^{n-1+c_j} \cr
0  \cr 
\end{pmat} -
\pmatset{4}{12pt} 
\pmatset{5}{12pt}
\begin{pmat}[{}]
\normalizecell{\cal{K}_{\cal{T}}} & \normalizecell{\cal{K}_{\cal{TF}}} \cr
\normalizecell{\cal{K}_{\cal{FT}}} & \normalizecell{\cal{K}_{\cal{F}}}  \cr 
\end{pmat}
\begin{pmat}[{}]
\rm{U}_{\T}^{n,j} \cr
\rm{U}_{\F}^{n,j}  \cr 
\end{pmat} 
\right),
\label{SDIRK_compact}
$
}
\end{equation}
where $\dofs{F}{T}{}^{n-1+c_j} := \dofs{F}{T}{}(t_{n-1}+c_j \Delta t)$. Notice that the last stage $(i=s+1)$ only requires a mass matrix inversion. Finally, we set $\dofs{U}{T}{}^{n+1} := \dofs{U}{T}{}^{n,s+1}$. For all $i \in \{1,...,s\}$, each stage in $\eqref{SDIRK_compact}$ can be rewritten as
\begin{equation}
\pmatset{4}{12pt} 
\pmatset{5}{12pt}
\begin{pmat}[{}]
\mass{}{} + a_{*} \Delta t \cal{K}_{\T}   & a_{*} \Delta t \cal{K}_{\T\F} \cr
a_{*} \Delta t \cal{K}_{\F\T} & a_{*} \Delta t \cal{K}_{\F} \cr 
\end{pmat}
\label{cp_dirk}
\begin{pmat}[{}]
\rm{U}_{\T}^{n,i} \cr
\rm{U}_{\F}^{n,i}  \cr 
\end{pmat}
=
\pmatset{4}{12pt} 
\pmatset{5}{12pt}
\begin{pmat}[{}]
\rm{B}_{\T}^{n,i} \cr
\rm{B}_{\F}^{n,i}  \cr 
\end{pmat},
\end{equation} 
with
\begin{equation}
\pmatset{4}{12pt} 
\pmatset{5}{12pt}
\begin{pmat}[{}]
\rm{B}_{\T}^{n,i} \cr
\rm{B}_{\F}^{n,i}  \cr 
\end{pmat} := 
\pmatset{4}{12pt} 
\pmatset{5}{12pt}
\begin{pmat}[{}]
\cal{M}_{\T}\rm{U}_{\T}^{n-1} + a_* \Delta t \rm{F}_{\T}^{n-1+c_i} + \Delta t {\displaystyle \sum_{j=1}^{i-1} a_{ij}} \left(\rm{F}_{\T}^{n-1+c_j} -\cal{K}_{\T}\rm{U}_{\T}^{n,j} -\cal{K}_{\T\F}\rm{U}_{\F}^{n,j}\right) \cr
-\Delta t {\displaystyle \sum_{j=1}^{i-1} a_{ij}} \left(\cal{K}_{\F\T}\rm{U}_{\T}^{n,j} + \cal{K}_{\F}\rm{U}_{\F}^{n,j}\right)\cr 
\end{pmat}.
\end{equation} 
A key observation is that the cell-cell submatrix in \eqref{cp_dirk} (associated with all the fluid and solid cell unknowns) is block-diagonal. Hence, a static condensation procedure can be performed, leading to a global problem coupling only the face unknowns. This strategy enhances computational efficiency by eliminating locally all the cell unknowns. Therefore, \eqref{cp_dirk} is solved by performing sequentially the following two solves:
\begin{subequations}
\begin{align}
& 
\begin{aligned}
a_* \Delta t\Big(\cal{K}_{\F} - a_* \Delta t\cal{K}_{\F\T} \left(\mass{}{} + a_* \Delta t \cal{K}_{\T}\right)^{-1} & \cal{K}_{\T\F}\Big)\rm{U}_{\F}^{n,i} =\\ 
& \rm{B}_\F^{n,i} - a_* \Delta t\cal{K}_{\F\T}\left(\mass{}{} + a_* \Delta t \cal{K}_{\T}\right)^{-1}\rm{B}_\T^{n,i},
\end{aligned}\\
& \rm{U}_{\T}^{n,i} = \left(\mass{}{} + a_* \Delta t \cal{K}_{\T}\right)^{-1} \left(\rm{B}_\T^{n,i} - a_* \Delta t \cal{K}_{\T\F}\rm{U}_{\F}^{n,i} \right).
\label{eq:sc2}
\end{align}
\end{subequations} 
Owing to the singly diagonal structure of the Butcher tableau, all the stages share the same system matrix allowing to re-use matrix factorizations. Moreover, for a constant time step, matrix factorizations can be re-used throughout the whole time integration. Notice also that the second step \eqref{eq:sc2} is local and embarrassingly parallel.

In this work, we employ $s$-stage Singly Diagonally Implicit Runge--Kutta (SDIRK) methods of order $(s{+}1)$, with $s \in \{2,3\}$. The associated Butcher tableaux are reported in \eqref{SDIRK} where the SDIRK$(3,4)$ scheme corresponds to the values $\nu := \frac{1}{\sqrt{3}} \cos\left(\frac{\pi}{18}\right) + \frac{1}{2}$ and $\xi := \frac{1}{6(2\nu - 1)^2}$:
\begin{equation}
\renewcommand{\arraystretch}{1.25} 
\begin{tabular}{c|cc}
$\frac{1}{4}$ & $\frac{1}{4}$ & $0$   \\ 
$\frac{3}{4}$ & $\frac{1}{2}$ & $\frac{1}{4}$  \\[0.1cm] 
\hline        & $\frac{1}{2}$ & $\frac{1}{2} $
\end{tabular}
\hspace{2cm}
\renewcommand{\arraystretch}{1.25} 
\begin{tabular}{c|ccc}
$\nu$      & $\nu$             & $0$         & $0$        \\ 
$\frac{1}{2}$ & $\frac{1}{2}-\nu$ & $\nu$    & $0$        \\ 
$1-\nu$    & $2\nu$            & $1-4\nu$ & $\nu$ \\ 
\hline        & $\xi$                & $1-2\xi$ & $\xi$
\end{tabular}
\label{SDIRK}
\end{equation}

\subsection{Explicit Runge--Kutta (ERK) schemes}

Given $\dofs{U}{T}{}^{n-1}$ from the previous time step or the initial condition, the time discretization of \eqref{compact_AR} by ERK$(s)$ schemes is as follows: We solve sequentially for all $i \in \{1,...,s+1\}$,
\begin{equation}
\resizebox{0.925\textwidth}{!}{%
$
\pmatset{4}{12pt} 
\pmatset{5}{12pt}
\begin{pmat}[{}]
\normalizecell{\mass{}{}} & 0 \cr
0 & 0  \cr 
\end{pmat}
\begin{pmat}[{}]
\rm{U}_{\T}^{n,i} \cr
\rm{U}_{\F}^{n,i}  \cr 
\end{pmat} = 
\pmatset{4}{12pt} 
\pmatset{5}{12pt}
\begin{pmat}[{}]
\normalizecell{\mass{}{}} & 0 \cr
0 & 0  \cr 
\end{pmat}
\begin{pmat}[{}]
\rm{U}_{\T}^{n-1} \cr
\rm{U}_{\F}^{n-1}  \cr 
\end{pmat} 
+ \Delta t {\displaystyle \sum_{j=1}^{i-1} a_{ij}}
\left( 
\begin{pmat}[{}]
\rm{F}_{\T}^{n-1+c_j} \cr
0  \cr 
\end{pmat} -
\pmatset{4}{12pt} 
\pmatset{5}{12pt}
\begin{pmat}[{}]
\normalizecell{\cal{K}_{\cal{T}}} & \normalizecell{\cal{K}_{\cal{TF}}} \cr
\normalizecell{\cal{K}_{\cal{FT}}} & \normalizecell{\cal{K}_{\cal{F}}}  
\cr 
\end{pmat}
\begin{pmat}[{}]
\rm{U}_{\T}^{n,j} \cr
\rm{U}_{\F}^{n,j}  \cr 
\end{pmat} 
\right),
$
}
\label{ERK_compact} 
\end{equation}
recalling that $\dofs{F}{T}{}^{n-1+c_j} := \dofs{F}{T}{}(t_{n-1}+c_j \Delta t)$. Finally, we set $\dofs{U}{T}{}^{n+1} := \dofs{U}{T}{}^{n,s+1}$. This time, we exploit the fact that the face-face submatrix $\cal{K}_\F$ is block-diagonal and we eliminate locally all the face unknowns. Thus, we rewrite \eqref{ERK_compact} only in terms of the cell unknowns (this procedure is called dG rewriting in the context of HDG methods). We first set $\dofs{U}{T}{}^{n,1} = \dofs{U}{T}{}^{n-1}$. Then, for all $i \in \{2,...,s+1\}$, we first solve:  
\begin{subequations}
\begin{equation}
\cal{K}_{\F}\rm{U}_{\F}^{n,i-1} = -\cal{K}_{\F\T}\rm{U}_{\T}^{n,i-1},
\end{equation}
followed by
\begin{equation}
\mass{}{} \rm{U}_{\T}^{n,i} = 
\mass{}{} \rm{U}_{\T}^{n-1} + \Delta t {\displaystyle \sum_{j=1}^{i-1} a_{ij}} \left( \rm{F}_{\T}^{n-1+c_j} - \left( \cal{K}_\T - \cal{K}_{\T\F}\cal{K}_{\F}^{-1}\cal{K}_{\F\T}\right) \dofs{U}{\T}{}^{n,j}\right).
\end{equation}
\label{eq:ex_solve}
\end{subequations}
We emphasize that both solves in \eqref{eq:ex_solve} are block-diagonal.

In this work, we use $s$-stage Explicit Runge--Kutta (ERK) schemes of order $s$, with $s \in \{2,3,4\}$. The corresponding Butcher tableaux are reported in \eqref{ERK}. We recall that ERK(s) schemes are subject to a CFL stability condition. \rev{From a theoretical viewpoint, stability under a CFL condition is established in \cite{EK_2024} for $s=3$ (and $s=4$ on simplicial meshes), and for $s=2$ under a strengthened 4/3-CFL condition; stability is also expected to hold under a standard CFL condition for $s=4$.} This point is further discussed numerically in Section 5.1 below.
\begin{equation}
\renewcommand{\arraystretch}{1.25} 
\begin{tabular}{c|cc}
$0$           & $0$           & $0$     \\ 
$\frac{1}{2}$ & $\frac{1}{2}$ & $0$ \\ 
\hline        & $0$ & $1$
\end{tabular}
\hspace{2cm}
\renewcommand{\arraystretch}{1.25} 
\begin{tabular}{c|ccc}
$0$           & $0$           & $0$           & $0$ \\ 
$\frac{1}{2}$ & $\frac{1}{2}$ & $0$           & $0$ \\ 
$1$           & $-1$          & $2$           & $0$ \\ 
\hline        & $\frac{1}{6}$ & $\frac{2}{3}$ & $\frac{1}{6}$
\end{tabular}
\hspace{2cm}
\renewcommand{\arraystretch}{1.25} 
\begin{tabular}{c|cccc}
$0$           & $0$           & $0$           & $0$           & $0$\\
$\frac{1}{2}$ & $\frac{1}{2}$ & $0$           & $0$           & $0$\\ 
$\frac{1}{2}$ & $0$           & $\frac{1}{2}$ & $0$           & $0$\\ 
$1$           & $0$           & $0$           & $1$           & $0$ \\ 
\hline        & $\frac{1}{6}$ & $\frac{1}{3}$ & $\frac{1}{3}$ & $\frac{1}{6}$
\end{tabular}
\label{ERK}
\end{equation}

\section{Numerical results}
 
In this section, we present 2D numerical results obtained with the HHO discretization of the elasto-acoustic problem described above. Our goal is threefold. Our first goal is to analyze the CFL condition of the explicit schemes and to perform a computational efficiency study comparing explicit and implicit schemes, in order to identify the most effective time-stepping strategy. To this purpose, we consider manufactured solutions. Our second goal is to assess the accuracy of the method using a classical test case with a Ricker wavelet as the initial condition, comparing the numerical results to a semi-analytical solution. Moreover, we provide a further comparison between explicit and implicit schemes using this second test case. Our third goal is to showcase the geometric flexibility of the method. To this purpose, we consider a geophysical application involving the propagation of a seismic wavefield through complex media. 

Recall that, for explicit time-stepping, we use an equal-order HHO discretization with $O(1)$ Least-Squares stabilization whereas, for implicit time-stepping, we use a mixed-order HHO discretization with Lehrenfeld-Schöberl $O(\frac{1}{h})$-stabilization. To solve the resulting linear systems, both a direct LU solver and an iterative Biconjugate Gradient (BiCG) solver preconditioned by an incomplete LU (ILU) factorization are used. The latter is considered so as to enable a fair comparison with the explicit schemes.

The implementation is carried out in the open-source software \texttt{disk++}, available at \url{https://github.com/wareHHOuse/diskpp}, and further detailed in \cite{CDE_2018}.

\subsection{CFL and efficiency study on a sinusoidal test case}

To define the discretization level we introduce two computational parameters: the spatial refinement level $\ell$, so that $h = 2^{-\ell}$, and the time refinement level $n$, so that $\Delta t = 0.1 \times 2^{-n}$. We consider three types of meshes: Cartesian, simplicial and polygonal. Some examples are shown in Figure \ref{meshes} for $\ell = 3$, with the solid subdomain mesh on the left side in red, and the fluid subdomain mesh on the right side in blue. 
\begin{figure}[!htb]
\centering 
\includegraphics[width=0.38\textwidth]{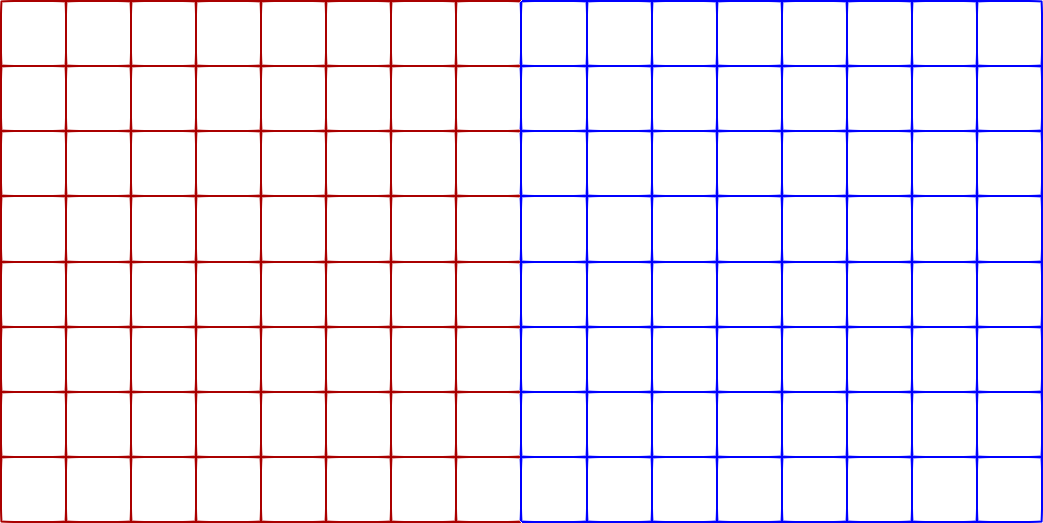}
\hspace{0.4cm}
\includegraphics[width=0.38\textwidth]{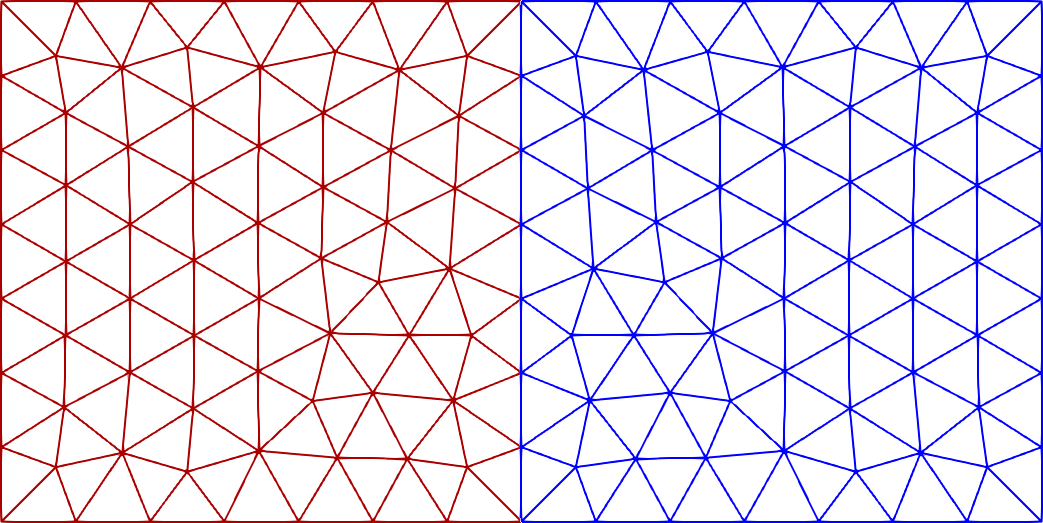}\\[0.6cm]
\includegraphics[width=0.38\textwidth]{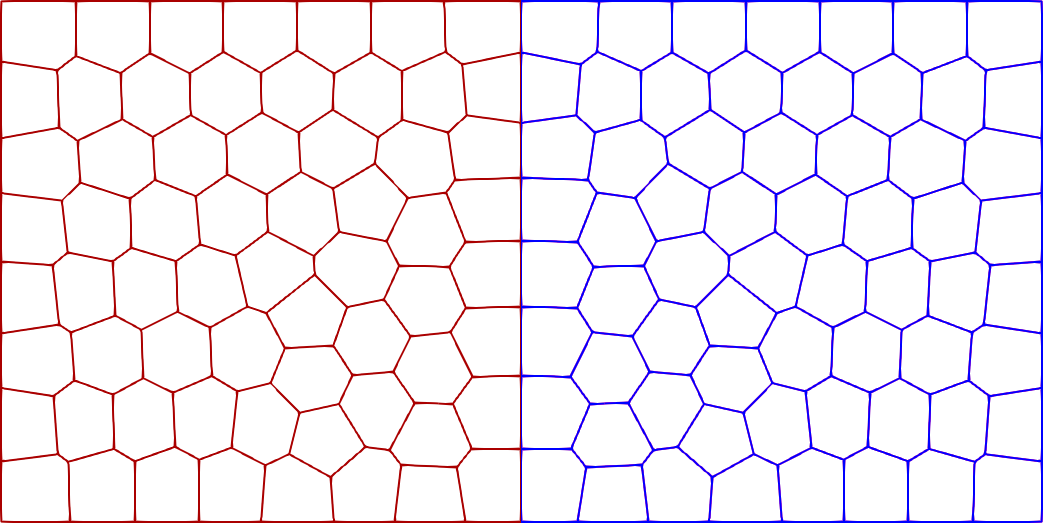}
\caption{Cartesian, simplicial, and polygonal meshes for $\ell = 3$.}
\label{meshes}
\end{figure}
\vspace{-0.15cm}

\paragraph{Test case setting.} The fluid and solid media have the same density and similar wave speeds. Specifically, we consider a simulation time $T_{\rm{f}} := 1$ and
\begin{itemize}
\item $\domain{f} := (0,1)\times(0,1)$, with density $\rho^\sc{f} := 1$, compressibility modulus $\kappa := 1$, so that the velocity of the pressure waves is $c_{\sc{p}}^\sc{f} = 1$;
\item $\domain{s} := (-1,0)\times(0,1)$, with density $\rho^\sc{s} := 1$, and Lamé parameters so that $c_\sc{p}^\sc{s} := \sqrt{3}$ and $c_\sc{s} = 1$.
\end{itemize}
The sinusoidal solution is expressed in terms of the potential $u$ (acoustic) and the displacement $\bd{u} := (u_x,u_y)$ (elastic) so that
\begin{subequations}
\begin{alignat}{3}
p &:= \partial_t u, &\qquad\qquad \bd{m} &:= \nabla u, & \qquad \qquad &\text{in } \domain{f},\\
\bd{v} & := \partial_t \bd{u}, & \qquad \qquad \bbm{C}^{-1}\bbm{s} & := \nabla_{\rm{sym}} \bd{u}, &\qquad\qquad & \text{in } \domain{s}.
\end{alignat}
\end{subequations}   
If not specified, the source terms, the (non)homogeneous Dirichlet boundary conditions, and the initial conditions are defined according to
\begin{subequations}
\begin{align}
u(t,x,y) & := x^2\sin(\omega \pi x)\sin(\omega \pi y)\sin(\theta \pi t), \\
u_x(t,x,y) := u_y(t,x,y) & = x^2\cos(\omega \tfrac{\pi}{2} x)\sin(\omega \pi y)\cos(\theta \pi t),
\end{align}
\label{non_poly_sol}
\end{subequations}
\hspace{-0.175cm}where $\omega$ and $\theta$ are the spatial and temporal frequencies, respectively. We consider two main configurations: a dominant spatial evolution with $\omega := 5$ and $\theta := \sqrt{2}$ and a dominant temporal evolution with $\omega := 1$ and $\theta := 10$.

\paragraph{Stability study of explicit time schemes.} We first study the stability of ERK$(s)$ schemes with $s \in \{2,3,4\}$. We consider additional weighting parameters, $\eta^\sc{f}$ and $\eta^\sc{s}$, which scale the stabilization bilinear forms appearing on the right-hand sides of \eqref{stab}. The stability condition for the ERK$(s)$ schemes can be written in the following form:
\begin{equation}
c_{\sharp} \frac{\Delta t}{h} \leq \sc{CFL$^*$}(s,k,\eta^{\sc{f}},\eta^{\sc{s}}),
\label{CFL}
\end{equation}
where \sc{CFL$^*$} denotes the critical value for stability and $c_{\sharp}$ is the largest velocity in the domain. To determine the value of CFL$^*$, we use an empirical procedure based on energy conservation. To this purpose, we slightly modify the test case so that the energy of the system remains constant. Specifically, we consider the functions in \eqref{non_poly_sol} only to prescribe the initial condition, take homogeneous Dirichlet boundary conditions and impose null source terms. Stability is assessed based on the relative energy increase, limited to $\varepsilon := 5\%$ of the initial or previous time-step energy. At the time step $t^n$, the energy is evaluated as 
\begin{equation}
E^n := \frac{1}{2} \|\bd{m}_{\cal{T}}(t^n)\|^2_{\bd{L}^2(\rho^\sc{f};\domain{f})} + \frac{1}{2} \|p_\cal{T}(t^n) \|^2_{L^2(\frac{1}{\kappa};\domain{f})} + \dfrac{1}{2} \|\bd{v}_{\cal{T}}(t^n)\|^2_{\bd{L}^2(\rho^\sc{s}; \domain{s})} + \frac{1}{2} \|\bbm{s}_{\cal{T}}(t^n)\|^2_{\bbm{L}^2(\bbm{C}^{-1};\domain{s})} .
\end{equation}
The value of CFL$^*$ is determined by reducing iteratively the total number of time steps by $\delta := 1\%$ after each stable simulation. The value of CFL$^*$ is then identified by bracketing the transition between stability and instability. The procedure is summarized in \hyperref[algo_CFL]{\Cref{algo_CFL}}.
\begin{algorithm}
\caption{Bounding of the CFL Stability Limit}
\begin{algorithmic}[1]
\While{simulation is stable}
\For{$n = 1$ to $N$}
\State Compute energy $E_n$
\State Compute relative energy increase $\displaystyle \Delta E := \underset{n}{\max} \left(\frac{|E_n - E_0|}{E_0},\frac{|E_n - E_{n-1}|}{E_{n-1}}\right)$
\If{$\Delta E > \varepsilon$}
\State Flag $N$ as the first unstable time step
\State \textbf{Break}
\EndIf
\EndFor
\State Flag $N$ as the last stable time step
\State Decrease $N$ by $\delta N$
\EndWhile
\end{algorithmic}
\label{algo_CFL} 
\end{algorithm}

In \hyperref[CFL_decoupled]{\Cref{CFL_decoupled}}, we first analyze the stability of the pure acoustic and elastic problems to determine the optimal values of $\eta^\sc{f}$ and $\eta^\sc{s}$ in each subdomain, \textit{i.e.} the values that maximize CFL$^*$. We obtain $\eta^\sc{f}_\star = 0.8$ and $\eta^\sc{s}_\star = 1.5$. These values are retained in the rest of the paper. We also notice that the stability coefficient CFL$^*$ behaves as $\rm{min}(\eta+c_1, \frac{1}{\eta}+c_2)$ for suitable constants $c_1$ and $c_2$, consistently with the observations reported in \cite{MEKG_2025} for the spectral radius of the stiffness matrix.
\begin{figure}[!htb]
\centering
\includegraphics[width=0.495\textwidth]{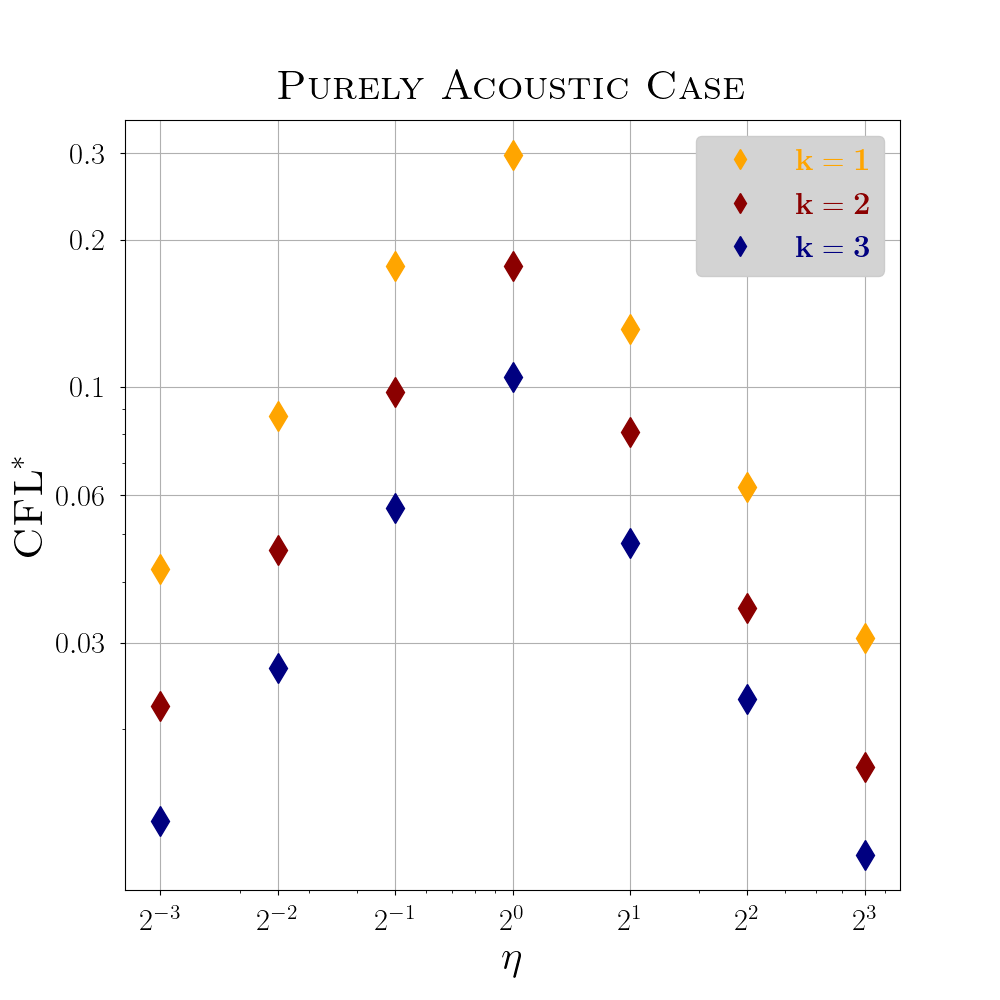}
\includegraphics[width=0.495\textwidth]{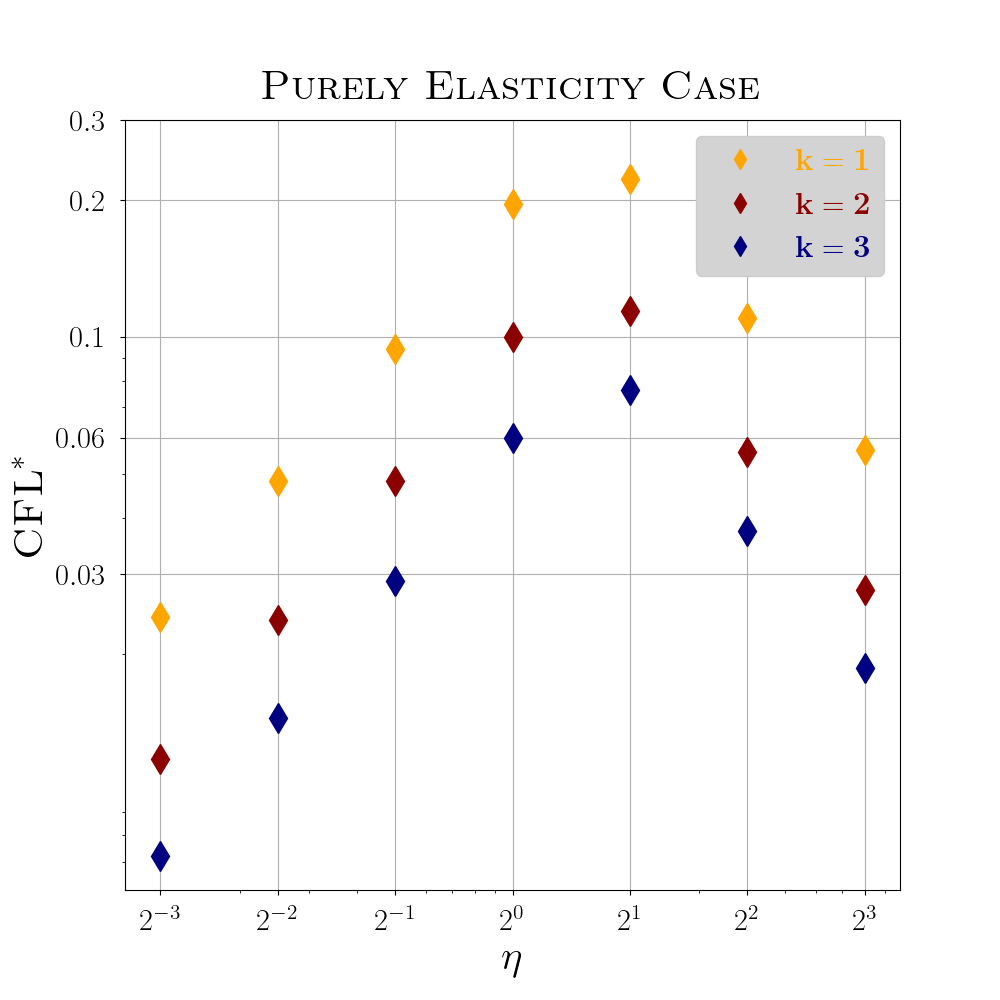}
\caption{CFL$^*$ in the equal-order setting for the pure acoustic and elastic cases for $k \in \{1,2,3\}$.}
\label{CFL_decoupled}
\end{figure}

In \hyperref[tab::ratio_ERK]{\Cref{tab::ratio_ERK}}, we report the CFL$^*$ values determined using \hyperref[algo_CFL]{\Cref{algo_CFL}}. These tables also provide the ratios of the CFL$^*$ coefficients between the different number of stages of the ERK$(s)$ schemes and the different polynomial orders, respectively. It can be observed that these ratios slightly increase with $s$, showing that increasing the number of stages improves stability. Moreover, as expected, CFL$^*$ decreases with the polynomial degree, essentially as $k^{-1}$. Finally, in \hyperref[tab::ratio_geo]{\Cref{tab::ratio_geo}}, we analyze the impact of the mesh geometry on the value of CFL$^*$. This value is only mildly sensitive to the mesh geometry, although meshes with higher face counts seem to slightly enhance the scheme stability.
\begin{table}[!htb]
\centering
\begin{tabular}{|c|ccc|ccc|ccc|}
\cline{2-10} 
\multicolumn{1}{c|}{} & \multicolumn{3}{c|}{$k=1$} & \multicolumn{3}{c|}{$k=2$} & \multicolumn{3}{c|}{$k=3$}\\
\cline{2-10} 
\multicolumn{1}{c|}{} & $s=2$ & $s=3$ & $s=4$ & $s=2$ & $s=3$ & $s=4$ & $s=2$ & $s=3$ & $s=4$ \\ 
\hline
CFL$^*$ & 0.205 & 0.253 & 0.282 & 0.099 & 0.123 & 0.138 & 0.063 & 0.079 & 0.087 \\
\sc{Ratio wrt $s$} & 1 & 1.23 & 1.38 & 1 & 1.24 & 1.39 & 1 & 1.25 & 1.38 \\
\sc{Ratio wrt $k$} & 1 & 1 & 1 & 0.48 & 0.49 & 0.49 & 0.31 & 0.31 & 0.31 \\
\hline   
\end{tabular}
\caption{CFL$^*$ coefficient (and ratios thereof) for ERK$(s)$, $s \in \{2,3,4\}$ and $k \in \{1,2,3\}$.}
\label{tab::ratio_ERK}
\end{table}

\begin{table}[!htb]
\centering
\begin{tabular}{|cc|ccc|ccc|ccc|}
\cline{3-11} 
\multicolumn{2}{c|}{} & \multicolumn{3}{c|}{$s=2$} & \multicolumn{3}{c|}{$s=3$} & \multicolumn{3}{c|}{$s=4$} \\ 
\hline
\multicolumn{2}{|c|}{\sc{Meshes}} & $\triangle$ & $\square$ & {\large $\varhexagon$} & $\triangle$ & $\square$ & {\large $\varhexagon$} & $\triangle$ & $\square$ & {\large $\varhexagon$}\\
\hline 
\multicolumn{1}{|c}{\multirow{2}{*}{$k=1$}} & CFL$^*$ & 0.191 & 0.205 & 0.264 & 0.238 & 0.253 & 0.329 & 0.265 & 0.282 & 0.363 \\
\multicolumn{1}{|c}{} & \sc{Ratio} & 1 & 1.07 & 1.38 & 1 & 1.06 & 1.38 & 1 & 1.06 & 1.37 \\
\hline
\multicolumn{1}{|c}{\multirow{2}{*}{$k=2$}} & CFL$^*$ & 0.106 & 0.099 & 0.136 & 0.133 & 0.123 & 0.170 & 0.147 & 0.138 & 0.188 \\
\multicolumn{1}{|c}{} & \sc{Ratio} & 1 & 0.93 & 1.28 & 1 & 0.92 & 1.28 & 1 & 0.94 & 1.28 \\
\hline
\multicolumn{1}{|c}{\multirow{2}{*}{$k=3$}} & CFL$^*$ & 0.072 & 0.063 & 0.082 & 0.090 & 0.079 & 0.102 & 0.100 & 0.087 & 0.115 \\
\multicolumn{1}{|c}{} & \sc{Ratio} & 1 & 0.88 & 1.13 & 1 & 0.88 & 1.13 & 1 & 0.87 & 1.15 \\
\hline   
\end{tabular}
\caption{CFL$^*$ coefficient (and ratios thereof) for ERK$(s)$, $s \in \{2,3,4\}$ and $k \in \{1,2,3\}$.}
\label{tab::ratio_geo}
\end{table}

\paragraph{Efficiency study.} We now assess the efficiency of the explicit and implicit schemes using the analytical solution \eqref{non_poly_sol}. Accuracy is quantified by computing the error in the $L^2$-norm of the dG variables at the final time step over the computational domain. The efficiency of the various schemes is compared in terms of error versus CPU time. For a fair comparison, implicit schemes employ an iterative solver with a tolerance defined as $\sc{Tol}_{\ell} := 2^{-\ell(k+1)}\sc{Tol}_0$ so that at each level of space refinement, the convergence criterion for the iterative solver is decreased appropriately. Moreover, for both implicit and explicit schemes of order $s$, the time steps are defined as $\Delta t_{\ell} := 2^{-\ell\frac{k+1}{s+1}} \Delta t_0$ so as to ensure that, at each level of space refinement, the space and time discretization errors remain balanced. The reference parameters $\sc{Tol}_0$ and $\Delta t_0$ are determined empirically. Specifically, we compute the error associated with various time step values and observe a stagnation of the total error once the time discretization error becomes negligible compared to the spatial discretization error. We then define $\Delta t_0$ as the largest time step for which this stagnation is observed. The reference tolerance $\sc{Tol}_0$ is selected slightly smaller than the spatial discretization error to ensure that the iterative solver does not dominate the total error. This procedure, which is only considered for $\ell=0$, yields an optimized implicit time step for solving the problem. The time step for explicit schemes is instead dictated by the CFL condition \eqref{CFL} for all $\ell \ge 0$.

In \hyperref[efficiency_ERK_SDIRK]{\Cref{efficiency_ERK_SDIRK}} and \hyperref[efficiency_champion]{\Cref{efficiency_champion}}, we first consider the test case in which the spatial discretization error dominates. In \hyperref[efficiency_ERK_SDIRK]{\Cref{efficiency_ERK_SDIRK}}, our objective is to compare explicit schemes on the one hand, and implicit schemes on the other hand. The first observation, valid for both explicit and implicit schemes, is that increasing the polynomial degree significantly improves efficiency. The second observation is that, in the regime where the spatial error dominates, all the explicit schemes exhibit similar efficiency, and the same holds for implicit schemes. Nevertheless, ERK(2) and SDIRK(3,4) appear to be slightly more efficient than the other schemes within their respective category. In \hyperref[efficiency_champion]{\Cref{efficiency_champion}}, the left panel compares the efficiency of ERK$(2)$ and SDIRK$(3,4)$ schemes. We observe that the implicit scheme is more efficient than the explicit one. This gain is due to the larger time steps allowed by the implicit scheme and the elimination of cell unknowns by static condensation, which compensates for the additional cost of solving the linear system. In the right panel, we evaluate the efficiency of ERK$(4)$ on meshes composed of different cell geometries: simplices, polygons, and squares. We observe that, overall, the efficiency remains similar across the three mesh types. However, at low polynomial orders, polygonal meshes appear to be slightly more efficient, as their more favorable CFL condition allows for somewhat larger time steps. Conversely, at higher orders, the large number of faces in polygonal meshes becomes a drawback. Indeed, the better CFL condition no longer offsets the additional computational cost due to the handling of extra faces. Hence, while the overall efficiency remains comparable, we can conclude that polygonal meshes are preferable at low order, whereas at higher order, it is advantageous to favor cell geometries with fewer faces, such as simplices.
\ifHAL
\begin{figure}[!htb]
\centering
\includegraphics[width=0.475\textwidth]{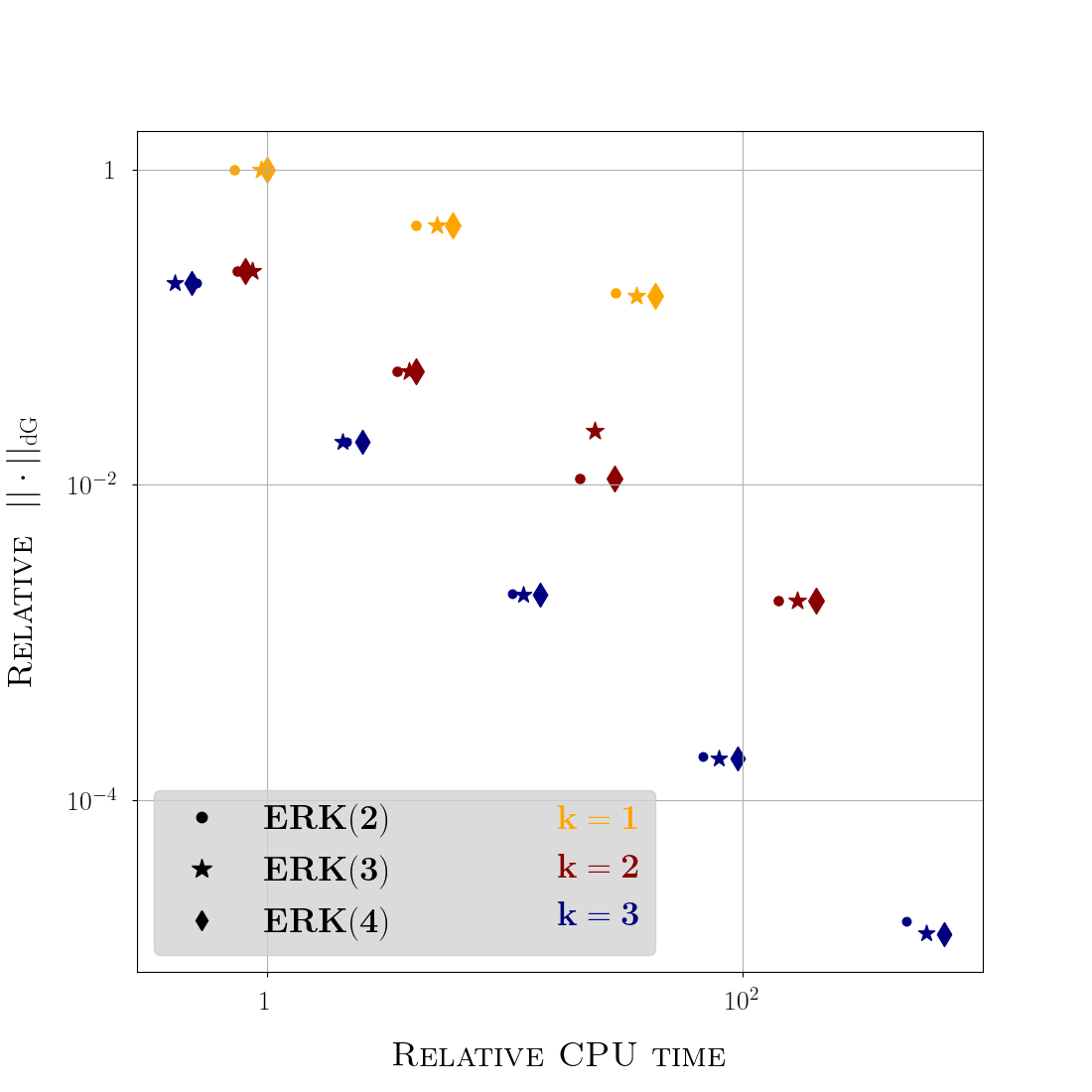}
\includegraphics[width=0.475\textwidth]{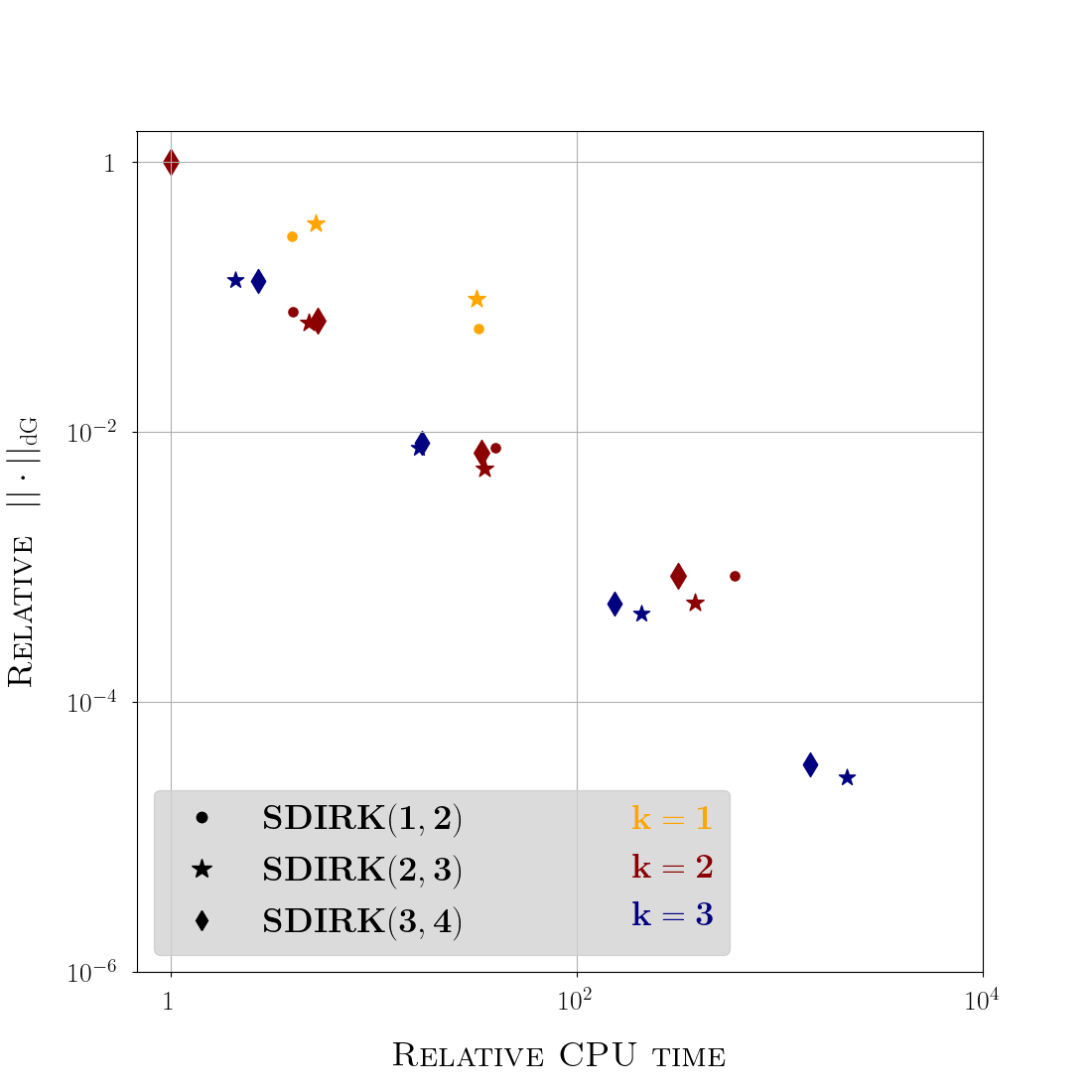}
\caption{Test case with dominant spatial error. \textbf{Left panel:} Efficiency comparison for ERK$(s)$ for $s \in \{2,3,4\}$. \textbf{Right panel:} Efficiency for SDIRK$(s,s+1)$ for $s \in \{1,2,3\}$.}
\label{efficiency_ERK_SDIRK}
\end{figure}
\begin{figure}[!htb]
\centering
\includegraphics[width=0.475\textwidth]{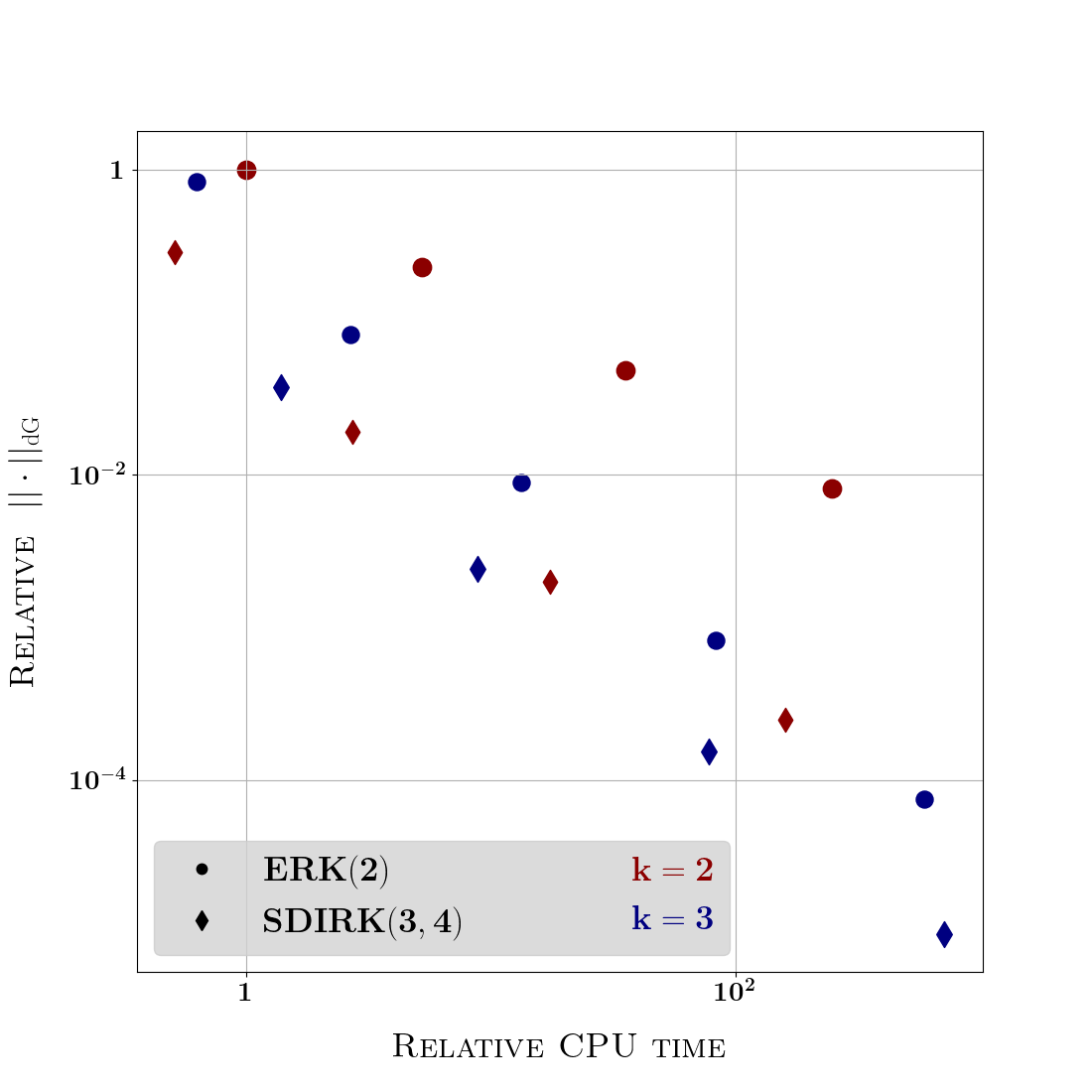}
\includegraphics[width=0.475\textwidth]{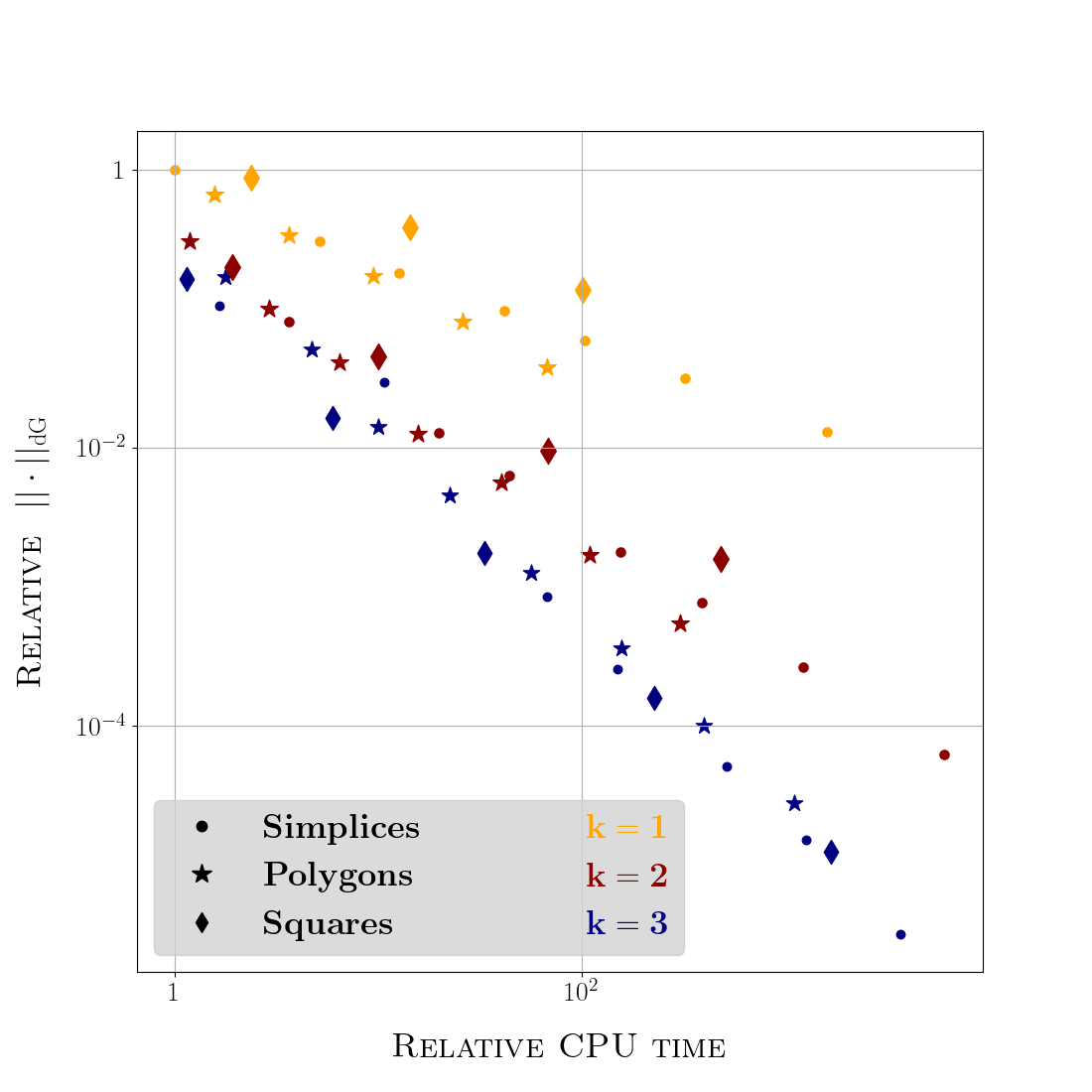}
\caption{Test case with dominant spatial error. \textbf{Left panel:} Efficiency comparison for ERK$(2)$ and SDIRK$(3,4)$. \textbf{Right panel:} Efficiency comparison for ERK$(4)$ on meshes composed of different cell geometries.}
\label{efficiency_champion}
\end{figure}  
\fi

In \hyperref[efficiency_time_dominant]{\Cref{efficiency_time_dominant}}, we now consider the test case where the time error is dominant. In this setting, the time steps used for the implicit schemes are smaller than the stability limit enforced by \eqref{CFL}. Therefore, there is no advantage in using implicit schemes. Regarding the cell geometry, since the stability constraint is no longer the limiting factor, reducing the number of faces is the most effective approach, regardless of the approximation order.
\ifHAL
\begin{figure}[!htb]
\centering
\includegraphics[width=0.495\textwidth]{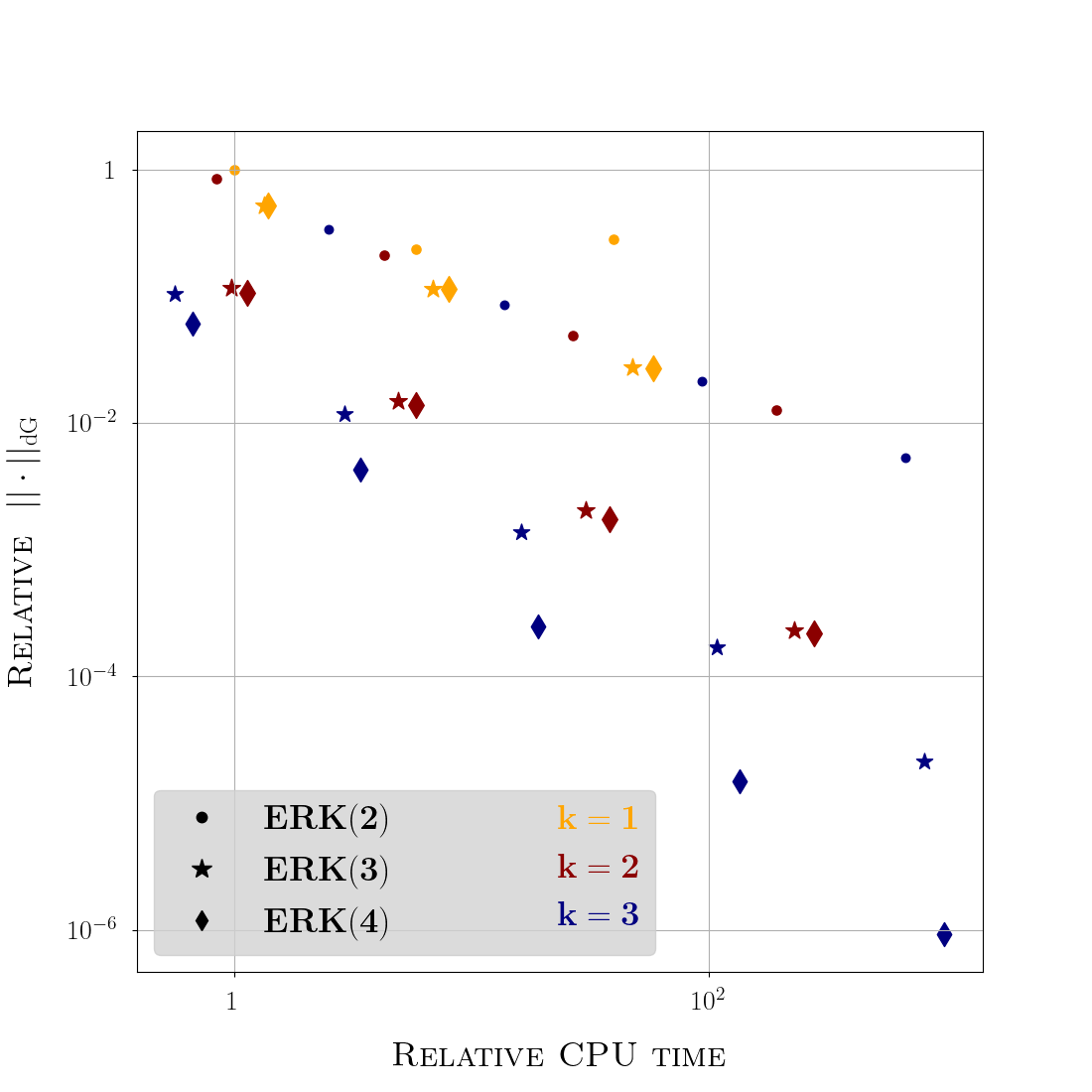}
\includegraphics[width=0.495\textwidth]{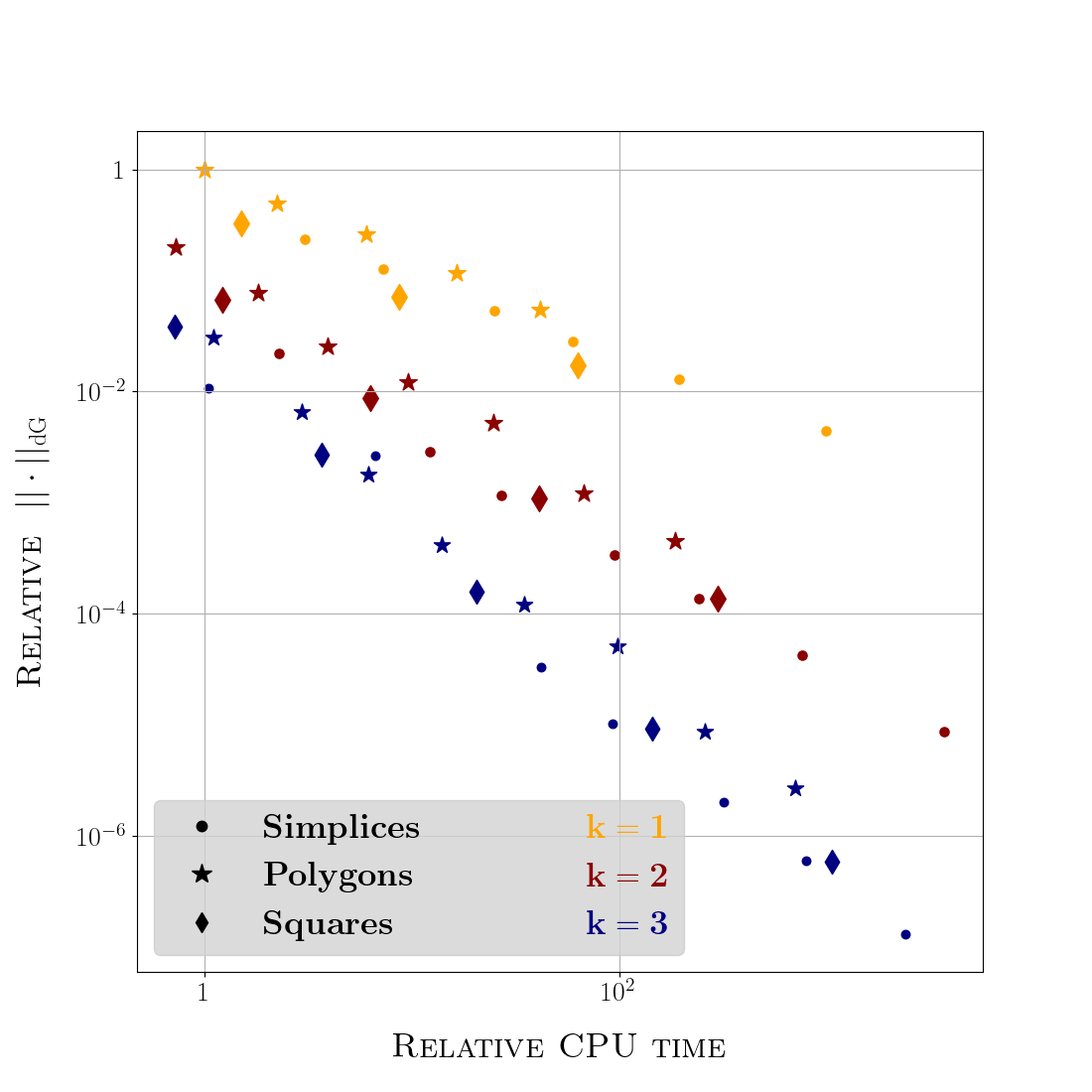}
\caption{Test case with dominant time error. \textbf{Left panel:} Efficiency comparison for ERK$(s)$ for $s \in \{2,3,4\}$. \textbf{Right panel:} Efficiency comparison for ERK$(4)$ on meshes composed of different cell geometries.}
\label{efficiency_time_dominant}
\end{figure}
\fi

\FloatBarrier
\subsection{Accuracy study and efficiency on a Ricker wavelet}

Here, we investigate the accuracy of the proposed method to simulate the propagation of an elasto-acoustic wave through a bilayered, stratified medium. Referring to \hyperref[domain]{\Cref{domain}}, let $H$ denote the height of the whole domain, $H^\sc{f}$ and $H^\sc{s}$ the height of the fluid and solid subdomain respectively, and $L$ the width of both subdomains. We define $\domain{f} := (0, L) \times (0, H^\sc{f})$ for the fluid subdomain and $\domain{s} := (0, L) \times (-H^\sc{s}, 0)$ for the solid subdomain. Homogeneous Dirichlet boundary conditions are applied, the source terms are set to zero, and the perturbation to quiescence is due to the initial velocity field that is a Ricker wavelet centered at the point $(x_c, y_c) \in \Omega^\sc{f}$ (indicated in purple in Figure \ref{domain}), expressed as
\begin{equation}
\bd{m_0}(x, y) := \theta \exp\bigg(\!\!-\pi^2\frac{r^2}{\lambda^2}\bigg)(x-x_c, y-y_c)^\dagger,
\label{ricker_ic}
\end{equation}
where $\theta := 1 \left[\rm{Hz}\right]$, $\lambda := \frac{c^\sc{f}_\sc{p}}{f_c}[\rm{m}]$ with $f_c := 10\left[\rm{Hz}\right]$, and $r^2 := (x-x_c)^2 + (y-y_c)^2$. 

Three virtual sensors record the properties of the propagating elasto--acoustic wavefield: $\cal{S}^{\sc{f}}$, $\cal{S}^{\sc{s}}$ and $\cal{S}^{\sc{i}}$, located in the fluid subdomain, in the solid subdomain and at the interface, respectively. The acoustic sensor $\cal{S}^{\sc{f}}$ records the acoustic pressure $p_{\Tf}$ and the fluid velocity $\bd{m}_{\Tf}$, whereas the elastic sensor $\cal{S}^{\sc{s}}$ records the elastic velocity $\bd{v}_{\Ts}$ and the stress tensor $\bbm{s}_{\Ts}$. The sensor located at the interface records the four quantities $p_{\Ff}$, $\bd{m}_{\Tf}$, $\bd{v}_{\Fs}$ and $\bbm{s}_{\Ts}$.
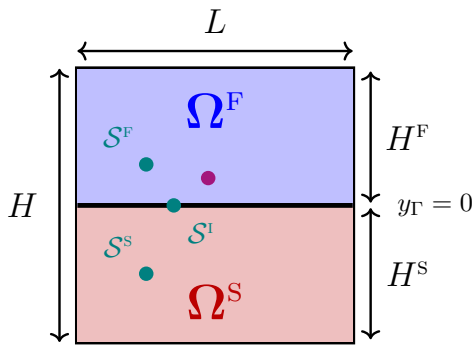
\begin{figure}[!htb]
\centering 
\ifHAL
\renewcommand{\a}{0.55}
\fi
\resizebox{\a\textwidth}{!}{
\begin{tikzpicture}[scale=4]
\draw[line width = 2, color = black] (0,0) rectangle (1,1);
\fill[color=ceared!25] (0, 0) rectangle (1, 0.5);
\fill[color=blue!25] (0, 0.5) rectangle (1,1);
\draw[line width = 2, color = black] (0,0.5) -- (1,0.5);
\draw[<->, line width = 1] (-0.065,0) -- (-0.065,1);
\draw[<->, line width = 1] (0,1.065) -- (1,1.065);
\draw[<->, line width = 1] (1.065,0.01) -- (1.065,0.49);
\draw[<->, line width = 1] (1.065,0.51) -- (1.065,0.99);
\node at (0.5,0.85)  {\huge \color{blue}   $\bd{\domain{f}}$};
\node at (0.5,0.15)  {\huge \color{ceared} $\bd{\domain{s}}$};
\node at (0.5,1.175)  {\Large $L$};
\node at (-0.2,0.5) {\Large $H$};
\node at (1.2,0.25) {\Large $H^\sc{s}$};
\node at (1.2,0.75) {\Large $H^\sc{f}$};
\node at (1.3,0.5)  {\normalsize $y_{\G}=0$};
\fill[color_dofs]  (0.475,0.6) circle (0.75pt);
\fill[white] (0.25,0.65) circle (0.75pt);
\fill[white] (0.25,0.25) circle (0.75pt);
\fill[white] (0.35,0.50) circle (0.75pt);
\node at (0.15,0.75) {\large \color{black} $\cal{S}^{\sc{f}}$};
\node at (0.45,0.4) {\large \color{black} $\cal{S}^{\sc{i}}$};
\node at (0.15,0.35) {\large \color{black} $\cal{S}^{\sc{s}}$};
\end{tikzpicture}
}
\caption{General setting for the accuracy tests with a Ricker wavelet as an initial condition. White dots indicate the position of the three sensors and the magenta dot indicates the localization of the center of the initial source.}
\label{domain}
\end{figure} 

\subsubsection{Academic material properties}

We first consider an academic case in which the acoustic and elastic media exhibit identical densities and propagate the acoustic waves and the elastic S-waves at the same velocity as compressional acoustic waves, \textit{i.e.},
\begin{equation}
\rho^\sc{f} := \rho^\sc{s} := 1, \qquad c_{\sc{p}}^\sc{s} := \sqrt{3}, \qquad c_{\sc{p}}^\sc{f} := c_{\sc{s}}^\sc{s} := 1.
\label{academic_properties}
\end{equation}  
Concerning the geometry, we set $L := H := 1~[\rm{m}]$, $H^\sc{f} := H^\sc{s} := 0.5~[\rm{m}]$, and the pulse is located at $x_c := 0~[\rm{m}]$, $y_c := 0.125~[\rm{m}]$ in the fluid subdomain. The simulation time is set to $T_{\rm{f}} := 1 ~[\rm{s}]$.

As observed below equations \eqref{weak_form_acoustic_eq}-\eqref{weak_form_elastic_eq}, the coupling conditions on ${\Gamma}$ are weakly imposed. Here, we take a closer look at what happens at the interface so as to illustrate numerically the convergence results established in \cite{MEKG_2025}. To this purpose, we use the sensor $\cal{S}^{\sc{i}} := (-0.3,0)~[\rm{m}]$ and we set, for all $t \in \overline{J}$,
\begin{subequations}
\label{errors}
\begin{align}
\label{error1}
\Delta_{\sc{kin}}(t,\cal{S}^{\sc{i}}) & := |(\bd{v}_{\Fs}(t,\cal{S}^{\sc{i}}) - \bd{m}_{\Tf}(t,\cal{S}^{\sc{i}})) \cdot \bd{n}_\G|, \\
\Delta_{\sc{dyn}}(t,\cal{S}^{\sc{i}}) & := \|p_{\Ff}(t,\cal{S}^{\sc{i}})\bd{n}_\G - \bbm{s}_{\Ts}(t,\cal{S}^{\sc{i}}) \cdot \bd{n}_\G \|,
\label{error2}
\end{align}
\end{subequations}
The numerical results are obtained using SDIRK$(3,4)$ and the time step is set to $\Delta t = 0.1 \times 2^{-8}$ (this value is small enough so that space discretization errors dominate). \hyperref[coupling_condition]{\Cref{coupling_condition}} reports the evolution at the discrete time nodes of the errors defined in \eqref{errors}, for various mesh refinement levels $\ell$. The left column shows the kinematic errors \eqref{error1}, whereas the right column shows the dynamic errors \eqref{error2}. The top row corresponds to the polynomial degree $k=1$ and the bottom row to $k=3$. As expected, the numerical error on the pointwise satisfaction of the coupling conditions decreases significantly when increasing the mesh refinement and/or the polynomial degree 
\begin{figure}[!htb]
\centering
\includegraphics[width=0.495\textwidth]{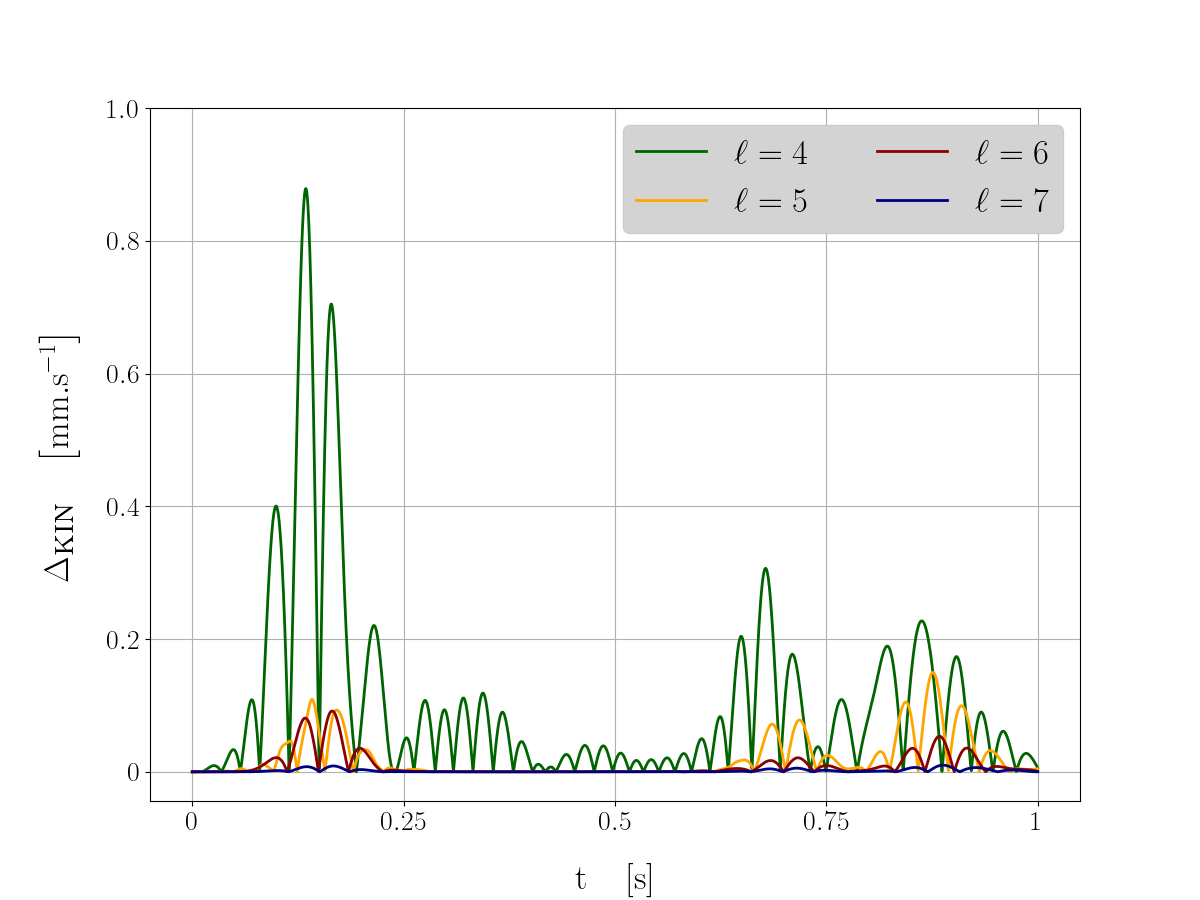}
\includegraphics[width=0.495\textwidth]{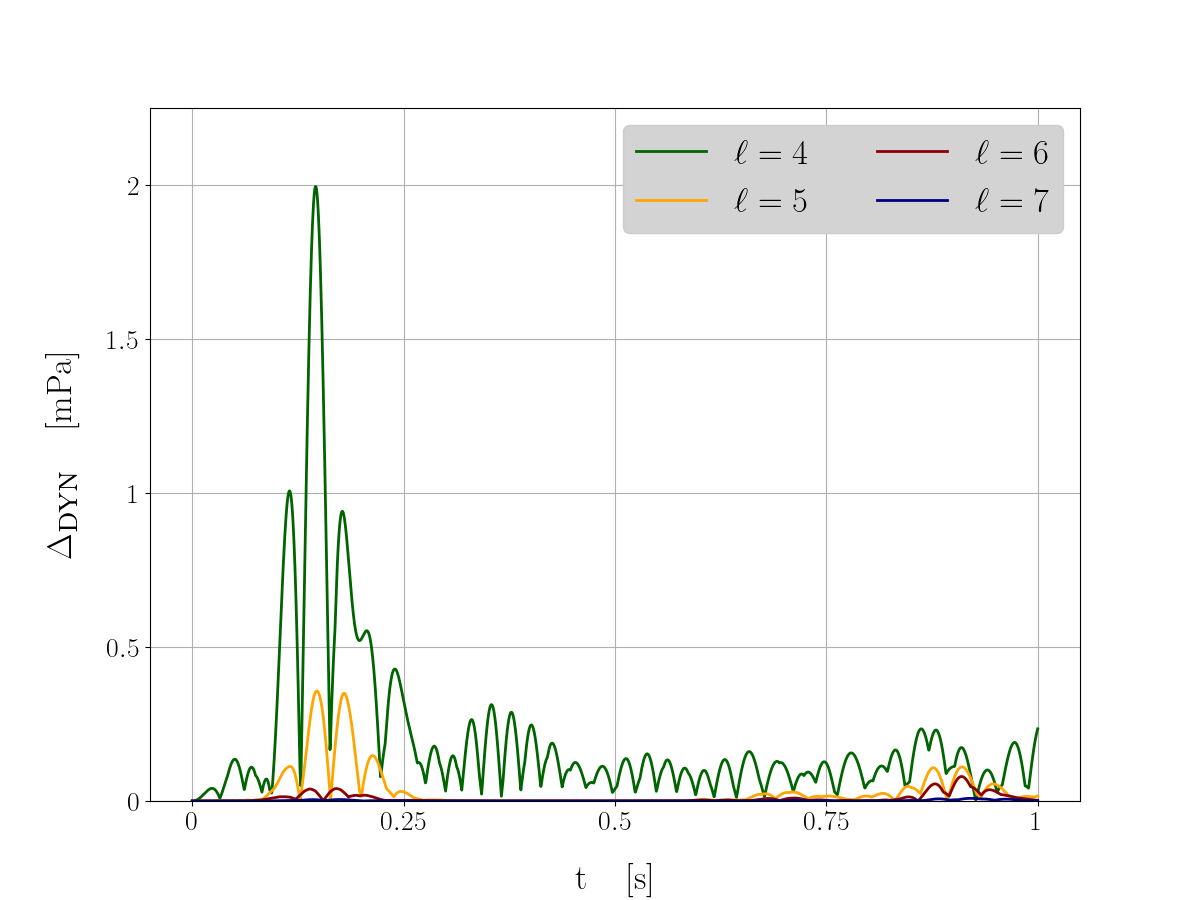}
\includegraphics[width=0.495\textwidth]{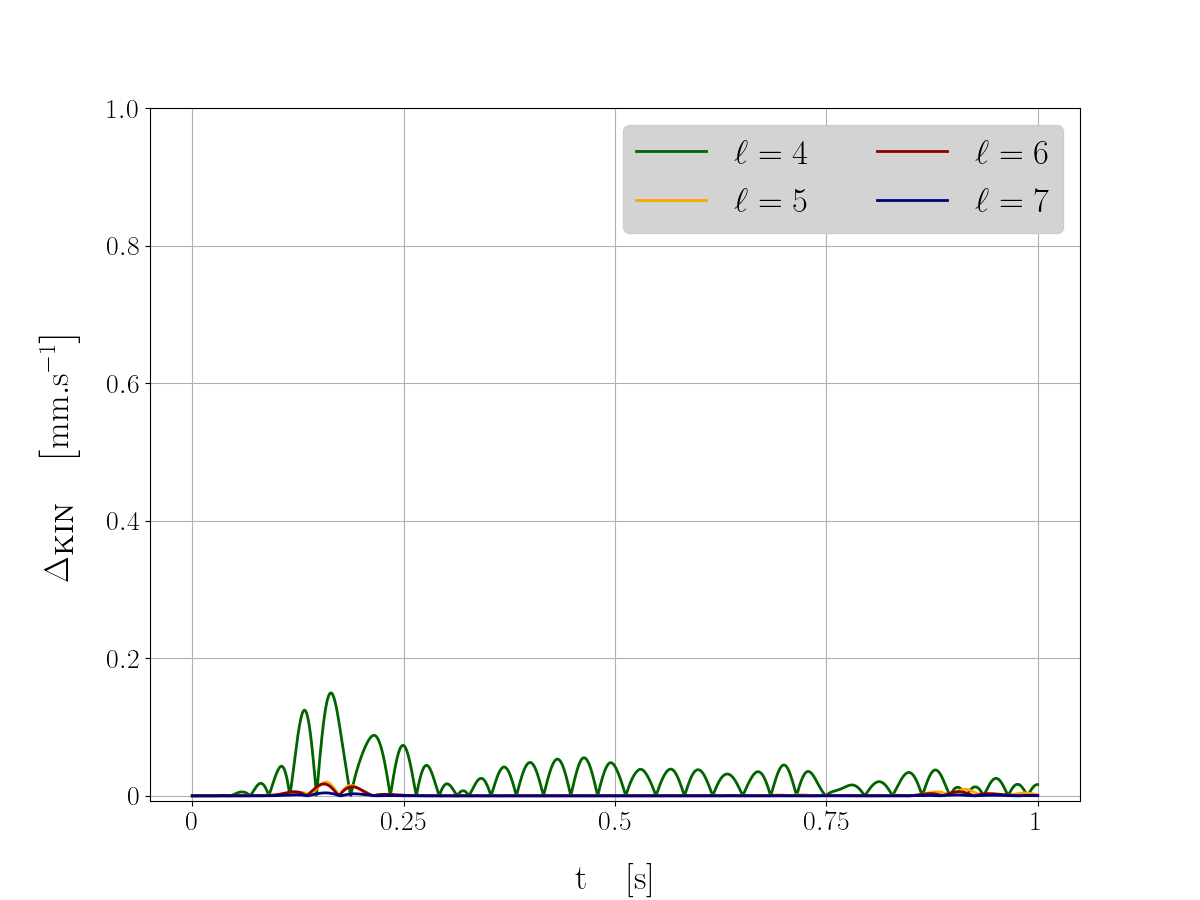}
\includegraphics[width=0.495\textwidth]{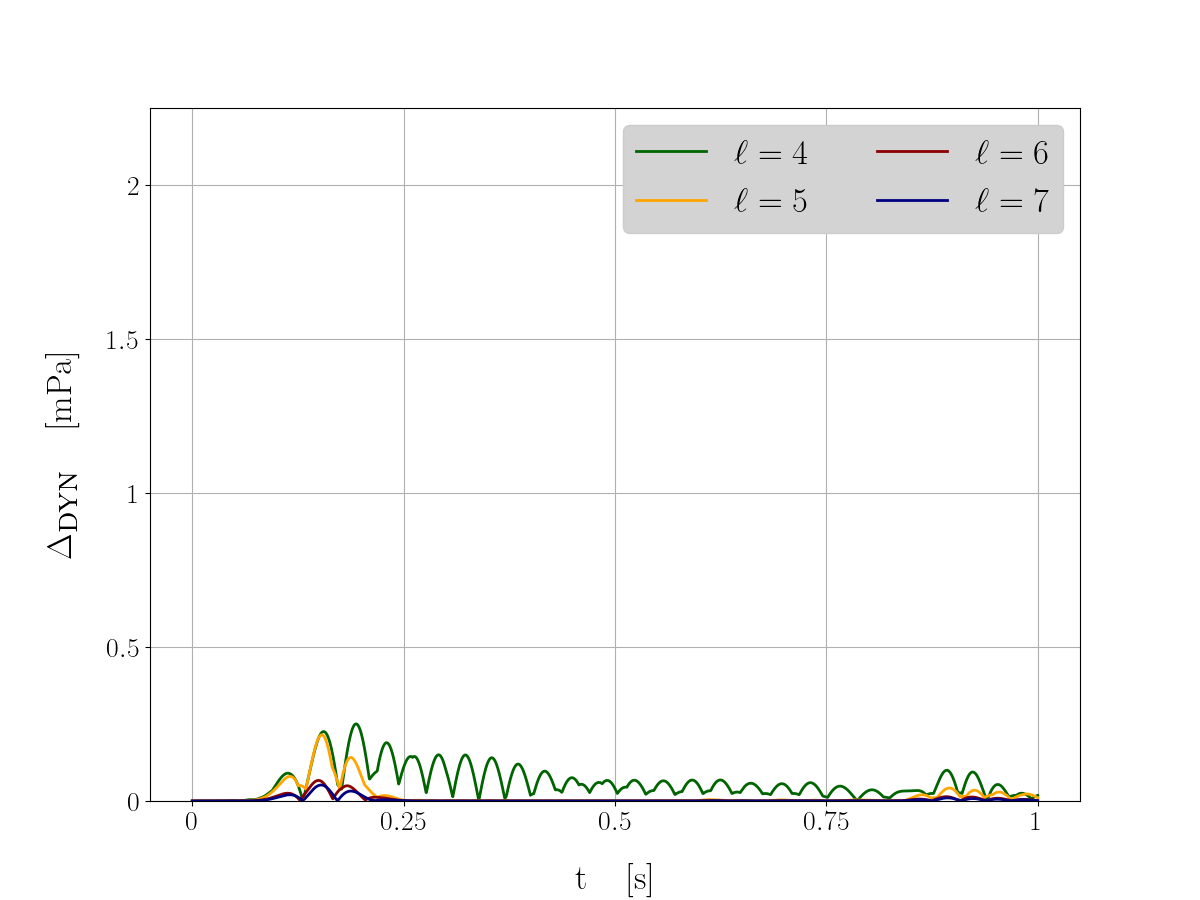}
\caption{Errors on the coupling conditions as a function of the time predicted by SDIRK$(3,4)$ for $k=1$ (top row) and $k=3$ (bottom row). \textbf{Left:} Kinematic errors (see \eqref{error1}). \textbf{Right:} Dynamic errors (see \eqref{error2}).}
\label{coupling_condition}
\end{figure}

\subsubsection{Realistic (strongly contrasted) material properties}

We now consider a test case with a strong contrast of material properties, specifically focusing on hard rock (granite) overlaid by a water layer with the following material properties:
\begin{equation}
\begin{alignedat}{2}
\rho^\sc{f} & := 1025~\rm{kg.m}^{-3}, \qquad & c_{\sc{p}}^{\sc{f}} & := 1500~\rm{m.s}^{-1}, \\
\rho^\sc{s} & := 2690~\rm{kg.m}^{-3}, \qquad & c_{\sc{p}}^\sc{s} & := 6000~\rm{m.s}^{-1}, \qquad c_{\sc{s}}^{\sc{s}} := 3000~\rm{m.s}^{-1}.
\end{alignedat}
\label{granite_water}
\end{equation}
The geometric setup is now $L := \rev{3.8}~[\rm{km}]$, $H := \rev{2.6}~[\rm{km}]$, $H^\sc{f} := \rev{1}~[\rm{km}]$ and $H^\sc{s} := \rev{1.6}~[\rm{km}]$. \rev{The initial condition is a pressure pulse, still defined using the right-hand side of~\eqref{ricker_ic} and} located at $x_c := 0~[\rm{km}]$, $y_c := \rev{0.2}~[\rm{km}]$ in the fluid subdomain. The simulation time is set to $T_{\rm{f}} := \rev{0.5}~[\rm{s}]$. The mesh is composed of \rev{square} elements with size of $h_x = \rev{h_y = 23.6}~[\rm{m}]$.

The results are obtained using SDIRK$(3,4)$ and the time step is set to $\Delta t = 0.1 \times 2^{-5}$. 
\rev{A reference solution is computed using the 2D spectral element solver of the SEM3D package (a co-developed and optimized spectral element code, \url{https://github.com/sem3d/SEM}), which is restricted to \rev{quadrilateral} meshes. This solver relies on a second-order (primal) formulation in time and an explicit Newmark scheme for time integration. The space discretization uses a polynomial degree $k=6$ (i.e., $7$ GLL nodes in each  direction), a mesh size $h=14~\rm{[m]}$, and a time step $\Delta t = 3.5 \times 10^{-5}~[\rm{s}]$.}
\hyperref[ricker]{\Cref{ricker}} displays the particularly rich structure of the pressure field in the fluid subdomain and of the Euclidian norm of the velocity field in the solid subdomain \rev{at the times $t=0.36$~[s] and $t=0.50$~[s]}. We observe the presence of conical and interface (Rayleigh, Scholte) waves, whose accurate capture is commonly used to assess the accuracy of a numerical method. \rev{Since the size of the computational box is relatively small, we observe that some reflections of the P-wave have occurred at $t=0.50$~[s] at the boundary due to the use of Dirichlet conditions.}
\ifHAL
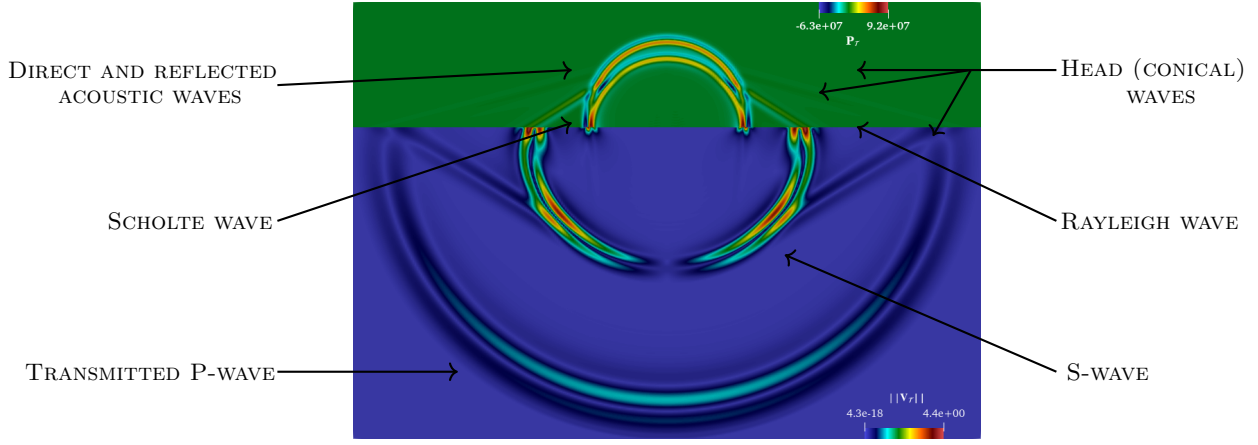
\begin{figure}[!htb]
\centering
\resizebox{\textwidth}{!}{
\begin{tikzpicture}
\node[inner sep=0pt] (img) at (-10,0)  {\includegraphics[width=0.495\textwidth]{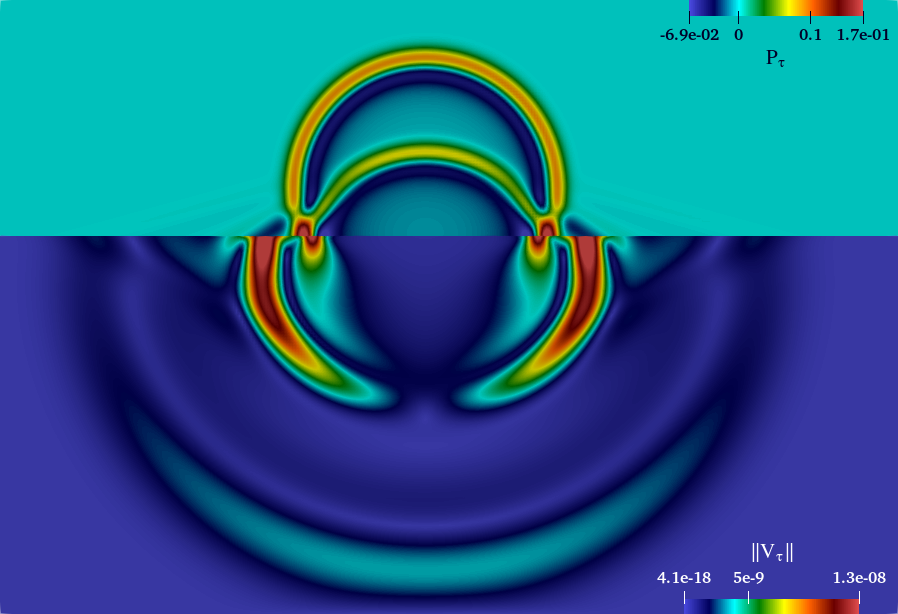}};
\node[inner sep=0pt] (img) at (-1.75,0) {\includegraphics[width=0.495\textwidth]{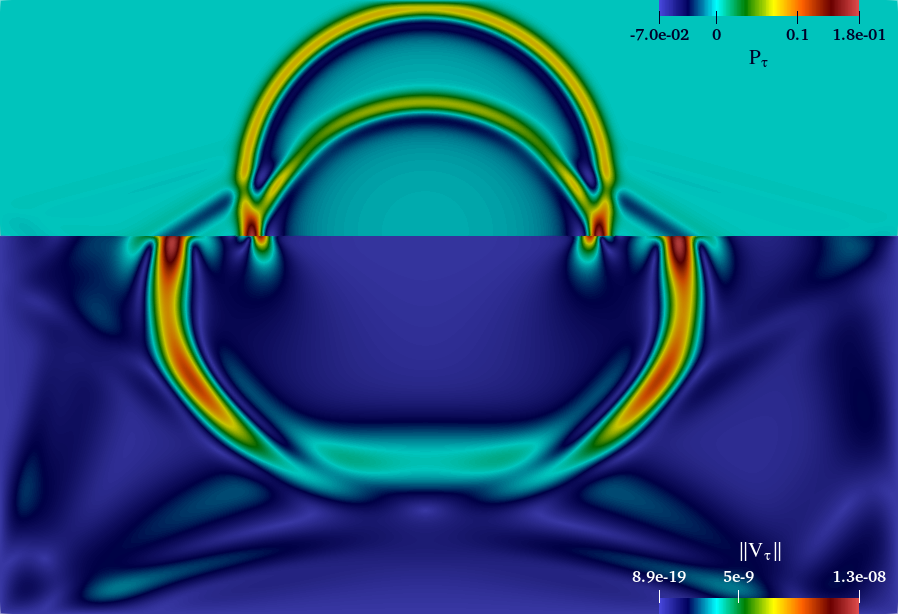}};
\node[fill=none, text=black, font=\footnotesize] at (-7.5,1.375) {\scshape S-wave};
\node[fill=none, text=black, font=\footnotesize] at (-12.1,2.6) {\scshape Direct and reflected};
\node[fill=none, text=black, font=\footnotesize] at (-12.3,2.325) {\scshape acoustic waves};
\node[fill=none, text=black, font=\footnotesize] at (-13+0.2,1.5) {\scshape Transmitted};
\node[fill=none, text=black, font=\footnotesize] at (-13+0.2,1.25) {\scshape P-wave};
\node[fill=none, text=white, font=\footnotesize] at (-1.95,-0.5) {\scshape Scholte wave};
\node[fill=none, text=black, font=\footnotesize] at (1.4-0.3,2.00-0.1) {\scshape Head (conical)};
\node[fill=none, text=black, font=\footnotesize] at (1.4-0.25,1.75-0.125) {\scshape waves};
\node[fill=none, text=black, font=\footnotesize] at (-4.4,2.6) {\scshape Rayleigh wave};
\draw[->, line width=1.25pt, black] (-12.25,2.1) -- (-10.5,2.2);
\draw[->, line width=1.25pt, black] (-12.25,2.1) -- (-10.5,1.35);
\draw[->, line width=1.25pt, black] (-12.75,1) -- (-11,-2.35);
\draw[->, line width=1.25pt, black] (-7.5,1) -- (-9,-0.35);
\draw[->, line width=1.25pt, white] (-2+0.05,-0.25) -- (-3.375+0.025,0.55);
\draw[->, line width=1.25pt, white] (-2+0.05,-0.25) -- (-0.65+0.075,0.55);
\draw[->, line width=1.25pt, black] (-5,2.35) -- (0.65,0.6);
\draw[->, line width=1.25pt, black] (-5,2.35) -- (-4.6,0.6);
\draw[->, line width=1.25pt, black] (1,1.4) -- (0,1);
\draw[->, line width=1.25pt, black] (1,1.4) -- (-4,0.8);
\end{tikzpicture}}
\caption{Spatial distribution of the acoustic pressure (upper side) and the elastic velocity norm (lower side) at \rev{$t = 0.36~[\rm{s}]$ (left) and at $t = 0.50~[\rm{s}]$ (right).}}
\label{ricker}
\end{figure}
\fi
\begin{table}[!htb]
\centering 
\resizebox{\textwidth}{!}{
\rev{\begin{tabular}{|c|c|c|c|cc|cc|c|c|}
\hline
\multicolumn{1}{|c|}{\sc{Schemes}} & \rule{0pt}{3.5ex} \sc{Solver} & k & CFL$^*$ & $\widetilde{\Delta t}$ & \sc{Ratio} & CPU $[\rm{s}]$ & \sc{Ratio} & $\sc{Err}^{\sc{f}}$ & $\sc{Err}^{\sc{s}}$ \\[0.5ex]
\hline 
\multirow{2}{*}{SDIRK$(3,4)$} & \multirow{2}{*}{direct} & $1$ & \multirow{2}{*}{n/a} & 0.398 & 7.4 & 207 & 1 & 5.76e-02 & 2.79e-02 \\
& & $3$ & & 0.398 & 7.4 & 1095 & 1 & 2.33e-03 & 8.58e-04\\
\hline
\multirow{2}{*}{SDIRK$(3,4)$} & \multirow{2}{*}{iterative} & $1$ & \multirow{2}{*}{n/a} & 0.398 & 7.4 & 281 & 1.4 & 5.76e-02 & 2.79e-02\\
&  & $3$ & & 0.398 & 7.4 & 1352 & 1.2 & 2.39e-03 & 9.11e-04\\
\hline
\multirow{2}{*}{ERK$(4)$} & \multirow{2}{*}{n/a} & $1$ & 0.282 & 0.159 & 2.96 & 161 & 0.8 & 5.73e-02 & 2.77e-02\\
&  & $3$ & 0.087 & 0.054 & 1 & 3322  & 3.0 & 1.93e-03 & 5.10e-04\\
\hline
\end{tabular}}}
\caption{Times step $\widetilde{\Delta t}$ (normalized as CFL$^*$), CPU times and errors (see \eqref{error_ricker}) for ERK(4), SDIRK(3,4) and $k \in \{1,3\}$.}
\label{tab::cfl_cpu}
\end{table} 

In \hyperref[tab::cfl_cpu]{\Cref{tab::cfl_cpu}} and \hyperref[sensors_ricker]{\Cref{sensors_ricker}}, we evaluate the accuracy of the numerical predictions at the three sensors located at $\cal{S}^{\sc{f}} := \rev{(-0.5,0.2)}~[\rm{km}]$, $\cal{S}^{\sc{s}} := \rev{(-1.0, -0.5)}~[\rm{km}]$ and \rev{$\cal{S}^{\sc{i}} := (-0.4, 0)~[\rm{km}]$}. The sensors $\cal{S}^{\sc{f}}$ and $\cal{S}^{\sc{s}}$ are located within the two subdomains and positioned so as to capture only the conical waves, whereas the sensor $\cal{S}^{\sc{i}}$ located at the interface also captures the main phases. We compare ERK$(4)$ and SDIRK$(3,4)$ and to avoid the proliferation of results, we focus on the polynomial degrees $k \in \{1,3\}$. CPU times and relative errors with respect to the reference solution, measured at the sensors $\cal{S}^{\sc{f}}$ and $\cal{S}^{\sc{s}}$, are reported in \hyperref[tab::cfl_cpu]{\Cref{tab::cfl_cpu}}. The relative errors are defined as
\begin{equation}
\rev{\sc{Err}^{\sc{f}} := \frac{\|p_{\Tf}(t,\cal{S}^{\sc{f}}) - p_{\sc{ref}}(t,\cal{S}^{\sc{f}})\|_{{\ell}^2(J)}}{\|p_{\sc{ref}}(t,\cal{S}^{\sc{f}})\|_{{\ell}^2(J)}},
\qquad
\sc{Err}^{\sc{s}} :=  \frac{\|\bd{v}_{\Ts}(t,\cal{S}^{\sc{s}}) - \bd{v}_{\sc{ref}}(t,\cal{S}^{\sc{s}})\|_{\bd{\ell}^2(J)}}{\|\bd{v}_{\sc{ref}}(t,\cal{S}^{\sc{s}})\|_{\bd{\ell}^2(J)}},}
\label{error_ricker}
\end{equation}
\rev{where the discrete $\|{\cdot}\|_{\ell^2(J)}$-norm is sampled at the discrete time nodes.}
Notice that all the methods compared in \hyperref[tab::cfl_cpu]{\Cref{tab::cfl_cpu}} lead to comparable errors, so that the CPU comparison is meaningful concerning their efficiency. \rev{Moreover, the errors for $k=3$ are more than twenty times smaller (and below $10^{-3}$) compared to the errors for $k=1$.} For the explicit scheme, we consider time steps just slightly below the stability limit. Instead, for the implicit scheme, we can consider larger time steps. In each case, the time step $\Delta t$ is reported in nondimensional form as $\widetilde{\Delta t} := c_\sc{p}^\sc{s}\frac{\Delta t}{h}$, where $c_\sc{p}^\sc{s}$ is defined in \eqref{granite_water} and is the largest velocity in the domain. We observe in \hyperref[tab::cfl_cpu]{\Cref{tab::cfl_cpu}} that the difference is significant between implicit and explicit time schemes. In particular, for $k=3$, the CPU ratio reaches nearly \rev{3}, since much smaller time steps are required for explicit schemes \rev{(up to seven times smaller)}. This negatively impacts the overall CPU time, confirming as in the previous section the higher efficiency reached by implicit schemes. \rev{For the implicit scheme, we compare the use of a direct solver (with pre-factorization) and an iterative solver (ILU/BiCG) with a relative tolerance on the residual set to $10^{-2}$, achieving comparable results in terms of error and CPU time.}
Referring now to \hyperref[sensors_ricker]{\Cref{sensors_ricker}}, we notice that, for $k=1$, \rev{small discrepancies between the numerical and exact solutions are still visible in the profiles, whereas the predictions for $k=3$ overlap with the reference solution}.

\begin{figure}[!htb]
\centering
\ifHAL
\includegraphics[width=0.495\textwidth]{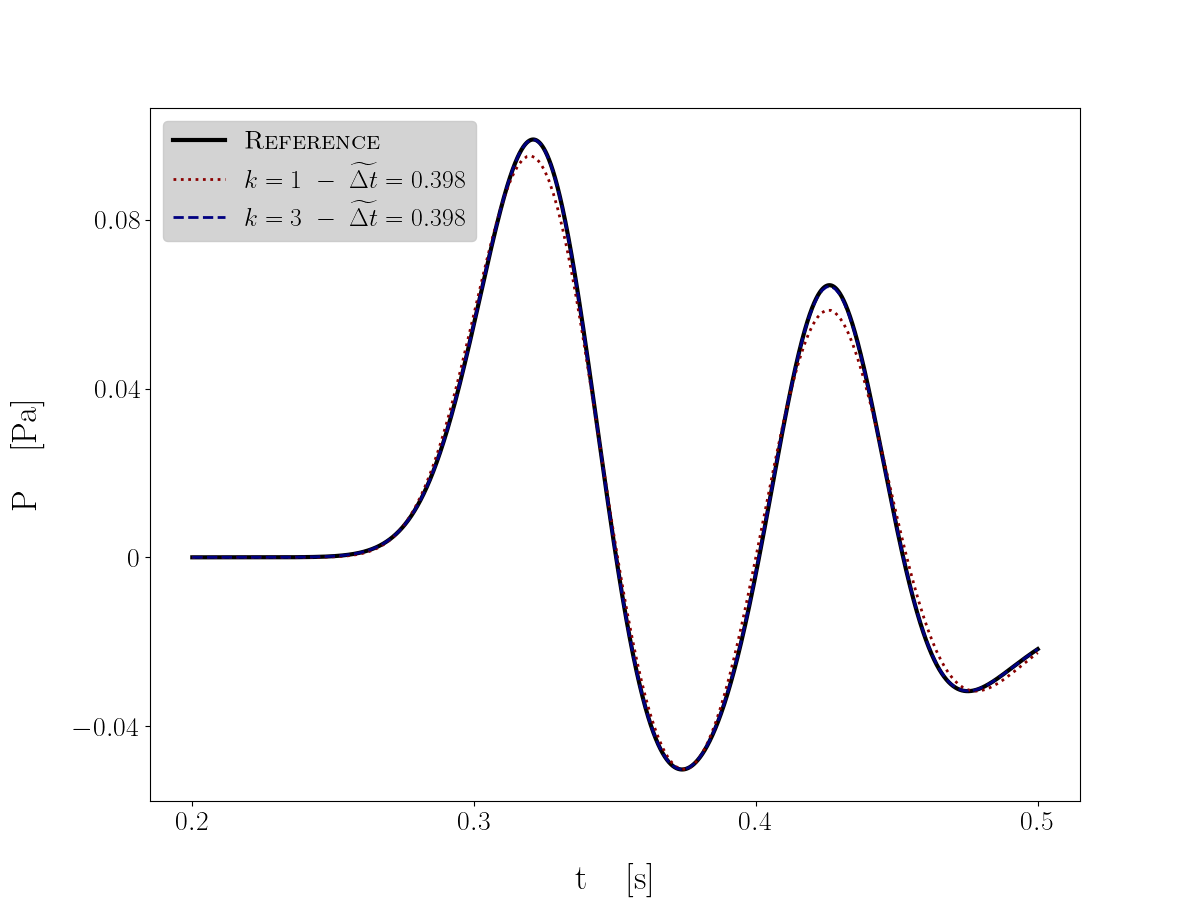}
\includegraphics[width=0.495\textwidth]{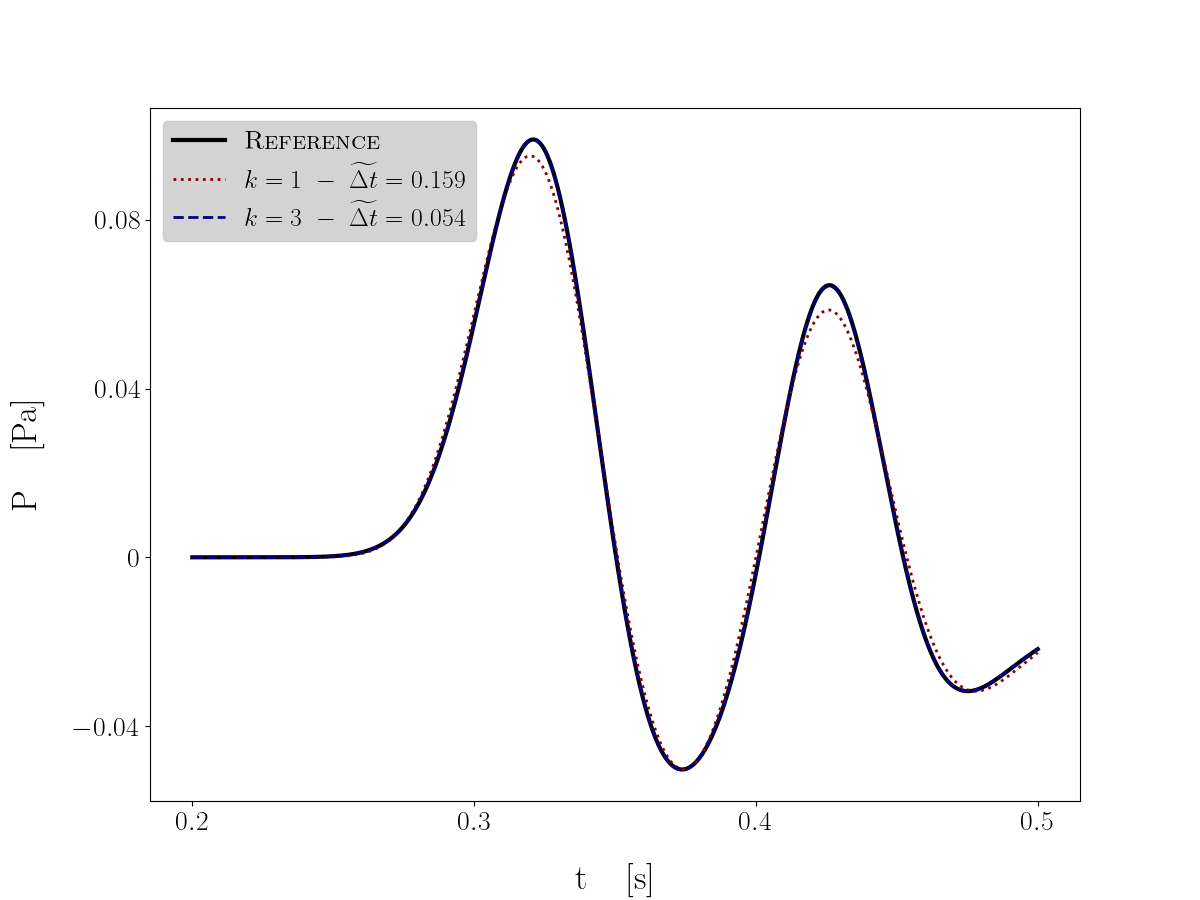}
\includegraphics[width=0.495\textwidth]{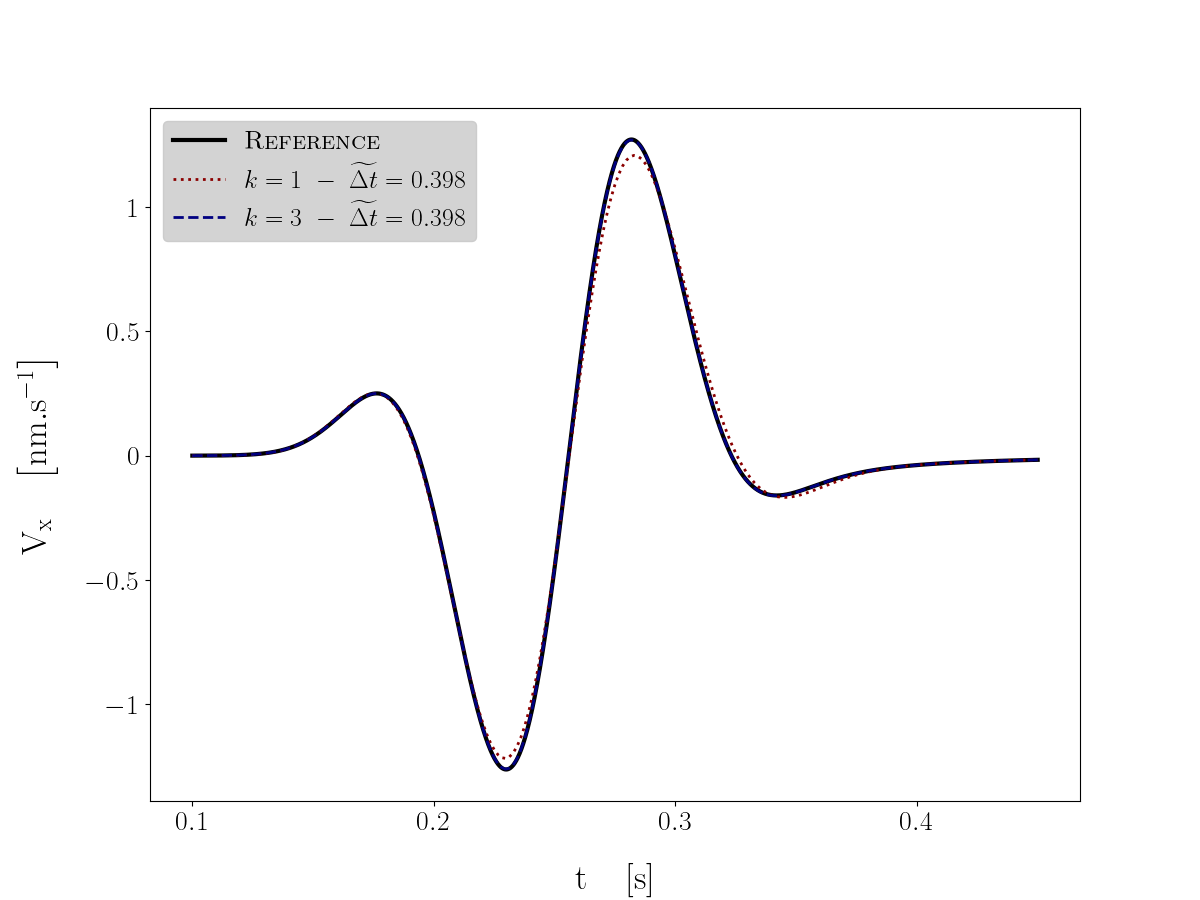}
\includegraphics[width=0.495\textwidth]{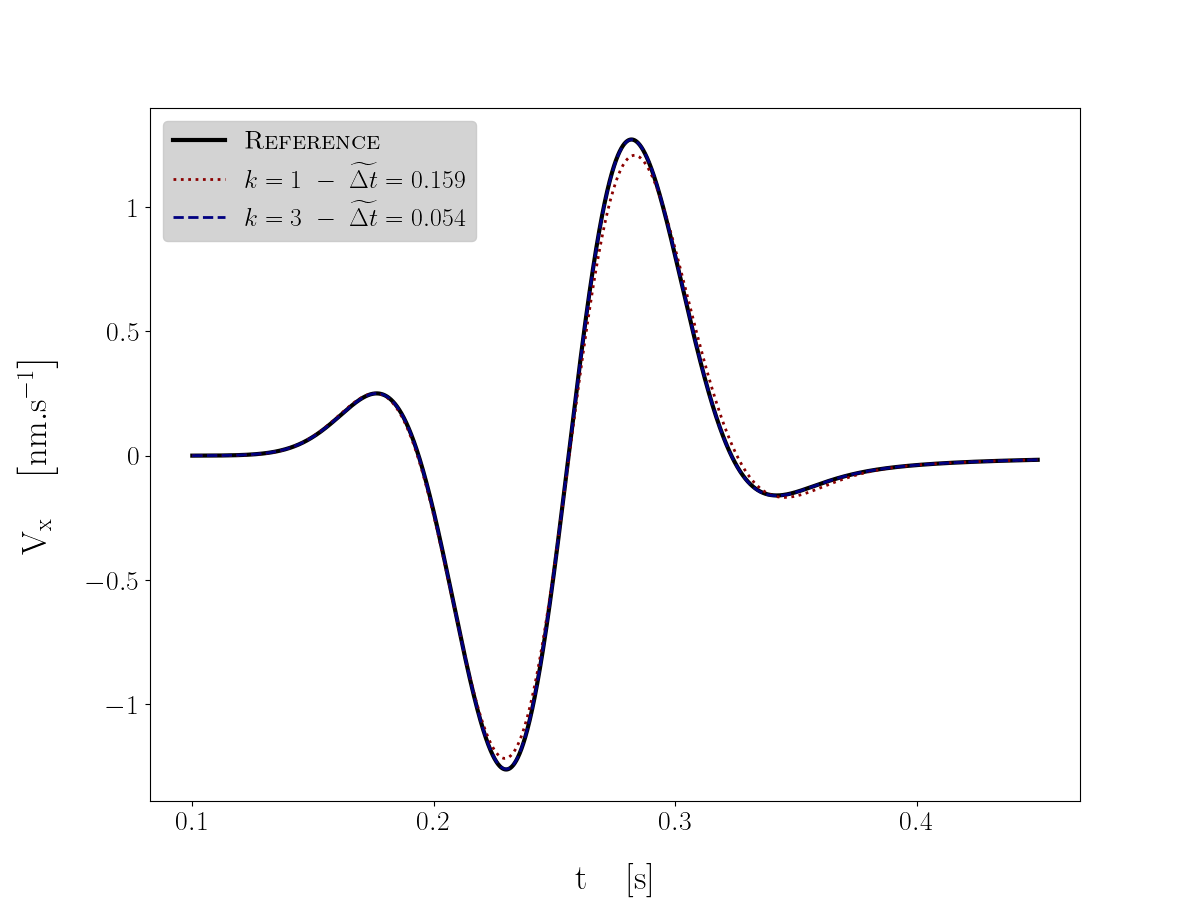}
\includegraphics[width=0.495\textwidth]{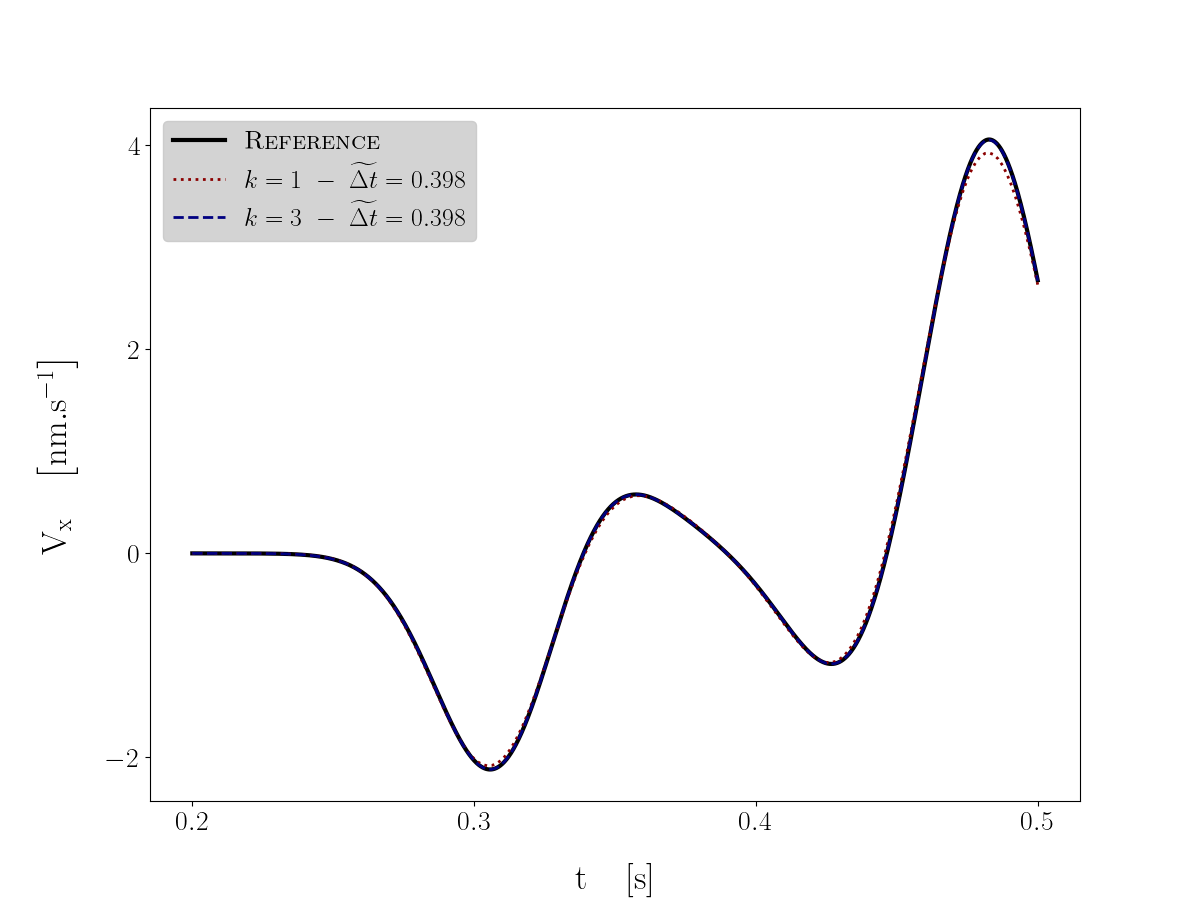}
\includegraphics[width=0.495\textwidth]{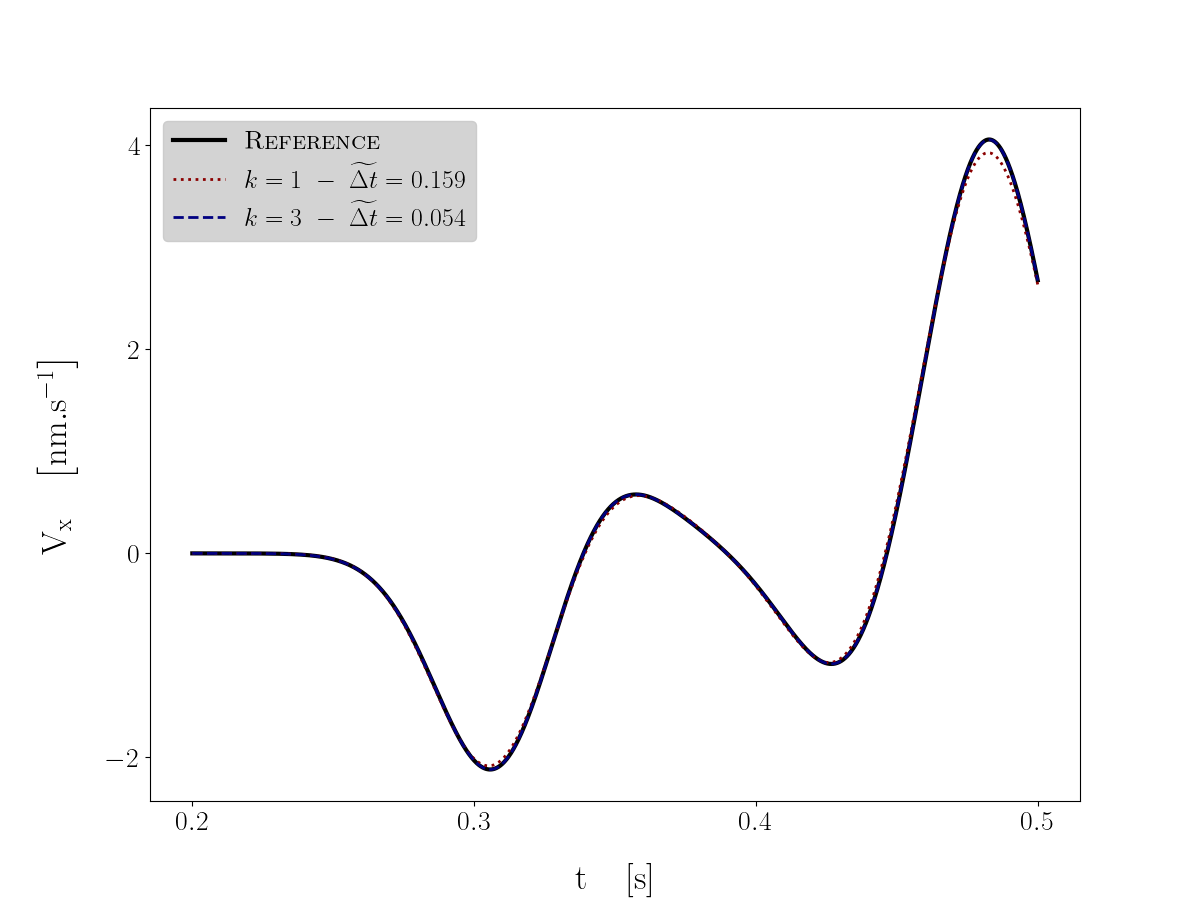}
\fi
\caption{Spatial distribution of the acoustic pressure at $\cal{S}^{\sc{f}}$ (top row), $x$-component of the velocity at $\cal{S}^{\sc{i}}$ (center row) and $\cal{S}^{\sc{s}}$ (bottom row) at times $t \in [0,0.5] [\rm{s}]$. \textbf{Left column:} SDIRK$(3,4)$, $k \in \{1,3\}$, and $\widetilde{\Delta t} = 0.398$. \textbf{Right column:} ERK$(4)$, $k=1$, $\widetilde{\Delta t}=0.159$, and $k=3$, $\widetilde{\Delta t} = 0.054$.}
\label{sensors_ricker}
\end{figure}

\subsection{Geometric flexibility of HHO: a realistic simulation}

Finally, we present a realistic test case of seismoacoustics coupling that is relevant in some geophysical applications. Namely, we are interested in the propagation of seismic waves that can be trapped in a sedimentary basin embedded in a hard rock, or bedrock, and such that large part of the energy can leak into the atmosphere. This phenomenon seems to have already been observed in \cite{HLVHCPMRB_2018, AASSE_2018}, and may be important for the monitoring of natural hazards or anthropogenic activities in areas where seismic arrays are sparse. 

We restrict ourselves to a 2D simulation, because the high-performance computing tools for a realistic 3D simulation are not yet available in the present implementation of \texttt{disk++}. This study is thus to be considered as a proof-of-concept to showcase the benefits of using HHO methods through an ease in meshing procedures (polytopal flexibility). First, we illustrate that there is no need to restrict the spatial discretization to a unique cell geometry. If this may not be that relevant in 2D (where a complete discretization with quadrangles is possible), it is clearly a nontrivial issue in 3D realistic situations \cite{DELAVAUD_2007,CMTBKKS_2010}. Second, we show that HHO allows one to perform straightforwardly simulations using a nonconforming mesh (that is, with hanging nodes).

\subsubsection{Test case setting}

The geological setting shown in the left panel of \hyperref[fig:bassin]{\Cref{fig:bassin}} is rather simplistic, especially concerning the spatial dimensions (less than 10~km in each direction, namely the domain extends from $-4~[\rm{km}]$ to $+4~[\rm{km}]$ laterally, and from $-6~[\rm{km}]$ to $1.4~[\rm{km}]$ vertically, with a solid-fluid boundary located around a height of $+1.1~[\rm{km}]$), whereas interactions of the seismic wavefield with the sedimentary basin-atmosphere can happen several tens of kilometers away from the source \cite{AASSE_2018}. Here, the source term is located under a small moutain (where coupling with the atmosphere is expected, see \cite{CKPPAS_2021}), and a sedimentary basin is located 1~km away. The seismic motion is introduced as an initial condition in velocity, whose spatial distribution is that of a Ricker wavelet (as in \eqref{ricker_ic}) with a central frequency of $f_c = 4~\rm{Hz}$. The center of the Ricker wavelet is located at $(0.563, 0.233)~[\rm{km}]$ (see {\Cref{fig:bassin}}). We consider a simulation time $T_{\rm{f}} = 1.5~[\rm{s}]$, and material properties are reported in the right panel of \hyperref[fig:bassin]{\Cref{fig:bassin}}.
\ifHAL
\renewcommand{\a}{0.8}
\fi
\begin{figure}[!htb]
\resizebox{\a\textwidth}{!}{
\begin{minipage}[c]{0.5\textwidth}
\ifHAL \centering \else \fi
\begin{tikzpicture}
\node[inner sep=0pt] (img) at (0,0) {\includegraphics[height=5.75cm]{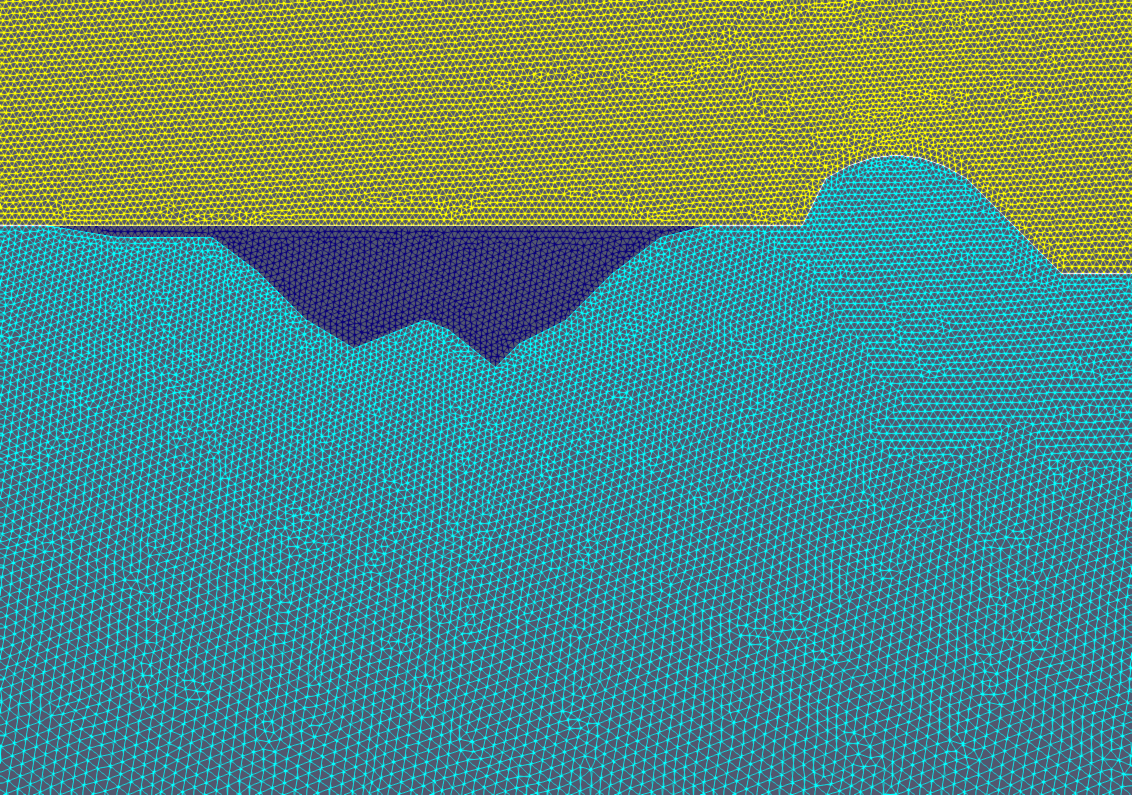}};
\node at (-2.65,2.5)  {\textbf{Atmosphere}};
\node at (-0.9,0.875) {\textbf{Sediments}};
\node at (-3,-2.50)   {\textbf{Bedrock}};
\node at (-1.6,1.5)   {\large \color{black} $\cal{S}^{\sc{s}}$};
\fill[white] (-1.95,1.125) circle (0.075);
\node at (-0.6,2.35) {\large \color{black} $\cal{S}^{\sc{f}}$};
\fill[white] (-1,2) circle (0.075);
\fill[violet] (2.25,-2) circle (0.2);
\end{tikzpicture}
\end{minipage}
\begin{minipage}[c]{0.45\textwidth}
\ifHAL
\hspace{1cm}
\renewcommand{\a}{1.25}
\fi
\resizebox{\a\textwidth}{!}{
\begin{tabular}{|c|c|c|c|}
\hline
\sc{Material}  &  \rule{0pt}{3ex} $\rho^{\sc{f/s}} ~\left[\frac{\rm{kg}}{\rm{m}^3}\right]$  & $c_{\sc{p}}^\sc{f/s} ~\left[ \frac{\rm{m}}{\rm{s}}\right]$  &  $c_{\sc{p}}^\sc{s} ~\left[ \frac{\rm{m}}{\rm{s}}\right]$  \\[1.25ex] \hline
\sc{Atmosphere}  & 1.225  & 343  &  n/a \\ \hline
\sc{Sediments}  & 1300  &  1600  &  900  \\ \hline
\sc{Bedrock}   &  2570  &  5350   &  3009  \\ \hline
\end{tabular}}
\end{minipage}}%
\caption{\textbf{Left panel:} Slightly zoomed-in view of the mesh (with triangles) of the computational domain: a bedrock (such as granite, in turquoise) embeds a sedimentary basin (in dark blue); both are overlaid by a homogeneous atmosphere (in yellow). Sensors (white dots): A barometric sensor $\cal{S}^{\sc{f}}$ is positioned in the atmosphere and a seismometer $\cal{S}^{\sc{s}}$ is positioned in the sedimentary basin. The center of the Ricker wavelet is indicated by a purple point. \textbf{Right panel:} Material properties.}
\label{fig:bassin}
\end{figure}

Our goal is to discretize the geometry with a classical meshing software (here, \texttt{gmsh}), without having to extensively tune the meshing process. Three meshes have therefore been generated: the canonical mesh delivered by \textit{gmsh} with only simplices (\hyperref[fig:bassin]{\Cref{fig:bassin}}, left panel); a second one with quadrangles only; and the last one, a hybrid mesh with both cell geometries (\hyperref[fig:meshes_zoom]{\Cref{fig:meshes_zoom}}, left panel). This latter hybrid mesh illustrates a common situation where, due to geometric complexity, it becomes difficult for the meshing software to maintain a unique geometrical shape for the mesh cells. As a result, undesired element types—typically simplices—are introduced locally within a mesh that was intended to be uniformly shaped. Dealing with this issue typically requires manually modifying the mesh to make it purely \rev{quadrilateral}, a process that can be cumbersome and that frequently results in elements with poor Jacobian conditioning. Moreover, from a coarse mesh, we also generate a nonconforming mesh with hanging nodes, following a one-step refinement in the atmosphere and in the basin (\hyperref[fig:meshes_zoom]{\Cref{fig:meshes_zoom}}, right panel). If both mixed and nonconforming meshes have already been used in other contexts \cite{HKC_2011}, we underline that these can be naturally dealt with using HHO by leveraging its polytopal flexibility, without any specific extra procedure.

All the conforming meshes are generated in order for the acoustic waves and the associated interface waves to be sufficiently sampled in space, so as to suffer from minimal grid dispersion. This corresponds to a length scale of around $14~[\rm{m}]$ along the Earth-atmosphere interface. The nonconforming mesh is generated starting from a characteristic size for the edge length of $28~[\rm{m}]$, and then refining in the atmosphere and the sedimentary basin only. The number of cells is $185\,923$ for the triangular mesh, $94\,019$ for the \rev{quadrilateral} mesh, $93\,490$ for the hybrid mesh, and $84\,204$ for the nonconforming mesh.
\begin{figure}[!htb]
\centering
\resizebox{\textwidth}{!}{
\begin{tikzpicture}
\node[anchor=south west,inner sep=0] (image) at (0,0) {
\includegraphics[width=0.495\linewidth]{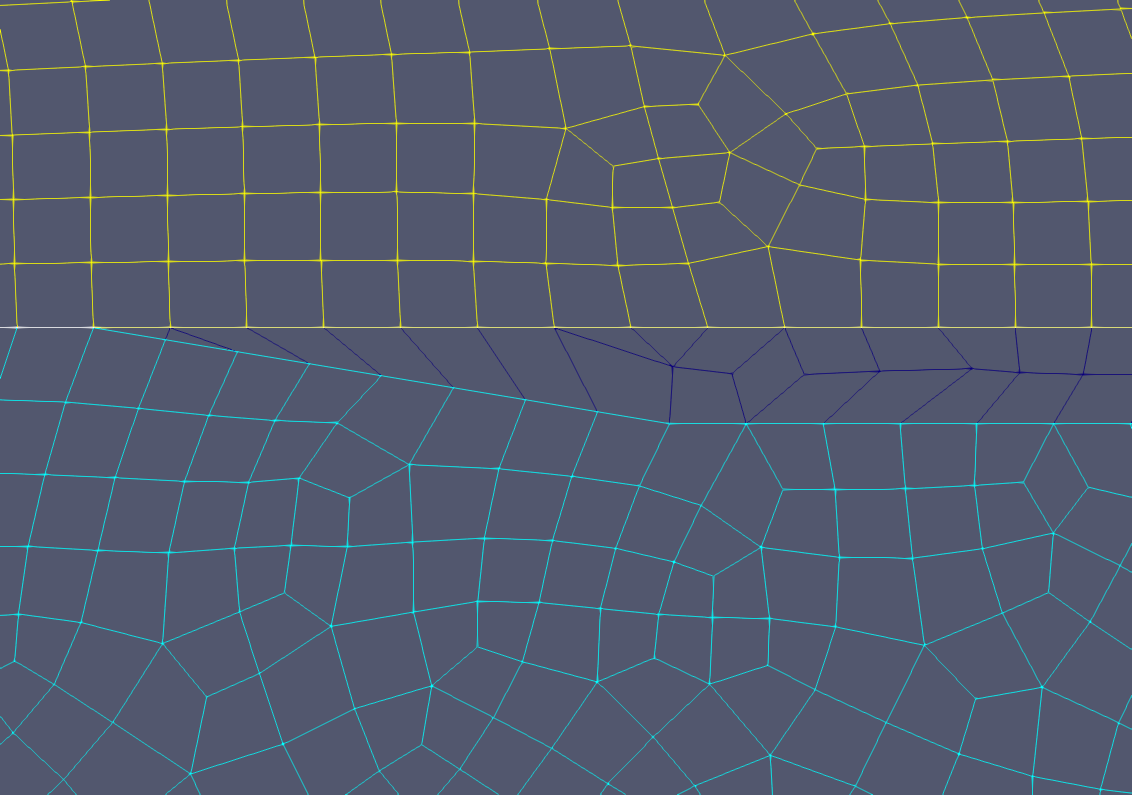}
};
\begin{scope}[x={(image.south east)},y={(image.north west)}]
\node[black] at (0.185,0.9) {\large \textbf{Atmosphere}};
\node[black] at (0.135,0.075) {\large \textbf{Bedrock}};
\node[black] at (0.825,0.525) {\large \textbf{Sediments}};
\end{scope}
\end{tikzpicture}
\hfill \hspace{1cm}
\begin{tikzpicture}
\node[anchor=south west,inner sep=0] (image) at (0,0) {
\includegraphics[width=0.495\linewidth]{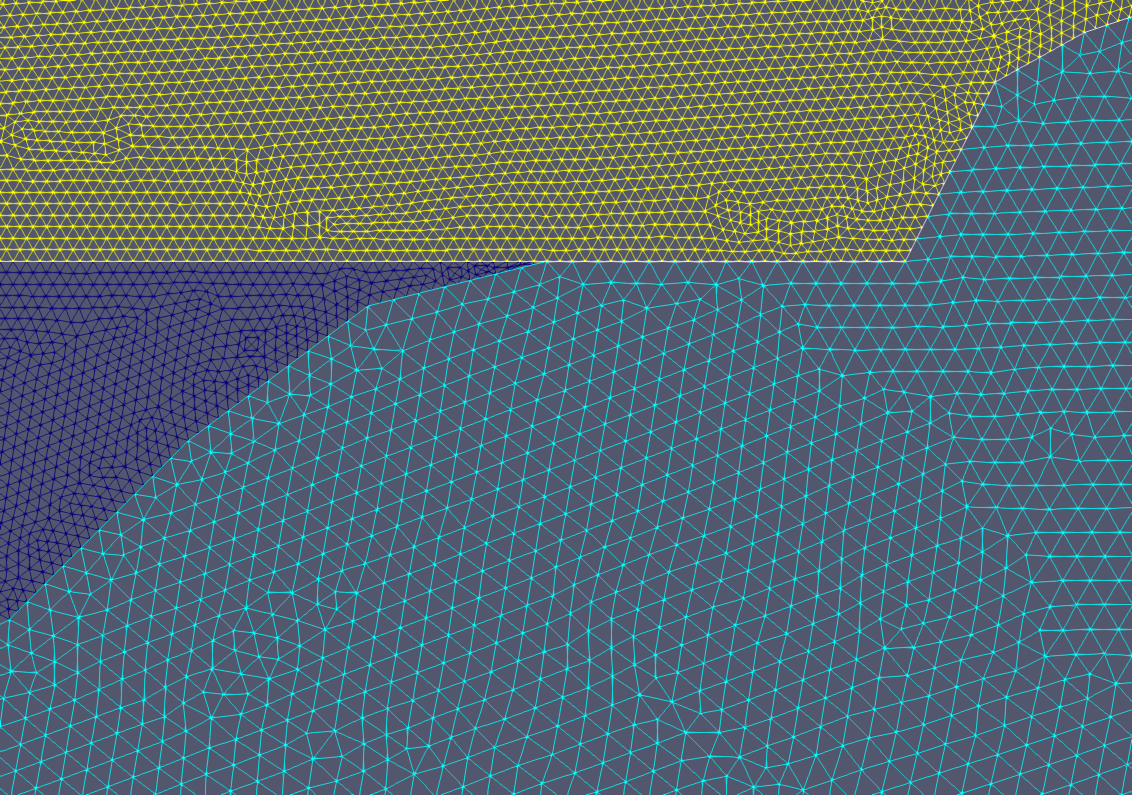}
};
\begin{scope}[x={(image.south east)},y={(image.north west)}]
\node[black] at (0.185,0.9) {\large \textbf{Atmosphere}};
\node[black] at (0.135,0.075) {\large \textbf{Bedrock}};
\node[black] at (0.155,0.625) {\large \textbf{Sediments}};
\end{scope}
\end{tikzpicture}}
\caption{\textbf{Left panel:} Hybrid mesh generated using \texttt{gmsh}: Most of the cells are quadrangles (more than 99$\%$), but with remaining triangles in some parts of the domain. \textbf{Right panel:} Nonconforming mesh with hanging nodes, after a one-pass refinement in the atmosphere and the basin.}
\label{fig:meshes_zoom}
\end{figure}
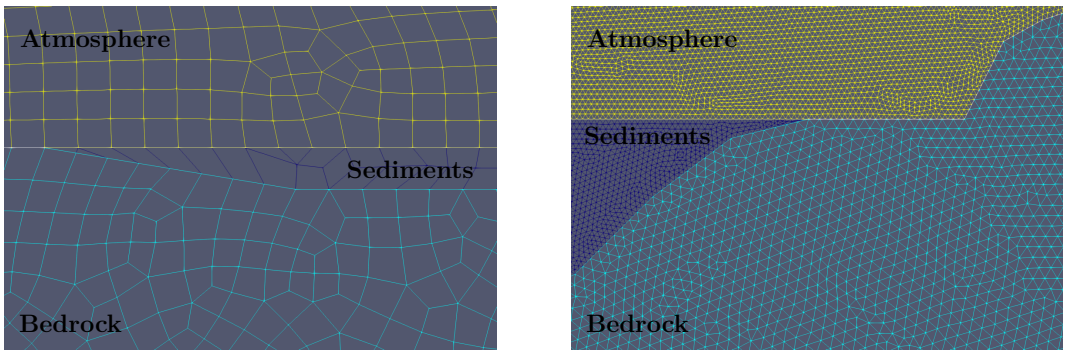

\subsubsection{Numerical results}

The results are obtained using SDIRK$(3,4)$ with polynomial degree $k=3$ and time step $\Delta t = 4.7 \times 10^{-3}~[\rm{s}]$. A reference solution is computed using again the 2D spectral element solver of the SEM3D package. We use here a polynomial degree $k=4$ \rev{or $k=8$} (i.e., $5$ \rev{or $9$} GLL nodes in each direction, respectively) and a time step $\Delta t = 0.61 \times 10^{-4}~[\rm{s}]$. In order to underline the significant reduction in dofs associated with static condensation in the case of an implicit time-stepping, these are reported in \rev{in the first and second lines of} \hyperref[nb_ddls]{\Cref{nb_ddls}} for the different meshes considered before and after static condensation of the cell dofs. The ratio of dofs reduction is about 75\% on all the meshes, a value that is consistent with the estimates in
\hyperref[tab::DOFS]{\Cref{tab::DOFS}}. 
\begin{table}[!h]
\begin{center}
\resizebox{\textwidth}{!}{
\begin{tabular}{c|c|c|c|c|}
\cline{2-5}
& \multicolumn{4}{|c|}{Type of mesh} \\ 
\cline{2-5}
& Triangular  &  \rev{Quadrilateral}  & Hybrid  &  Nonconforming \\ 
\hline
\multicolumn{1}{|l|}{Dofs before static condensation} & $8\,538\,346$ & $4\,562\,532$ &$4\,533\,700$ & $3\,761\,468$  \\ 
\multicolumn{1}{|l|}{Dofs after static condensation} & $1\,448\,120$ & $977\,560$ &$971\,416$ & $810\,344$  \\ \hline
\multicolumn{1}{|l|}{\rev{Relative error at solid sensor}} & 3.07e-3 & 5.84e-3 & 3.91e-3 & 3.59e-3 \\
\multicolumn{1}{|l|}{\rev{Relative error at fluid sensor}} & 4.88e-3 & 7.57e-3 & 8.11e-3 & 7.65e-3 \\ \hline
\end{tabular}
}
\caption{\rev{First and second lines:} Number of dofs for different meshes, before and after static condensation of the cell dofs. \rev{Third and fourth lines: relative errors at the solid and fluid sensors in the $\ell^2(J)$-norm with respect to the reference solution using 9 GLL nodes in each direction.}}
\label{nb_ddls}
\end{center}
\end{table}

\hyperref[fig:snapshot]{\Cref{fig:snapshot}} displays the two-dimensional pressure field in the fluid subdomain and the Euclidian norm of the velocity field in the solid subdomain at the times $t \in \{0, 0.16, 0.6, 1\}~[\rm{s}]$. As expected, seismoacoustics coupling happens where the mountain is located (snapshot at $t = 0.16~[\rm{s}]$): The pressure wave mainly results from the direct transmission of the solid velocity pulse above the mountain, but is also related to inhomogeneous waves that are guided along the interface (lateral head waves). The seismic energy also penetrates the right part of the sedimentary basin. This part of the energy is then trapped in the sedimentary basin (because of differences in impedances with the bedrock), goes back and forth and continuously couples with the atmosphere (snapshots at times $t \in \{0.6, 1\}~[\rm{s}]$). Long after its initial excitation, the basin still leaks energy into the atmosphere, as reported in \cite{AASSE_2018}.
\ifHAL
\begin{figure}[!htb]
\centering 
\includegraphics[width=0.49\linewidth]{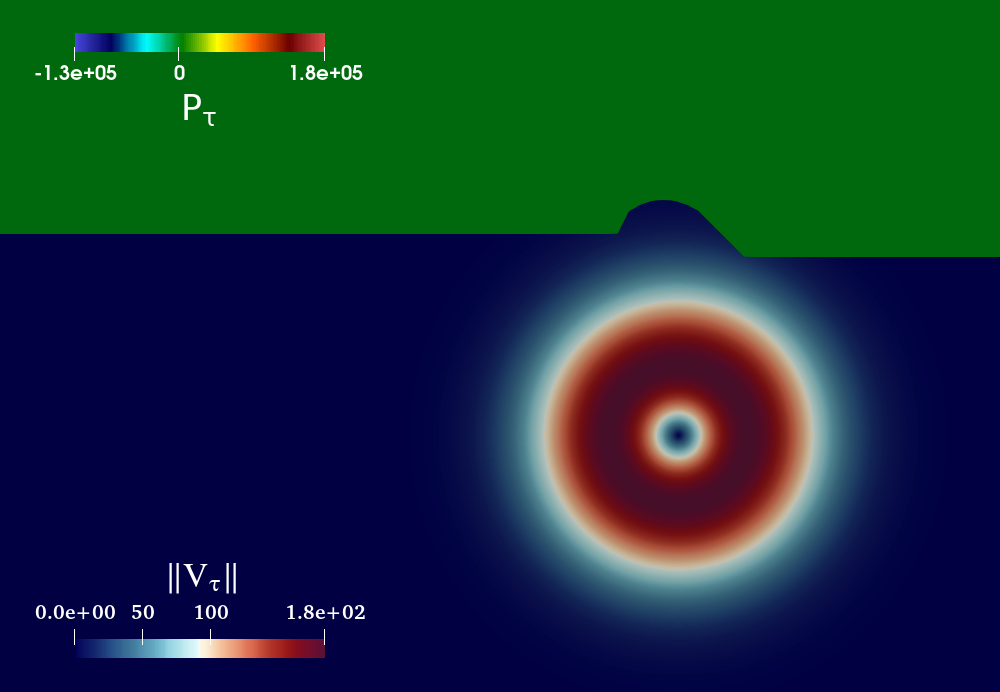}
\includegraphics[width=0.49\linewidth]{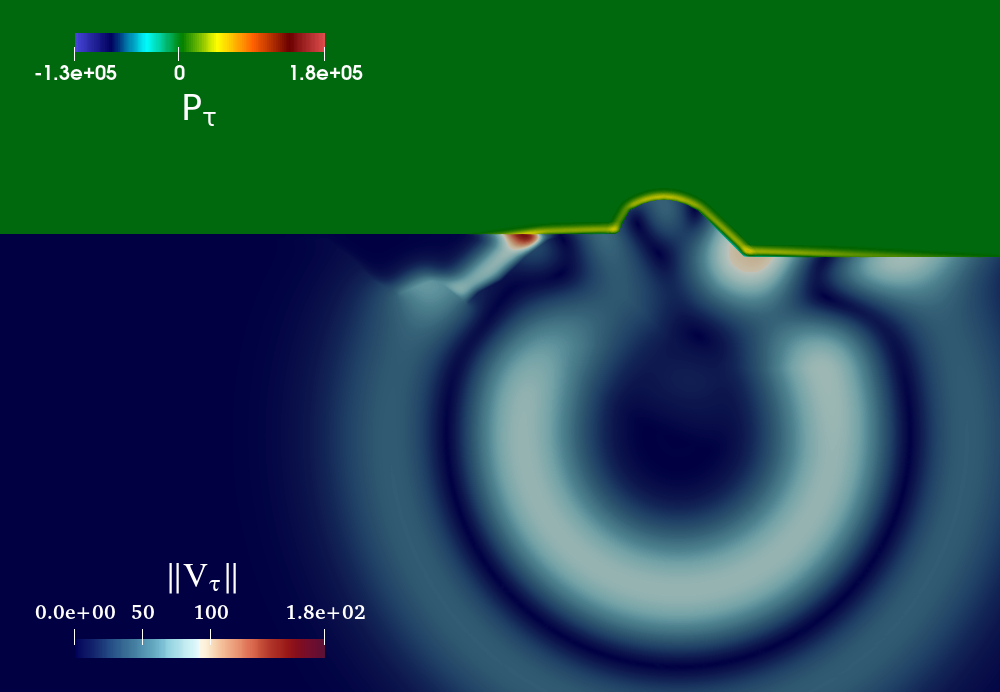}
\includegraphics[width=0.49\linewidth]{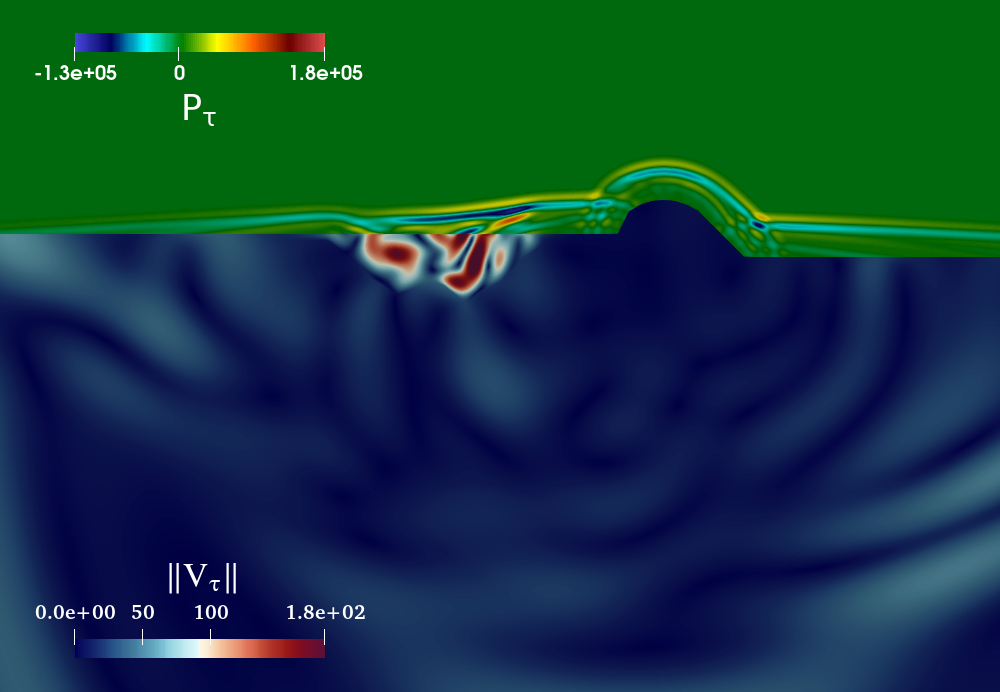}
\includegraphics[width=0.49\linewidth]{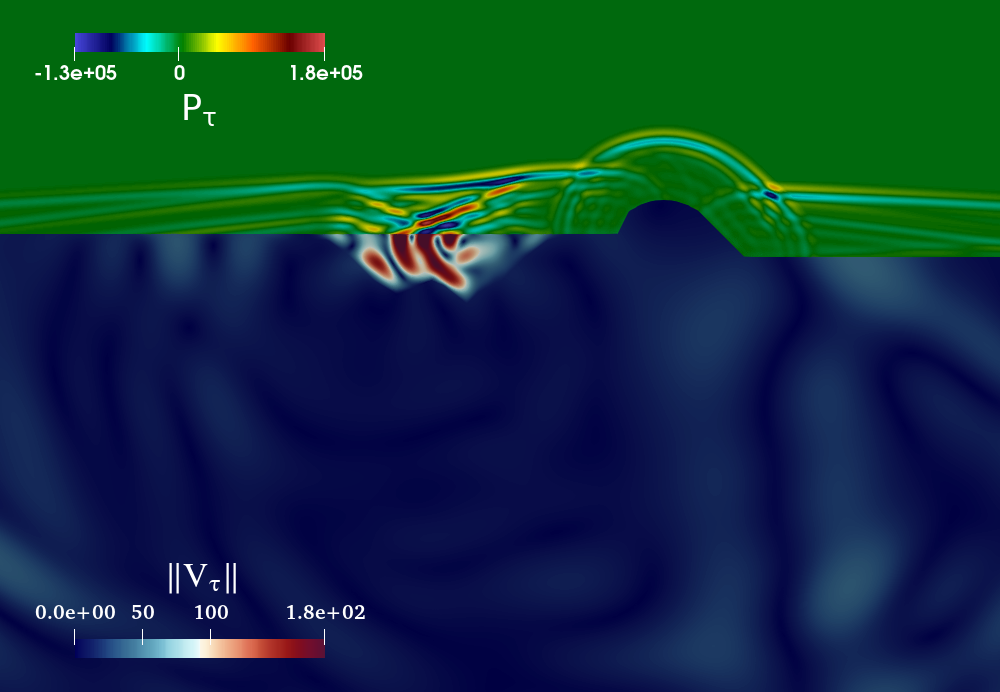} 
\caption{Spatial distribution of the acoustic pressure (atmosphere) and the elastic velocity norm (solid subdomain) at $t \in \{0, 0.16, 0.6, 1\}~[\rm{s}]$. A part of the seismic energy leaks into the atmosphere where the mountain is located; another part is trapped in the sedimentary basin, which continuously emits waves to the atmosphere.}
\label{fig:snapshot}
\end{figure}
\fi

In \hyperref[fig:trace]{\Cref{fig:trace}}, we compare the solution obtained using the HHO discretization with the reference solution \rev{using 5 GLL nodes} at the two sensors located at $\cal{S}^{\sc{f}} := (-0.65, 1.3) ~[\rm{km}]$ and $\cal{S}^{\sc{s}} := (-0.75, 1.08)~[\rm{km}]$. An important conclusion is that the numerical results obtained using the HHO discretization applied to the different meshes are undistinguishable, thereby confirming the robustness of the proposed method with respect to the mesh. Furthermore, the numerical predictions perfectly match the reference solution for which the mesh is (and can only be) \rev{quadrilateral}. \rev{To appreciate the errors more quantitatively, we report in the third and fourth lines of \hyperref[nb_ddls]{\Cref{nb_ddls}} the relative $\ell^2(J)$-errors at both sensors with respect to the reference solution using 9 GLL nodes in each direction. For comparison, these errors are, respectively, $4.68\times10^{-4}$ and $9.78\times10^{-4}$ when comparing the results of the spectral element solver using 5 and 9 GLL nodes.}
These results underline the accuracy of the developed HHO method for the simulation of wavefield propagation in complex media, and its strong potential for extension to large-scale simulations. 
\ifHAL
\begin{figure}[!htb]
\centering 
\includegraphics[width=0.49\linewidth]{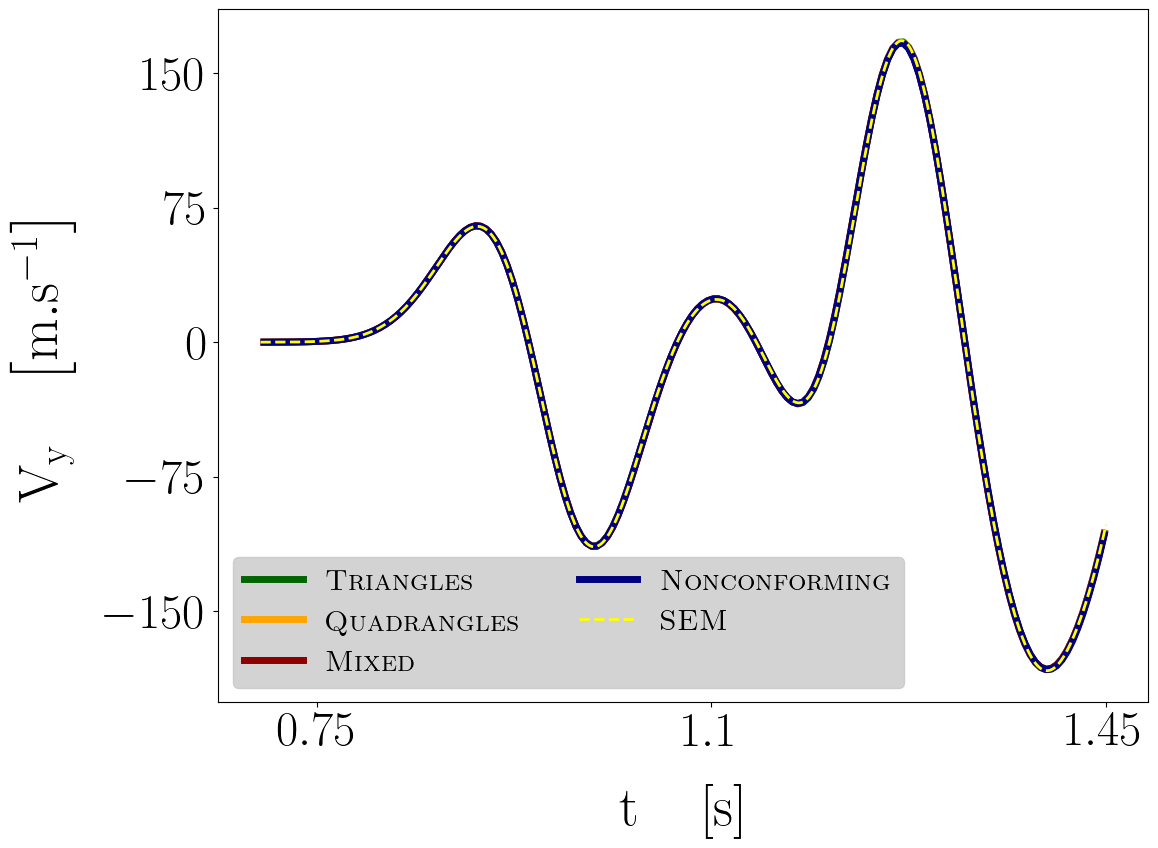}
\includegraphics[width=0.49\linewidth]{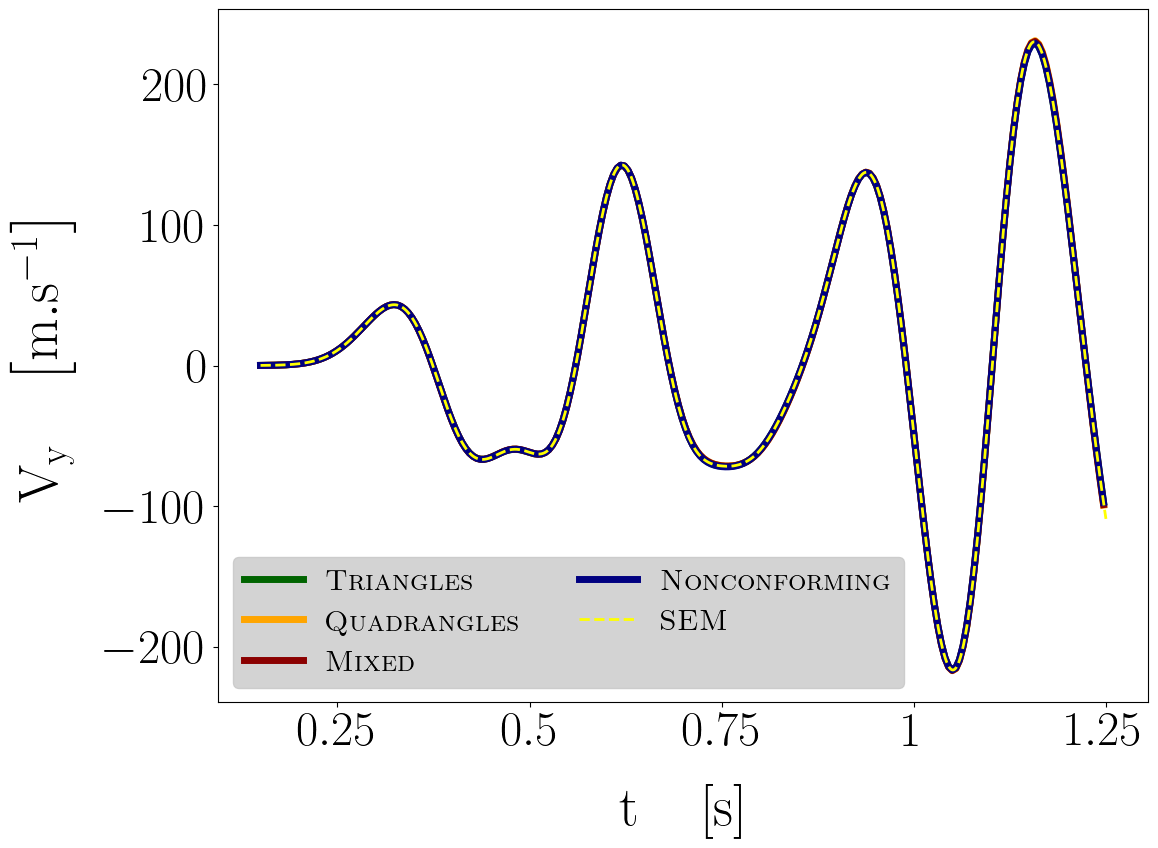}  
\caption{Evolution of the vertical component of the velocity as a fonction of time. \textbf{Left:} Record at the $S^\sc{f}$ sensor, in the atmosphere. \textbf{Right:} Record at the $S^\sc{s}$ sensor, in the sedimentary basin. }
\label{fig:trace}
\end{figure}
\fi
\FloatBarrier

\section{Conclusion} 

We proposed a HHO method for coupled elasto-acoustic wave propagation, combining first-order formulations in time with explicit and implicit Runge--Kutta schemes. Static condensation significantly reduces computational costs  in both cases. A CFL stability condition was derived numerically, and comparisons showed good performances, with implicit schemes favored. Numerical experiments confirmed the accuracy and robustness of the proposed methodology on hybrid and nonconforming meshes, including realistic seismic scenarios.
 
A natural extension of this work consists in coupling the HHO discretization with SEM, leveraging the geometric flexibility of HHO in complex regions while exploiting the efficiency of the SEM in the bulk where structured meshes can be used. This combined approach would enable realistic, large-scale 3D geophysical simulations to be performed. Transitioning to three dimensions will also require significant numerical and computational developments, including parallelization and possibly GPU acceleration to align with emerging exascale architectures. Such efforts will likely necessitate a renewed performance analysis, in the spirit of \cite{KSMW_2016}. In parallel, further improvements such as structure-preserving time integrators (see \cite{SCNPC_2017, SCNP_2021} for HDG space semi-discretizations) and superconvergent post-processing tailored to HHO can be explored to enhance long-time accuracy. The use of alternative polynomial bases will also be investigated as a possible means to improve accuracy and computational efficiency. 
Finally, an HHO-based approach combined with adaptive mesh refinement would be particularly well-suited for wave propagation problems. It could leverage the polyhedral flexibility of the method to allow for nonconforming, locally refined meshes, thereby enhancing accuracy in regions with high solution gradients while controlling computational costs.

\FloatBarrier 
\bibliographystyle{abbrv}
\bibliography{references}
\ifHAL
\renewcommand{\bibname}{References} 
\addcontentsline{toc}{section}{References}
\fi

\end{document}

%% file: preamble.tex
\ifHAL

\usepackage[utf8]{inputenc}           
\usepackage[T1]{fontenc}             
\usepackage[left=2.25cm, right=2.25cm, top=2cm, bottom=2.5cm]{geometry} 
\usepackage[english]{babel}  
\usepackage{lmodern}
\usepackage{authblk}
\usepackage[subfigure]{tocloft}      
\usepackage{titletoc}
\usepackage{textcomp}                
\usepackage{gensymb}                 
\usepackage{enumerate}               
\usepackage{enumitem}                
\usepackage{titlesec}                
\usepackage{color}                   
\usepackage{indentfirst}             
\usepackage{dsfont} 
\usepackage{fancyhdr}                
\usepackage{chngpage}                
\usepackage[symbol]{footmisc}        
\usepackage{amsmath,amsfonts,amssymb} 
\usepackage{stmaryrd}
\usepackage{mathtools}                
\usepackage{amsthm}                   
\usepackage{empheq,etoolbox}          
\usepackage{mathbbol}
\usepackage{mathrsfs}
\usepackage{yfonts}
\patchcmd{\subequations}{\theparentequation\alph{equation}}
{\theparentequation\alph{equation}}{}{}
\newtheorem{theorem}{Theorem}[section]

\newtheorem{lemma}[theorem]{Lemma}

\usepackage{nameref} 
\usepackage[hidelinks]{hyperref} 
\usepackage{pdfpages}
\usepackage{scalerel}
\titleformat*{\section}{\LARGE\bfseries}
\titleformat*{\subsection}{\Large\bfseries}
\titleformat*{\subsubsection}{\large\bfseries}
\titlespacing*{\section}{0pt}{0.25cm}{1em}
\titlespacing*{\subsection}{0pt}{0.75em}{1em}

\fi
\usepackage{multicol}
\usepackage{float}
\usepackage{algorithm}
\usepackage{algorithmicx}
\usepackage{algpseudocode}
\usepackage[dvipsnames]{xcolor}     
\definecolor{ceared}{HTML}{BB0000}   
\definecolor{color_dofs}{HTML}{A3167C}
\definecolor{reddishgray}{rgb}{0.6, 0.45, 0.45}
\usepackage{placeins}
\usepackage[font=small, skip=10pt]{caption} 
\usepackage{subcaption}                     
\usepackage{tikz}                           
\usepackage{tkz-euclide}                    
\usetikzlibrary{patterns, patterns.meta}    
\usepackage{graphicx}
\usepackage{pgfplots}
\usetikzlibrary{shapes, positioning}
\pgfplotsset{compat=1.16}
\usepackage{multirow}
\usepackage{pmat}                     
\usepackage{wasysym}                  
\usepackage{bbm}
\usepackage{cleveref}

\ifHAL
\fi

%% file: macros.tex
\renewcommand{\sc}{\textsc}
\renewcommand{\cal}[1]{\mathcal{#1}}
\renewcommand{\rm}{\mathrm}
\newcommand{\bb}[1]{\mathbb{#1}}
\newcommand{\bbm}{\mathbbm}
\newcommand{\bd}[1]{\boldsymbol{#1}}
\newcommand{\T}{\cal{T}}
\newcommand{\F}{\cal{F}}

\newcommand{\dofs}[3]{\rm{#1}_{\cal{#2}^\sc{#3}}}    
\newcommand{\mass}[2]{\cal{M}_{\T^\sc{#2}}^{#1}}  
\newcommand{\stab}[3]{{\Sigma}_{\cal{#2}^{\sc{#1}}\cal{#3}^{\sc{#1}}}} 
\newcommand{\stabdiag}[2]{{\Sigma}_{\cal{#2}^{\sc{#1}}}} 
\newcommand{\grad}[3]{\cal{G}_{\cal{#1}^{\sc{#3}}\cal{#2}^{\sc{#3}}}} 
\newcommand{\graddiag}[2]{\cal{G}_{\cal{#1}^{\sc{#2}}}} 
\newcommand{\straindiag}[1]{\cal{H}_{\cal{#1}^{\sc{s}}}} 
\newcommand{\strain}[2]{\cal{H}_{\cal{#1}^{\sc{s}}\cal{#2}^{\sc{s}}}} 
\newcommand{\coupling}{\cal{C}_{\F^{\Gamma}}} 

\newcommand{\G}{\Gamma}                   
\newcommand{\domain}[1]{\Omega^{\sc{#1}}} 

\newcommand{\Ts}{\cal{T}^\sc{s}}
\newcommand{\Tf}{\cal{T}^\sc{f}}

\newcommand{\sF}{\sc{f}}
\newcommand{\sS}{\sc{s}}

\newcommand{\cF}{\mathcal{F}}

\newcommand{\tih}{\tilde{h}}

\newcommand{\Ff}{\cal{F}^{\sF}}
\newcommand{\Fs}{\cal{F}^{\sS}}

\newcommand{\rev}[1]{{\color{black}{#1}}}

\newcommand{\normalizecell}[1]{\strut#1\strut}
\numberwithin{equation}{section}